\documentclass[a4paper, 11pt]{article}

\usepackage{fullpage}

\usepackage[english]{babel}
\usepackage[utf8]{inputenc}
\usepackage{amsmath, amsfonts, amssymb}
\usepackage{mathtools}
\usepackage{stmaryrd}
\usepackage{bm}
\usepackage{authblk}
\usepackage{cite}

\usepackage[colorlinks=true, allcolors=blue]{hyperref}
\usepackage{booktabs}
\usepackage{graphicx}
\usepackage{subcaption}

\usepackage[ruled,vlined]{algorithm2e}
\SetKwInput{KwInput}{input}
\SetAlFnt{\small}

\newcommand{\Real}{\mathbb{R}}

\makeatletter%
\@ifclassloaded{beamer}%
    {}%
    {

        \newtheorem{definition}{Definition}[section]
        
        \newtheorem{remark}[definition]{Remark}
    }%
\makeatother%

\title{Sketching for low-rank nonnegative matrix approximation: a numerical study}
\author[1,2]{Sergey Matveev}
\author[2]{Stanislav Budzinskiy}
\affil[1]{Faculty of Computational Mathematics and Cybernetics, Lomonosov Moscow State University}
\affil[2]{Marchuk Institute of Numerical Mathematics RAS}
\date{}

\begin{document}
\maketitle

\begin{abstract}
We propose new approximate alternating projection methods, based on randomized sketching, for the low-rank nonnegative matrix approximation problem: find a low-rank approximation of a nonnegative matrix that is nonnegative, but whose factors can be arbitrary. We calculate the computational complexities of the proposed methods and evaluate their performance in numerical experiments. The comparison with the known deterministic alternating projection methods shows that the randomized approaches are faster and exhibit similar convergence properties.
\end{abstract}

\tableofcontents

\clearpage
\section{Introduction}
Nonnegative functions and datasets arise in many areas of research and industrial applications. They come in different forms that include, but are not limited to, probability density functions, concentrations of substances in physics and chemistry, ratings, images, and videos. As the amounts of data grow, it becomes increasingly important to have reliable approximation techniques that permit fast post-processing of the compressed data in the low-parametric format. 

For matrices, this can be achieved with low-rank approximations. The best low-rank approximation problem is well-understood and has an exact solution in any unitarily invariant norm: the truncated singular value decomposition. However, there is no guarantee that the resulting low-rank matrix retains the nonnegativity of its elements.

A possible remedy can be found in nonnegative matrix factorizations (NMF) \cite{wang2012nonnegative, gillis2020nonnegative}, an ideology that explicitly enforces nonnegativity by searching for an approximate low-rank decomposition with nonnegative latent factors. A lot of progress has been made in this field; for instance, it has been extended to multi-dimensional tensors \cite{cichocki2009nonnegative, zhou2015efficient, lee2016nonnegative, shcherbakova2019nonnegative_lssc}.

The interpretative properties of NMF explain why it shines in such areas as data analysis \cite{fu2019nonnegative}. In other applications, however, the main goal is to achieve good compression of the data. This mostly concerns scientific computing in such areas as numerical solution of large-scale differential equations \cite{smirnov2016fast, matveev2016tensor, dolgov2021tensor, chertkov2021fokkerplanck, allmann2022parallel} and multivariate probability \cite{dolgov2020approximation, novikov2021tensor}. There, NMF can appear to be a bottleneck, since the nonnegative rank might be significantly larger than the usual matrix rank (not to mention the computational complexity of approximate NMF; the exact NMF is NP-hard \cite{vavasis2010complexity}). 

For these tasks, a recently proposed low-rank nonnegative matrix factorization (LRNMF) problem \cite{SongNgNonnegative2020} is more suitable. Given a matrix $X \in \Real^{m \times n}_{+}$ with nonnegative entries and a target rank $r$, the goal is to find a rank-$r$ approximation that is nonnegative, but whose factors can be arbitrary:
\begin{equation}\label{intro:eq:lrnmf1}
    \| X - U V^T \|_F^2 \to \min_{U \in \Real^{m \times r}, V \in \Real^{n \times r}} \quad \text{s.t.} \quad U V^T \in \Real^{m \times n}_{+}.
\end{equation}
Likely even more important for the applications in question is a closely related problem, where the nonnegative matrix $X$ itself is unknown; instead, one has a rank-$r$ approximation $X \approx U_0 V_0^T$ that contains negative elements, and the goal is to produce a nonnegative rank-$r$ matrix that is close to it:
\begin{equation}\label{intro:eq:lrnmf2}
    \| U_0 V_0^T - U V^T \|_F^2 \to \min_{U \in \Real^{m \times r}, V \in \Real^{n \times r}} \quad \text{s.t.} \quad U V^T \in \Real^{m \times n}_{+}.
\end{equation}

In \cite{SongNgNonnegative2020}, it was shown that \eqref{intro:eq:lrnmf1} and \eqref{intro:eq:lrnmf2} can be solved with alternating projections. Classically, this method is used to find a point in the intersection of two closed convex sets $C_1$ and $C_2$ in a Euclidean space:
\begin{equation}
\label{intro:eq:ap}
    y_k = \Pi_{C_1}(x_k), \quad x_{k+1} = \Pi_{C_2}(y_k),
\end{equation}
where $\Pi_{C}(x) = \arg\min_{y \in C} \| x - y \|$ denotes the projection of $x$ onto convex $C$ in the Euclidean norm. While $\Real^{m \times n}_{+}$ is indeed convex, the set $\mathcal{M}_r \subset \Real^{m \times n}$ of rank-$r$ matrices is not. The best rank-$r$ approximation in the Frobenius norm (which is Euclidean) still exists, though, and it was proved that the iterations converge to a point in $\mathcal{M}_r \cap \Real^{m \times n}_{+}$ \cite{SongNgNonnegative2020}. The following work \cite{SongEtAlTangent2020} introduced an approximate projection operator onto $\mathcal{M}_r$ that reduced the per-iteration complexity of the alternating projections to $O(mnr)$. 

In this paper, we consider the alternating projections for LRNMF \eqref{intro:eq:lrnmf1}-\eqref{intro:eq:lrnmf2} with approximate projections onto $\mathcal{M}_r$ that are based on randomized sketching techniques with the aim to further reduce the complexity of the iterations, as compared with \cite{SongEtAlTangent2020}. In Section~\ref{sec:methods}, we present in detail the algorithms from \cite{SongNgNonnegative2020, SongEtAlTangent2020} and describe three randomized approximate alternating projection approaches based on \cite{HalkoEtAlFinding2011a, TroppEtAlPractical2017, NakatsukasaFast2020}. Section~\ref{sec:exp} is devoted to numerical experiments that we use to evaluate and compare the performance of the different versions of alternating projections. The three examples are random matrices, an image, and a solution to the Smoluchowski equation. Our paper also has an Appendix, where we meticulously calculate the computational complexities of the presented algorithms.

\section{Methods}\label{sec:methods}
\subsection{Deterministic}
The first algorithm we describe directly follows the alternating projection framework \eqref{intro:eq:ap} as it computes the exact projections onto $\mathcal{M}_{\leq r} = \{ X \in \Real^{m \times n} : \mathrm{rank} X \leq r \}$ and $\Real^{m \times n}_{+}$ that minimize the Frobenius norm \cite{SongNgNonnegative2020}. For matrices, the best rank-$r$ approximation in any unitarily invariant norm (including the Frobenius norm) is delivered by the truncated singular value decomposition (SVD). We denote it by $\mathrm{SVD}_r$ so that
\begin{equation*}
    \Pi_{\mathcal{M}_{\leq r}}(X) = \mathrm{SVD}_r(X) = U_r \Sigma_r V_r^T = \arg\min_{Y \in \mathcal{M}_{\leq r}} \| X - Y \|_F,
\end{equation*}
where $U_r \in \Real^{m \times r}$ and $V_r \in \Real^{n \times r}$ are the truncated left and right singular factors, and $\Sigma_r \in \Real^{r \times r}$ is the truncated diagonal matrix of singular values. To find the nonnegative matrix that best approximates a given one, we simply need to set all of its negative elements to zero:
\begin{equation*}
    \Pi_{\Real^{m \times n}_{+}}(X) = \max(X, 0) = \{ \max(x_{ij}, 0) \} = \arg\min_{Y \in \Real^{m \times n}_{+}} \| X - Y \|_F.
\end{equation*}
Applied iteratively one by one, these projections give rise to Alg.~\ref{alg:ap} that we will call \textit{SVD}. The computational complexity of this approach is dominated by $\mathrm{SVD}_r$, which costs $O(mn^2)$ flops.
\begin{algorithm}[th]
\caption{Exact alternating projections (SVD) \cite{SongNgNonnegative2020}}
\label{alg:ap}
\KwData{Initial approximation $Y^{(0)} \in \Real^{m \times n}$ of rank $r$, number of iterations $s$}
\For{$i = 1, \ldots, s$}{
    $X^{(i)} \gets \max(Y^{(i-1)}, 0)$\;
    $[U_r, \Sigma_r, V_r] \gets \textsc{SVD}_r(X^{(i)})$\;
    $Y^{(i)} \gets U_r \Sigma_r V_r^T$\;
}
\Return{$Y^{(s)}$}
\end{algorithm}

The algorithm from \cite{SongEtAlTangent2020} cleverly uses the smooth manifold structure of $\mathcal{M}_{r}$ to trade the exact low-rank projection for faster iterations. If $Y \in \mathcal{M}_{r}$ and $U_r \in \Real^{m \times r}$ and $V_r \in \Real^{n \times r}$ are its singular factors, the tangent space to $\mathcal{M}_{r}$ at $Y$ can be described as
\begin{equation*}
    T_{Y}\mathcal{M}_{r} = \left\{ U_r A^T + B V_r^T~:~ A \in \Real^{n \times r}, \quad B \in \Real^{m \times r} \right\},
\end{equation*}
and the orthogonal projection onto it is computed according to
\begin{equation*}
    \Pi_{T_{Y}\mathcal{M}_{r}}(X) = U_r U_r^T X + (I - U_r U_r^T) X V_r V_r^T.
\end{equation*}
It is easy to see that all the matrices from the tangent space $T_{Y}\mathcal{M}_{r}$ have rank at most $2r$. This motivates the following approximate projection: let $Y \in \mathcal{M}_{r}$ be the last iterate and let $X = \Pi_{\Real^{m \times n}_{+}}(Y)$ be its nonnegative correction, then the new low-rank iterate is chosen as
\begin{equation*}
    \tilde{\Pi}_{\mathcal{M}_{\leq r}}(X) = \mathrm{SVD}_r(\Pi_{T_{Y}\mathcal{M}_{r}}(X)).
\end{equation*}
The resulting Alg.~\ref{alg:tap} (\textit{Tangent}) has asymptotic complexity $O(mnr)$ per iteration with the dominant term $6mnr$. Note that it would be more correct to talk about the real-algebraic variety $\mathcal{M}_{\leq r}$ and its tangent cones since it may occur that the rank of the iterate is below $r$, but such singular cases are, nevertheless, extremely rare (see \cite{schneider2015convergence}).
\begin{algorithm}
\caption{Alternating projections via tangent spaces (Tangent) \cite{SongEtAlTangent2020}}
\label{alg:tap}
\KwData{Initial approximation $Y^{(0)} \in \Real^{m \times n}$ of rank $r$, number of iterations $s$}
$[U_r, \Sigma_r, V_r] \gets \textsc{SVD}_r(Y^{(0)})$\;
\For{$i = 1, \ldots, s$}{
    $X^{(i)} \gets \max(Y^{(i-1)}, 0)$\;
    $G_1 \gets U_r^T X^{(i)} (I - V_r V_r^T) \in \Real^{r \times n}, \quad G_2 \gets (I - U_r U_r^T) X^{(i)} V_r \in \Real^{m \times r}$\;
    $[Q_1, R_1] \gets \textsc{QR}(G_1^T), \quad [Q_2, R_2] \gets \textsc{QR}(G_2)$\;
    $Z \gets \begin{bmatrix} U_r^T X^{(i)} V_r & R_1^T \\ R_2 & 0 \end{bmatrix} \in \Real^{2r \times 2r}$\;
    $[\tilde{U}_r, \Sigma_r, \tilde{V}_r] \gets \textsc{SVD}_r(Z)$\;
    $U_r \gets \begin{bmatrix} U_r & Q_2 \end{bmatrix} \tilde{U}_r, \quad V_r \gets \begin{bmatrix} V_r & Q_1 \end{bmatrix} \tilde{V}_r$\;
    $Y^{(i)} \gets U_r \Sigma_r V_r^T$\;
}
\Return{$Y^{(s)}$}
\end{algorithm}

\subsection{Randomized}
The family of sketching techniques for low-rank matrix approximation is an important part of rapidly developing randomized numerical linear algebra \cite{martinsson2020randomized}. The general idea consists in dimension reduction, which is achieved with the multiplication by a random test matrix that is sampled from a probability distribution of choice; the reduced matrix, the sketch, is then used to estimate the range and co-range of the initial matrix. We will focus on three classes of random test matrices $\Psi \in \Real^{n \times k}$ \cite{TroppEtAlPractical2017}:
\begin{itemize}
    \item iid standard Gaussian entries
    \begin{equation*}
        \psi_{ij} \sim \mathcal{N}(0,1);
    \end{equation*}
    \item iid Rademacher entries
    \begin{equation*}
        \psi_{ij} \sim \mathrm{Rad}, \quad \psi_{ij} = \begin{cases}
            1, & \text{ with probability } 1/2,\\
            -1, & \text{ with probability } 1/2;
        \end{cases}
    \end{equation*}
    \item iid Rademahcer entries on a sparse mask with density $\rho$
    \begin{equation*}
        \psi_{ij} \sim \mathrm{Rad}(\rho), \quad \psi_{ij} = \begin{cases}
            0, & \text{ with probability } 1-\rho,\\
            1, & \text{ with probability } \rho/2,\\
            -1, & \text{ with probability } \rho/2.
        \end{cases}
    \end{equation*}
\end{itemize}

We begin with the randomized truncated SVD algorithm from \cite[Alg.~5.1]{HalkoEtAlFinding2011a}. It applies a random test matrix $\Psi \in \Real^{n \times k}$ with $k \geq r$ on the right, computes the range of the sketch via QR decomposition, projects the initial matrix onto this range, and finally computes $\mathrm{SVD}_r$ of a fat $k \times n$ matrix:
\begin{equation*}
    \tilde{\Pi}_{\mathcal{M}_{\leq r}}(X) = Q \cdot \mathrm{SVD}_r(Q^T X), \quad [Q, R] = \mathrm{QR}(X \Psi).
\end{equation*}
It is also possible to combine this procedure with $p$ iterations of the power method that lead to better estimation of the range:
\begin{equation*}
    [Q, R] = \mathrm{QR}((X X^T)^p X \Psi).
\end{equation*}
In Alg.~\ref{alg:ap_hmt}, we present the alternating projection method \textit{HMT}$(p, k)$, which uses an equivalent form of the power method \cite[Alg.~4.4]{HalkoEtAlFinding2011a}. 
\begin{algorithm}[th]
\caption{Alternating projections via \cite{HalkoEtAlFinding2011a} (HMT)}
\label{alg:ap_hmt}
\KwData{Initial approximation $Y^{(0)} \in \Real^{m \times n}$ of rank $r$, co-range sketch size $k \geq r$, number of power method iterations $p$, test matrix generator $\textsc{TestMatrix}$, number of iterations $s$}
\For{$i = 1, \ldots, s$}{
    $X^{(i)} \gets \max(Y^{(i-1)}, 0)$\;
    $\Psi \gets \textsc{TestMatrix}(n, k) \in \Real^{n \times k}$\;
    $Z_1 \gets X^{(i)} \Psi \in \Real^{m \times k}$\;
    $[Q, R] \gets \textsc{QR}(Z_1)$\;
    \For{$j = 1, \ldots, p$}{
        $Z_2 \gets Q^T X^{(i)} \in \Real^{k \times n}$\;
        $[Q, R] \gets \textsc{QR}(Z_2^T)$\;
        $Z_1 \gets X^{(i)} Q \in \Real^{m \times k}$\;
        $[Q, R] \gets \textsc{QR}(Z_1)$\;
    }
    $Z_2 \gets Q^T X^{(i)} \in \Real^{k \times n}$\;
    $[U_r, \Sigma_r, V_r] \gets \textsc{SVD}_r(Z_2)$\;
    $Y^{(i)} \gets Q U_r \Sigma_r V_r^T$\;
}
\Return{$Y^{(s)}$}
\end{algorithm}

For the next variant of sketching-based alternating projections, we use a different randomized SVD algorithm \cite{TroppEtAlPractical2017}. Unlike the previous method, it uses two test matrices $\Phi \in \Real^{l \times m}$ and $\Psi \in \Real^{n \times k}$:
\begin{equation*}
    \tilde{\Pi}_{\mathcal{M}_{\leq r}}(X) = Q \cdot \mathrm{SVD}_r( (\Phi Q)^{\dagger} \Phi X), \quad [Q, R] = \mathrm{QR}(X \Psi).
\end{equation*}
The details of the corresponding \textit{Tropp}$(k, l)$ alternating projections are listed in Alg.~\ref{alg:ap_tropp}.
\begin{algorithm}[th]
\caption{Alternating projections via \cite{TroppEtAlPractical2017} (Tropp)}
\label{alg:ap_tropp}
\KwData{Initial approximation $Y^{(0)} \in \Real^{m \times n}$ of rank $r$, co-range sketch size $k \geq r$, range sketch size $l \geq k$, test matrix generator $\textsc{TestMatrix}$, number of iterations $s$}
\For{$i = 1, \ldots, s$}{
    $X^{(i)} \gets \max(Y^{(i-1)}, 0)$\;
    $\Psi \gets \textsc{TestMatrix}(n, k) \in \Real^{n \times k}, \quad \Phi \gets \textsc{TestMatrix}(l, m) \in \Real^{l \times m}$\;
    $Z \gets X^{(i)} \Psi \in \Real^{m \times k}$\;
    $[Q, R] \gets \textsc{QR}(Z)$\;
    $W \gets \Phi Q \in \Real^{l \times k}$\;
    $[P, T] \gets \textsc{QR}(W)$\;
    $G \gets T^{-1} P^T \Phi X^{(i)} \in \Real^{k \times n}$\;
    $[U_r, \Sigma_r, V_r] \gets \textsc{SVD}_r(G)$\;
    $Y^{(i)} \gets Q U_r \Sigma_r V_r^T$\;
}
\Return{$Y^{(s)}$}
\end{algorithm}

Finally, we consider the generalized Nystr\"om method \cite{NakatsukasaFast2020} that does not use the SVD whatsoever and is the base of \textit{GN}$(l)$ (see Alg.~\ref{alg:ap_gn}). Given two test matrices $\Phi \in \Real^{l \times m}$ and $\Psi \in \Real^{n \times r}$, it computes
\begin{equation*}
    \tilde{\Pi}_{\mathcal{M}_{\leq r}}(X) = (X \Psi R^{-1}) (Q^T \Phi X), \quad [Q, R] = \mathrm{QR}(\Phi X \Psi).
\end{equation*}
\begin{algorithm}
\caption{Alternating projections via \cite{NakatsukasaFast2020} (GN)}
\label{alg:ap_gn}
\KwData{Initial approximation $Y^{(0)} \in \Real^{m \times n}$ of rank $r$, range sketch size $l \geq r$, test matrix generator $\textsc{TestMatrix}$, number of iterations $s$}
\For{$i = 1, \ldots, s$}{
    $X^{(i)} \gets \max(Y^{(i-1)}, 0)$\;
    $\Psi \gets \textsc{TestMatrix}(n, r) \in \Real^{n \times r}, \quad \Phi \gets \textsc{TestMatrix}(l, m) \in \Real^{l \times m}$\;
    $Z \gets X^{(i)} \Psi \in \Real^{m \times r}$\;
    $W \gets \Phi Z \in \Real^{l \times r}$\;
    $[Q, R] \gets \textsc{QR}(W)$\;
    $V \gets (\Phi X^{(i)})^T Q \in \Real^{n \times r}$\;
    $U \gets Z R^{-1} \in \Real^{m \times r}$\;
    $Y^{(i)} \gets U V^T$\;
}
\Return{$Y^{(s)}$}
\end{algorithm}

In all three sketching-based alternating projection methods, the computational complexity of a single iteration is determined by matrix-matrix products, similarly to the Tangent approach. However, the constant in front of $mn$ can be reduced for the randomized algorithms (see Tab.~\ref{methods:tab:complexity}). Indeed, while with Gaussian sketching the lowest value that can be obtained is $6r$ (if we set $p = 0$, $k = r$, and $l = r$) as in Tangent, Rademacher and sparse Rademacher sketching can lead to smaller complexities. The detailed computational complexity analysis is carried out in the Appendix.
\begin{table}[h]\centering
\begin{tabular}{@{}lccc@{}}
    \toprule
     & HMT$(p, k)$ & Tropp$(k, l)$ & GN$(l)$ \\
    \midrule
    $\mathcal{N}(0, 1)$ & $(4p+4)k + 2r$ & $2r+2k+2l$ & $4r+2l$ \\
    $\mathrm{Rad}$ & $(4p+3)k + 2r$ & $2r+k+l$ & $3r+l$\\
    $\mathrm{Rad}(\rho)$ & $(4p+2+\rho)k + 2r$ & $2r+\rho k + \rho l$ & $2r + \rho r + \rho l$ \\
    \bottomrule\\
\end{tabular}
\caption{The dominant terms, divided by $mn$, of per-iteration complexities related to matrix-matrix products. For the Tangent method, this value is $6r$.}
\label{methods:tab:complexity}
\end{table} 
\section{Numerical experiments}\label{sec:exp}
\subsection{Random uniform matrices}
\begin{figure}[th]
\begin{subfigure}[b]{0.45\textwidth}
\centering
	\includegraphics[width=0.8\textwidth]{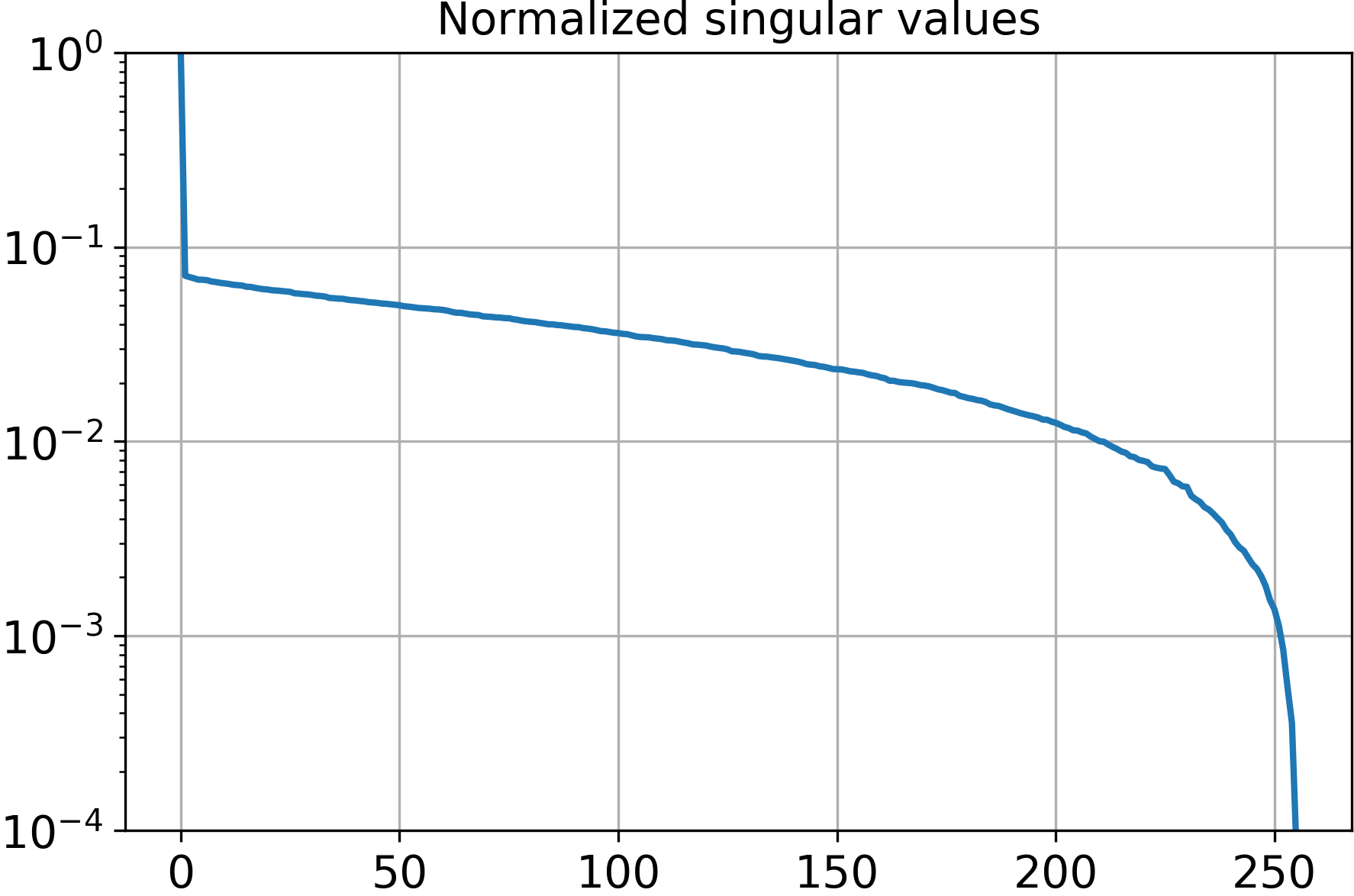}
	\caption{}
\end{subfigure}\hfill%
\begin{subfigure}[b]{0.45\textwidth}
\centering
	\includegraphics[width=0.8\textwidth]{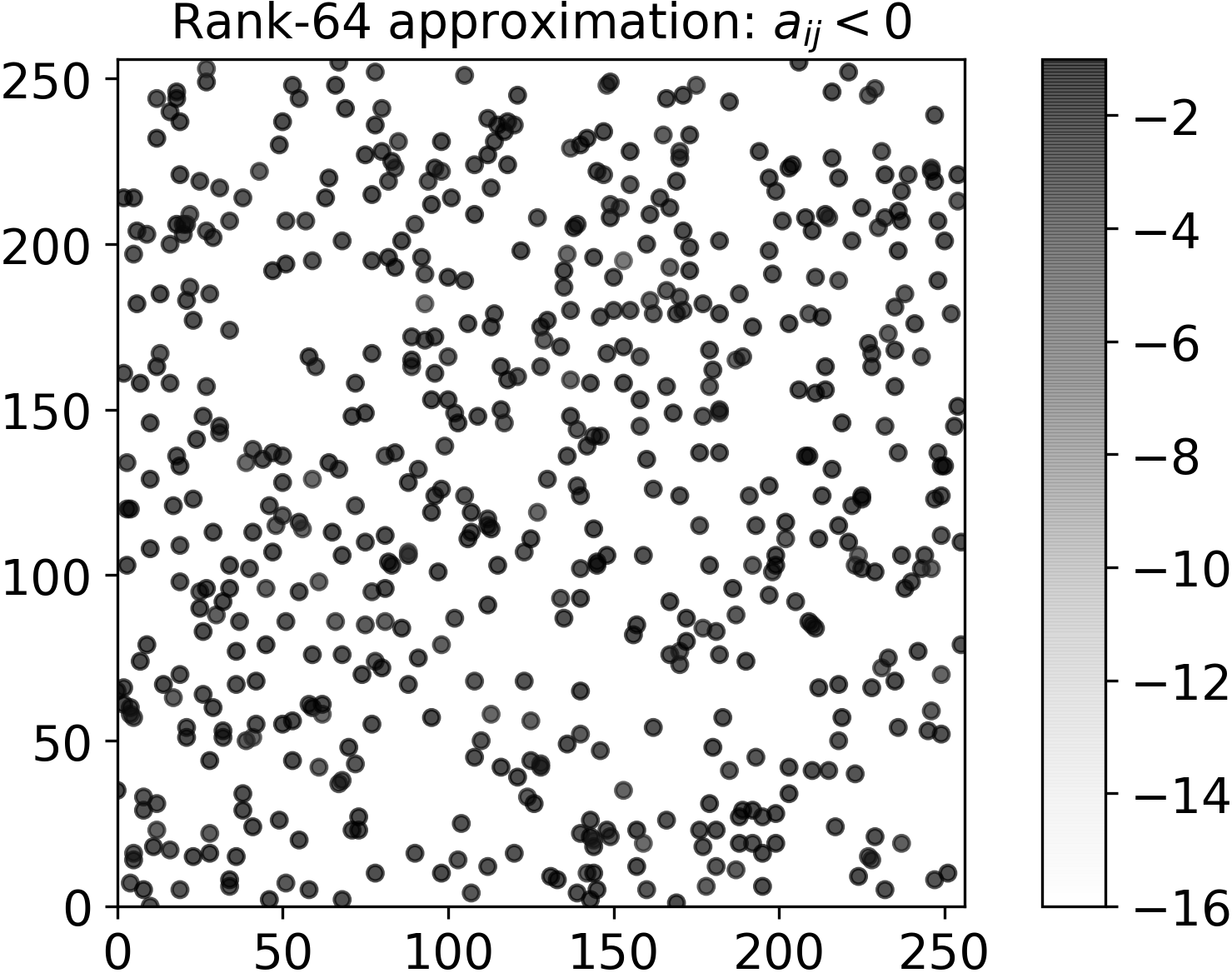}
	\caption{}
\end{subfigure}
\caption{Properties of a random $256 \times 256$ matrix with iid elements distributed uniformly on $[0, 1]$: normalized singular values~(a) and the magnitude (in log scale) of the negative elements of its best rank-$64$ approximation~(b).}
\label{num:fig:uniform_data}
\end{figure}

\begin{figure}[th]
\begin{subfigure}[b]{0.45\textwidth}
\centering
	\includegraphics[width=0.8\textwidth]{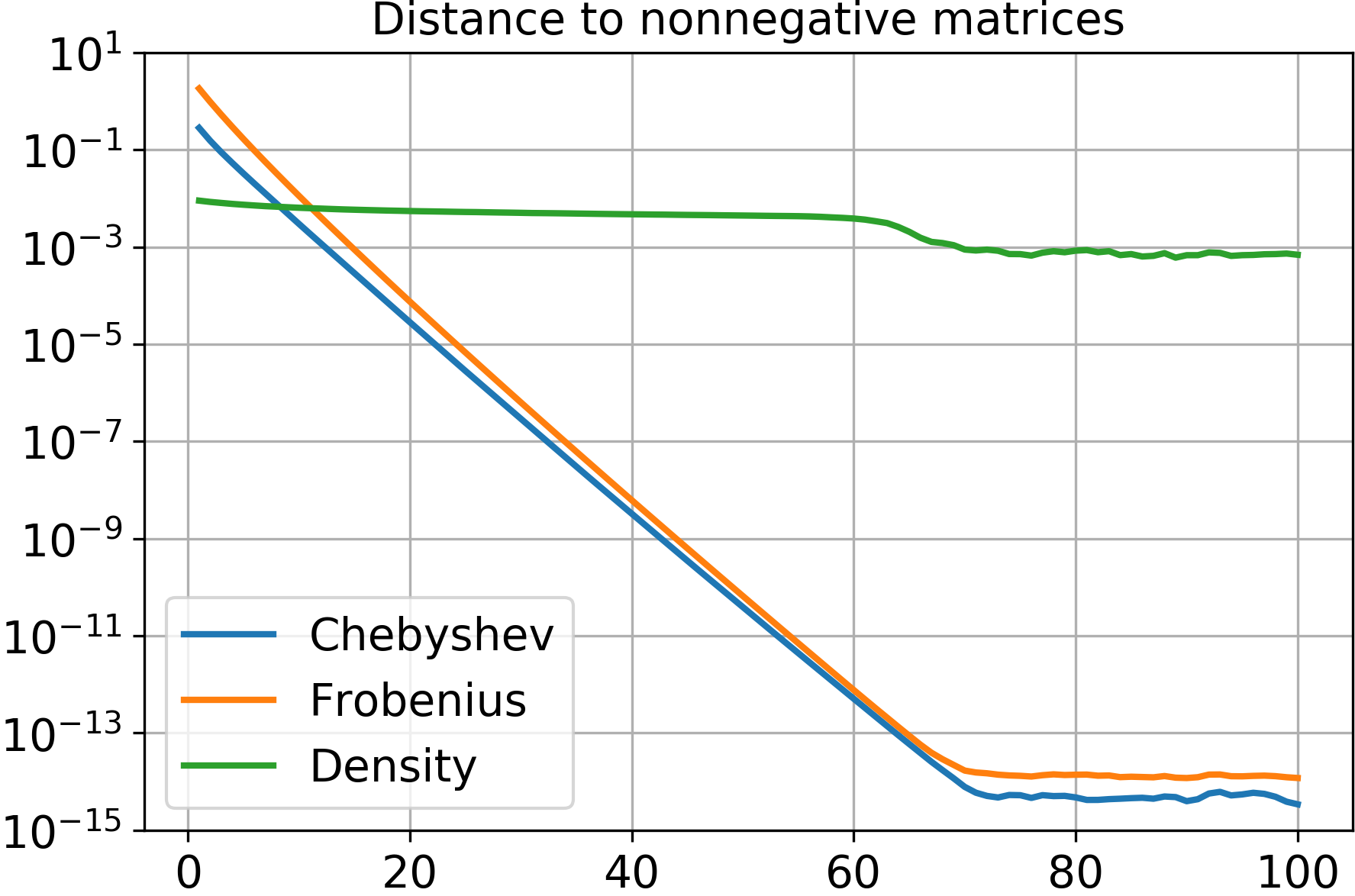}
	\caption{SVD}
\end{subfigure}\hspace{0.05\textwidth}
\begin{subfigure}[b]{0.45\textwidth}
\centering
	\includegraphics[width=0.8\textwidth]{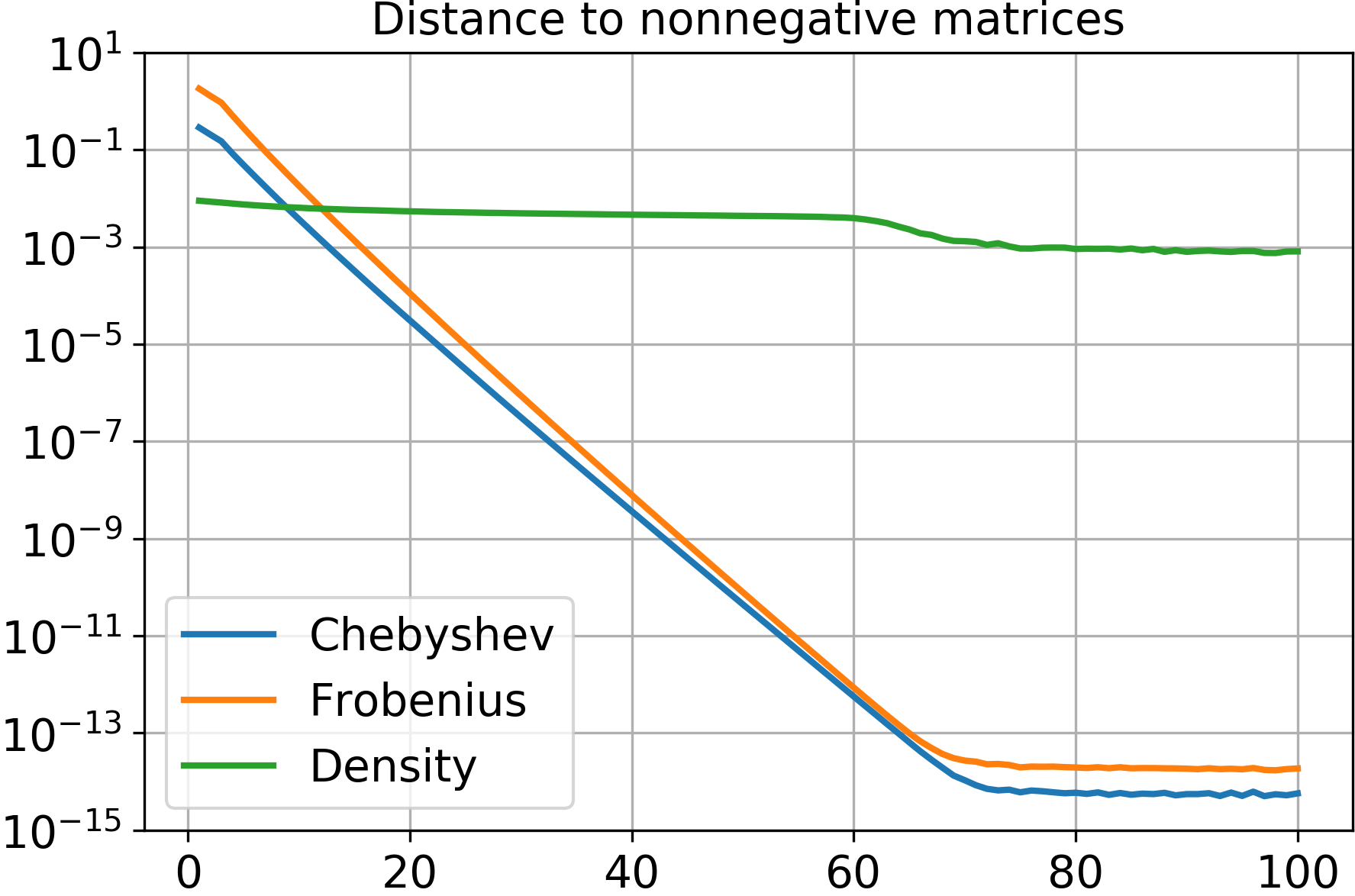}
	\caption{Tangent}
\end{subfigure}
\begin{subfigure}[b]{0.45\textwidth}
\centering
	\includegraphics[width=0.8\textwidth]{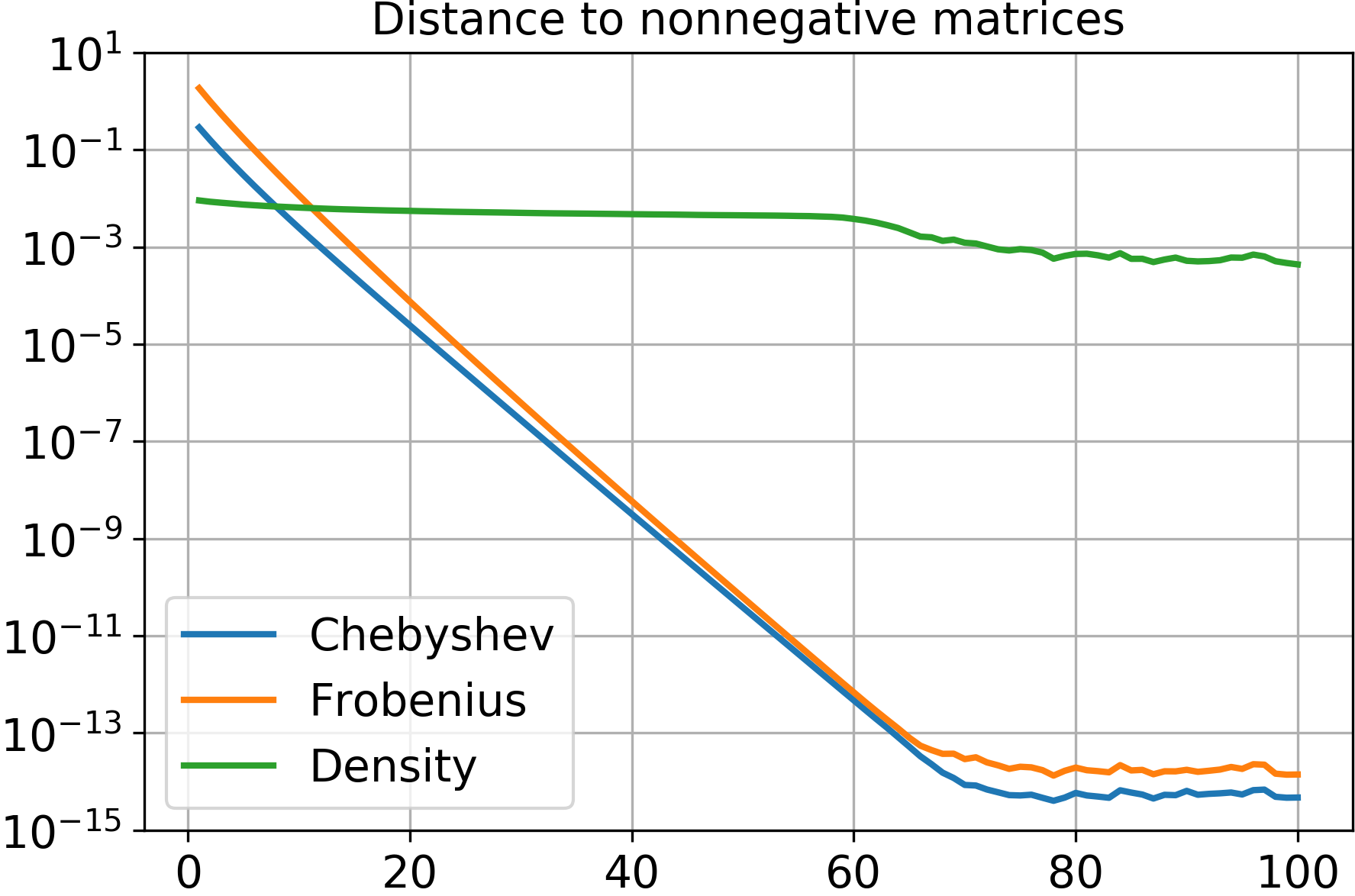}
	\caption{HMT(1, 70), $\mathcal{N}(0,1)$}
\end{subfigure}\hspace{0.05\textwidth}
\begin{subfigure}[b]{0.45\textwidth}
\centering
	\includegraphics[width=0.8\textwidth]{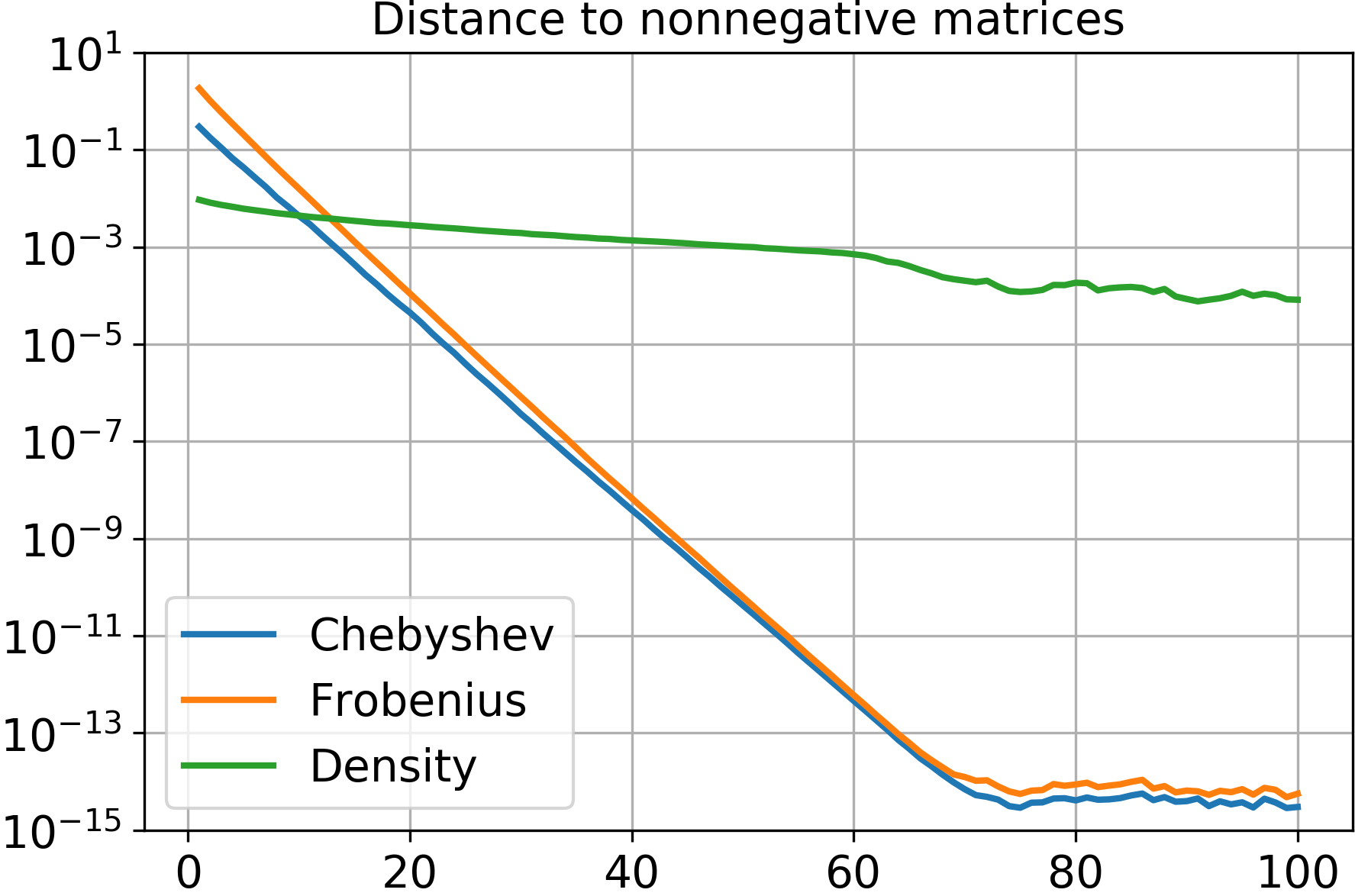}
	\caption{HMT(0, 70), $\mathrm{Rad}(0.2)$}
\end{subfigure}
\begin{subfigure}[b]{0.45\textwidth}
\centering
	\includegraphics[width=0.8\textwidth]{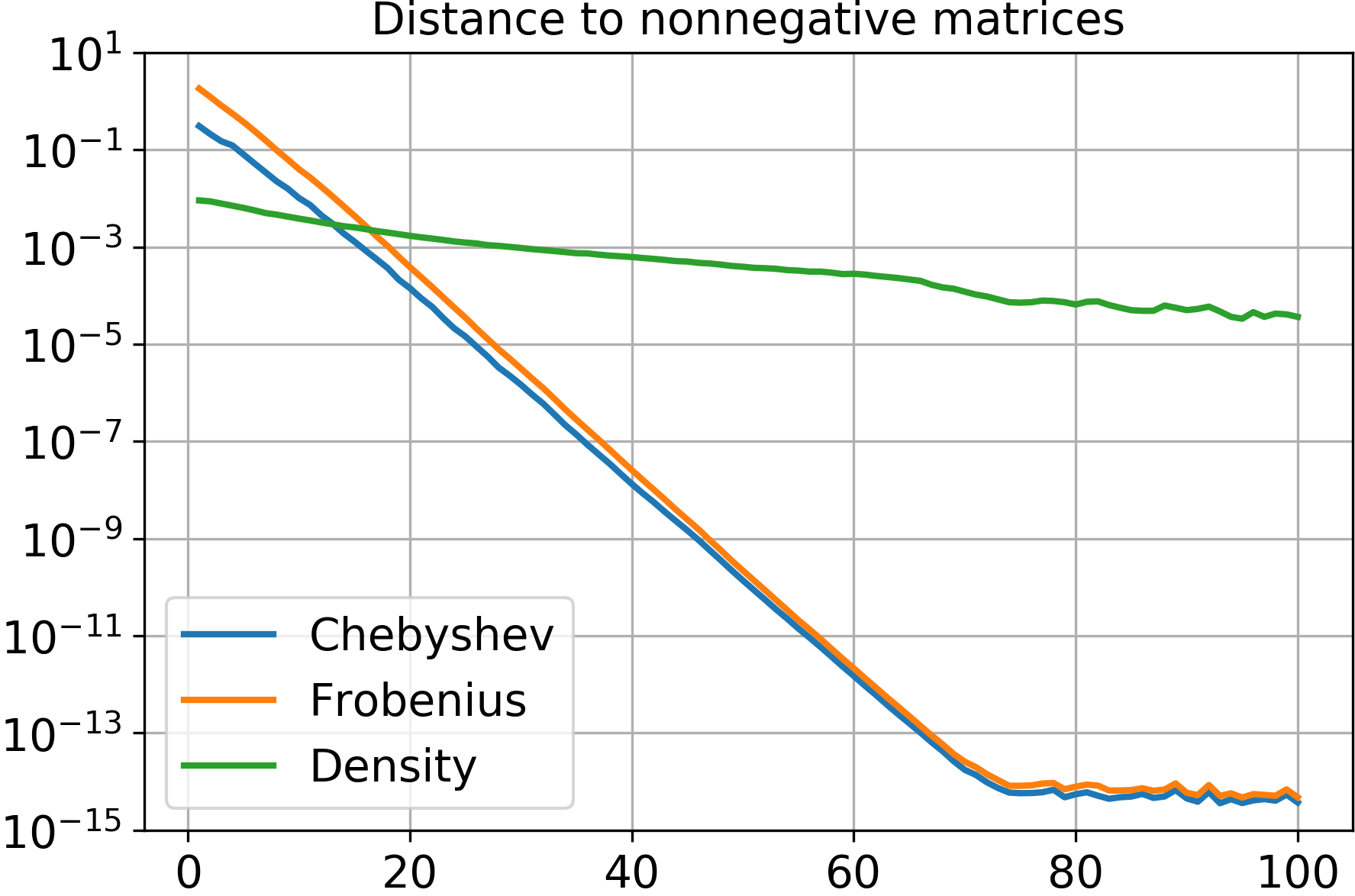}
	\caption{Tropp(70, 100), $\mathrm{Rad}(0.2)$}
\end{subfigure}\hspace{0.05\textwidth}
\begin{subfigure}[b]{0.45\textwidth}
\centering
	\includegraphics[width=0.8\textwidth]{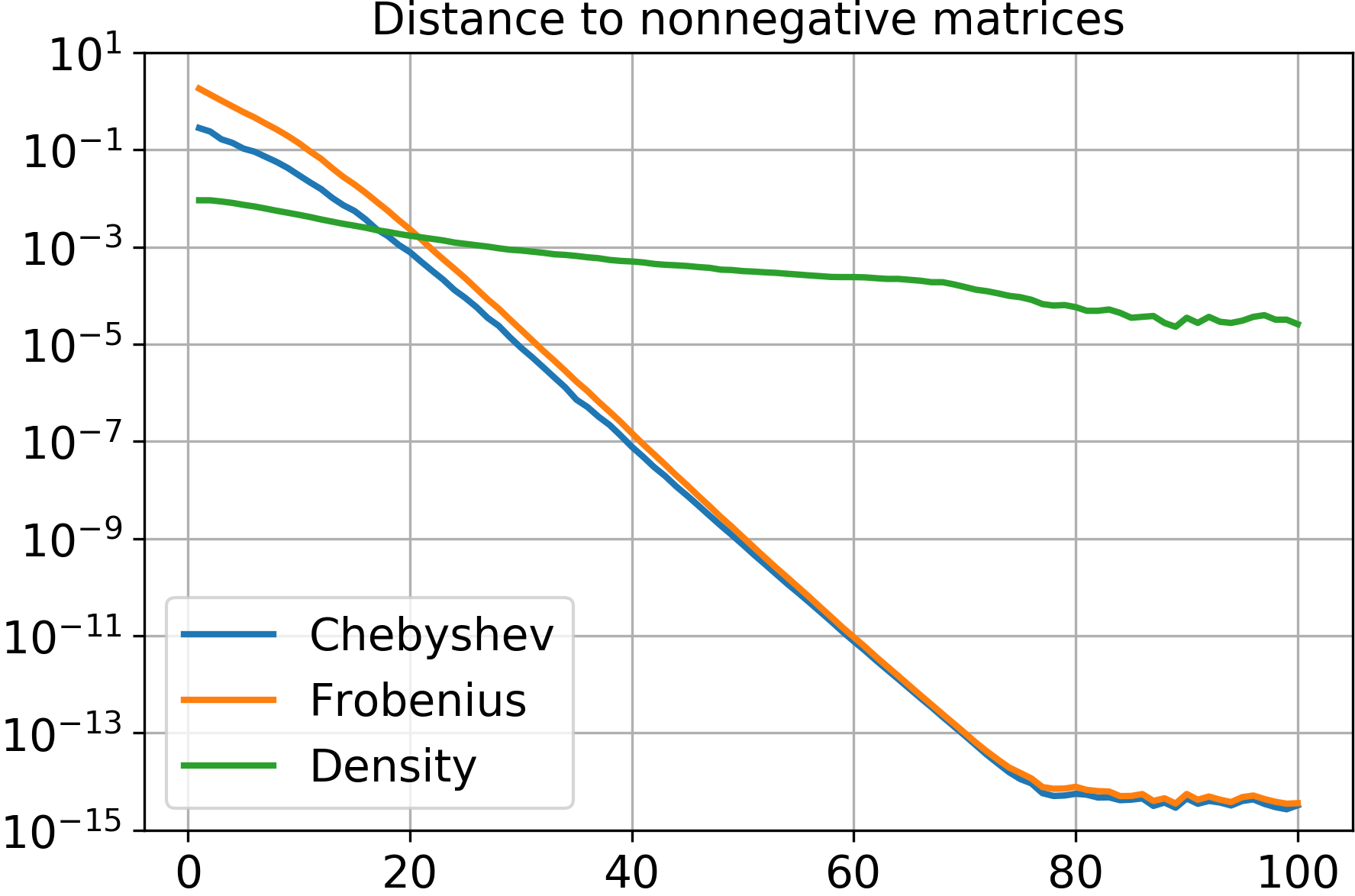}
	\caption{Tropp(70, 85), $\mathrm{Rad}(0.2)$}
\end{subfigure}
\begin{subfigure}[b]{0.45\textwidth}
\centering
	\includegraphics[width=0.8\textwidth]{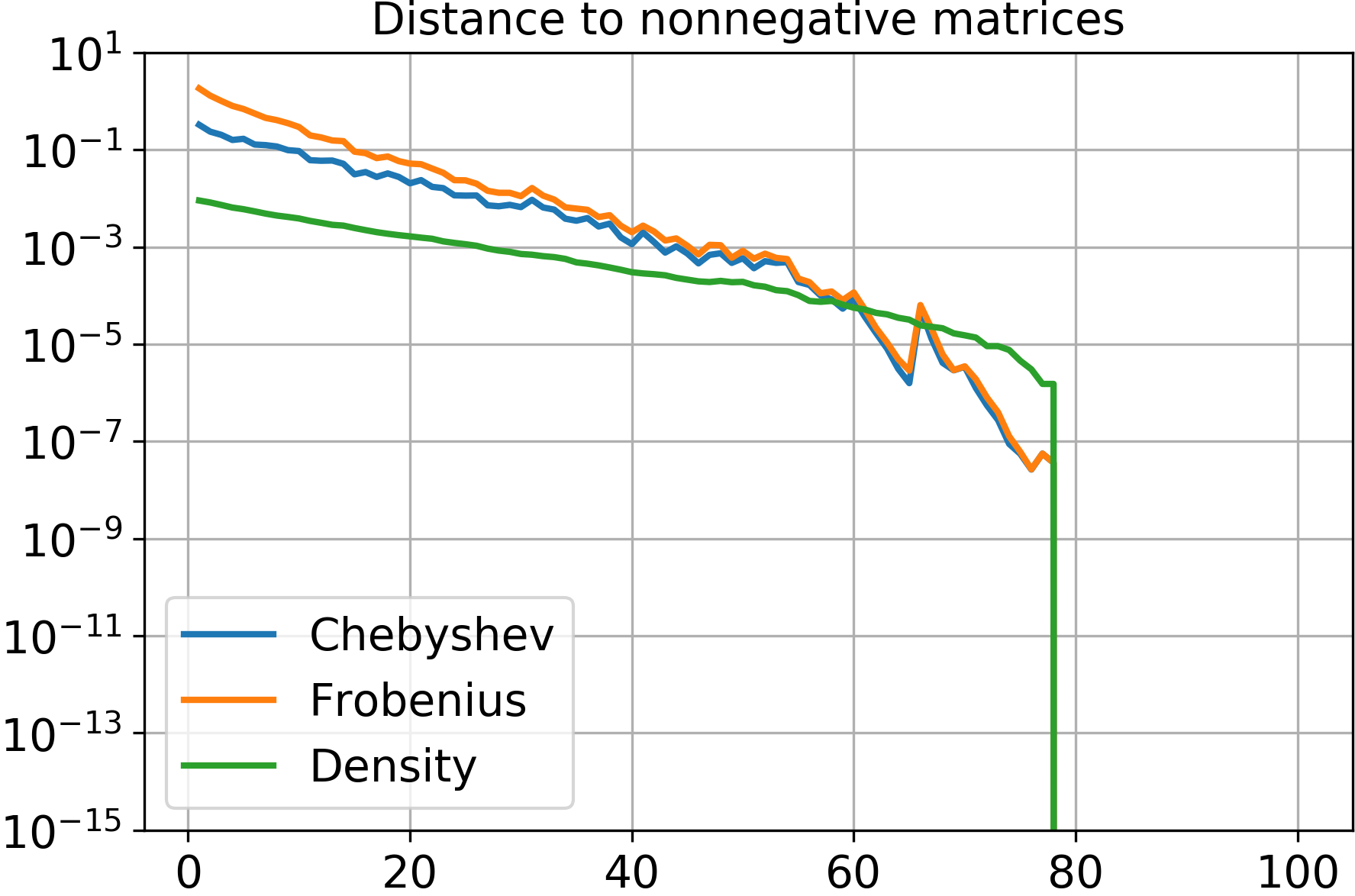}
	\caption{GN(150), $\mathrm{Rad}(0.2)$}
\end{subfigure}\hspace{0.05\textwidth}
\begin{subfigure}[b]{0.45\textwidth}
\centering
	\includegraphics[width=0.8\textwidth]{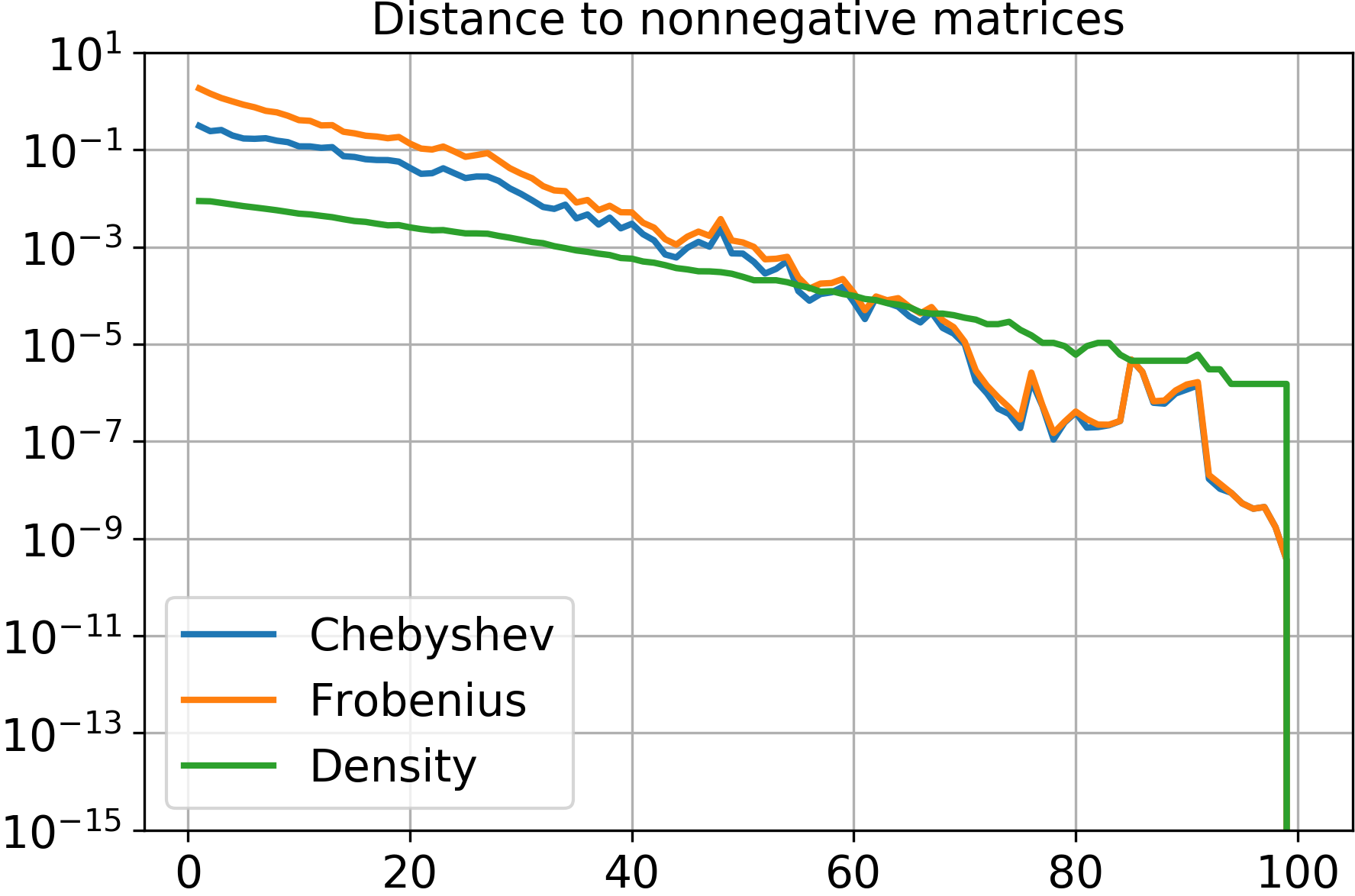}
	\caption{GN(120), $\mathrm{Rad}(0.2)$}
\end{subfigure}
\caption{Comparison of alternating projection methods for rank-$64$ nonnegative approximation of random $256 \times 256$ matrices with iid elements distributed uniformly on $[0, 1]$: the Frobenius and Chebyshev norms of the negative part and the density of negative elements over 100 iterations. The results are averaged over 10 trials.}
\label{num:fig:uniform_ap_comparison}
\end{figure}

In the first example, we consider random $256 \times 256$ matrices with independent identically distributed entries, distributed uniformly on $[0, 1]$, and try to approximate them with nonnegative rank-64 matrices. The best rank-64 approximation given by the truncated singular value decomposition contains many negative elements (see Fig.\ref{num:fig:uniform_data}), and we attempt to correct it using alternating projections. The results are presented in Tab.~\ref{num:tab:uniform} and Fig.~\ref{num:fig:uniform_ap_comparison}: the former contains the per-iteration computational complexities of each approach measured in flops, and the approximation errors in the Frobenius and Chebyshev (maximum) norms after 100 iterations; the latter shows the decay rate of the negative elements (we consider a value negative if it is below $-10^{-15}$). We see that the randomized approaches perform fewer operations per iteration than SVD and Tangent and have similar convergence properties. The only exception is GN: it keeps more large negative elements in the process and then abruptly makes them positive.

\begin{table}[h]\centering
\begin{tabular}{@{}llccc@{}}\toprule
    Method & Sketch & Flops per iter & Frobenius & Chebyshev \\
    \midrule
    Initial $\mathrm{SVD}_r$ & N/A & $3.5 \cdot 10^8$ & $3.07 \cdot 10^{-1}$ & $7.18 \cdot 10^{-1}$ \\
    \midrule
    SVD & N/A & $3.6 \cdot 10^8$ & $3.08 \cdot 10^{-1}$ & $7.17 \cdot 10^{-1}$ \\
    Tangent & N/A & $9.2 \cdot 10^7$ & $3.08 \cdot 10^{-1}$ & $7.19 \cdot 10^{-1}$ \\
    HMT(1, 70) & $\mathcal{N}(0,1)$ & $7.7 \cdot 10^7$ & $3.08 \cdot 10^{-1}$ & $7.14 \cdot 10^{-1}$ \\
    HMT(0, 70) & $\mathcal{N}(0,1)$ & $5.0 \cdot 10^7$ & $3.11 \cdot 10^{-1}$ & $7.14 \cdot 10^{-1}$ \\
    HMT(0, 70) & $\mathrm{Rad}$ & $4.4 \cdot 10^7$ & $3.10 \cdot 10^{-1}$ & $7.24 \cdot 10^{-1}$ \\
    HMT(0, 70) & $\mathrm{Rad}(0.2)$ & $4.0 \cdot 10^7$ & $3.10 \cdot 10^{-1}$ & $7.18 \cdot 10^{-1}$ \\
    Tropp(70, 100) & $\mathrm{Rad}(0.2)$ & $3.8 \cdot 10^7$ & $3.17 \cdot 10^{-1}$ & $7.47 \cdot 10^{-1}$ \\
    Tropp(70, 85) & $\mathrm{Rad}(0.2)$ & $3.6 \cdot 10^7$ & $3.30 \cdot 10^{-1}$ & $7.97 \cdot 10^{-1}$ \\
    GN(150) & $\mathrm{Rad}(0.2)$ & $2.0 \cdot 10^7$ & $3.40 \cdot 10^{-1}$ & $8.25 \cdot 10^{-1}$ \\
    GN(120) & $\mathrm{Rad}(0.2)$ & $1.8 \cdot 10^7$ & $3.60 \cdot 10^{-1}$ & $8.33 \cdot 10^{-1}$ \\
    \bottomrule\\
\end{tabular}
\caption{Comparison of alternating projection methods for rank-$64$ nonnegative approximation of random $256 \times 256$ matrices with iid elements distributed uniformly on $[0, 1]$: their computational complexities and relative errors in the Frobenius and Chebyshev norms after 100 iterations.}
\label{num:tab:uniform}
\end{table}

\clearpage
\subsection{Solution to Smoluchowski equation}
Our second example comes from the two-component Smoluchowski coagulation equation
\begin{equation}\label{num:eq:smol}
\begin{split}
    &\frac{\partial n(v_1, v_2, t)}{\partial t} = -n(v_1, v_2,t) \int_0^{\infty} \int_0^{\infty} K(u_1, u_2; v_1, v_2) n(u_1, u_2,t) du_1 du_2 \\
    & +\frac{1}{2}\int_0^{v_1}  \int_0^{v_2} K(v_1 - u_1, v_2 - u_2; u_1, u_2) n(v_1 - u_1, v_2 - u_2, t) n(u_1, u_2, t) du_1 du_2,
\end{split}
\end{equation}
which describes the evolution of the concentration function $n(v_1, v_2, t)$ of the two-component particles of size $(v_1, v_2)$ per unit volume. In the previous works \cite{smirnov2016fast, matveev2016tensor}, we showed that the corresponding initial-value problem can be solved by explicit time-integration in low-rank format for a wide range of coagulation kernels $K(u_1, u_2; v_1, v_2)$ and nonnegative initial conditions. This means that at every time instant $t$ the solution $n(v_1, v_2, t)$ is represented as a low-rank matrix, which accelerates computation. In some cases, analytical solutions are known \cite{fernandez2007exact}: the solution to Eq.~\eqref{num:eq:smol} with constant kernel 
\begin{equation*}
    K(u_1, u_2; v_1, v_2) \equiv K    
\end{equation*}
and the initial conditions
\begin{equation*}
    n(v_1, v_2, t = 0) = \sqrt{K} ab e^{-a v_1 - b v_2}    
\end{equation*}
is given by
\begin{equation}\label{num:eq:Analytical}
    n(v_1, v_2, t) = \sqrt{K} \frac{ab e^{-a v_1 - b v_2}} {(1 +  \sqrt{K} t /2)^2} I_0 \left(2 \sqrt{\frac{ab v_1 v_2 \sqrt{K} t}{ \sqrt{K} t + 2}}\right),
\end{equation}
where $a, b > 0$ are arbitrary positive numbers and $I_0$ is the modified Bessel function of order zero. It was proved in \cite{matveev2016tensor} that \eqref{num:eq:Analytical}, discretized on any equidistant rectangular grid, can be approximated with accuracy $\varepsilon$ by a matrix of rank $O(\log 1 / \varepsilon)$ that is independent of the grid.

\begin{figure}[th]
\begin{subfigure}[b]{0.45\textwidth}
\centering
	\includegraphics[width=0.75\textwidth]{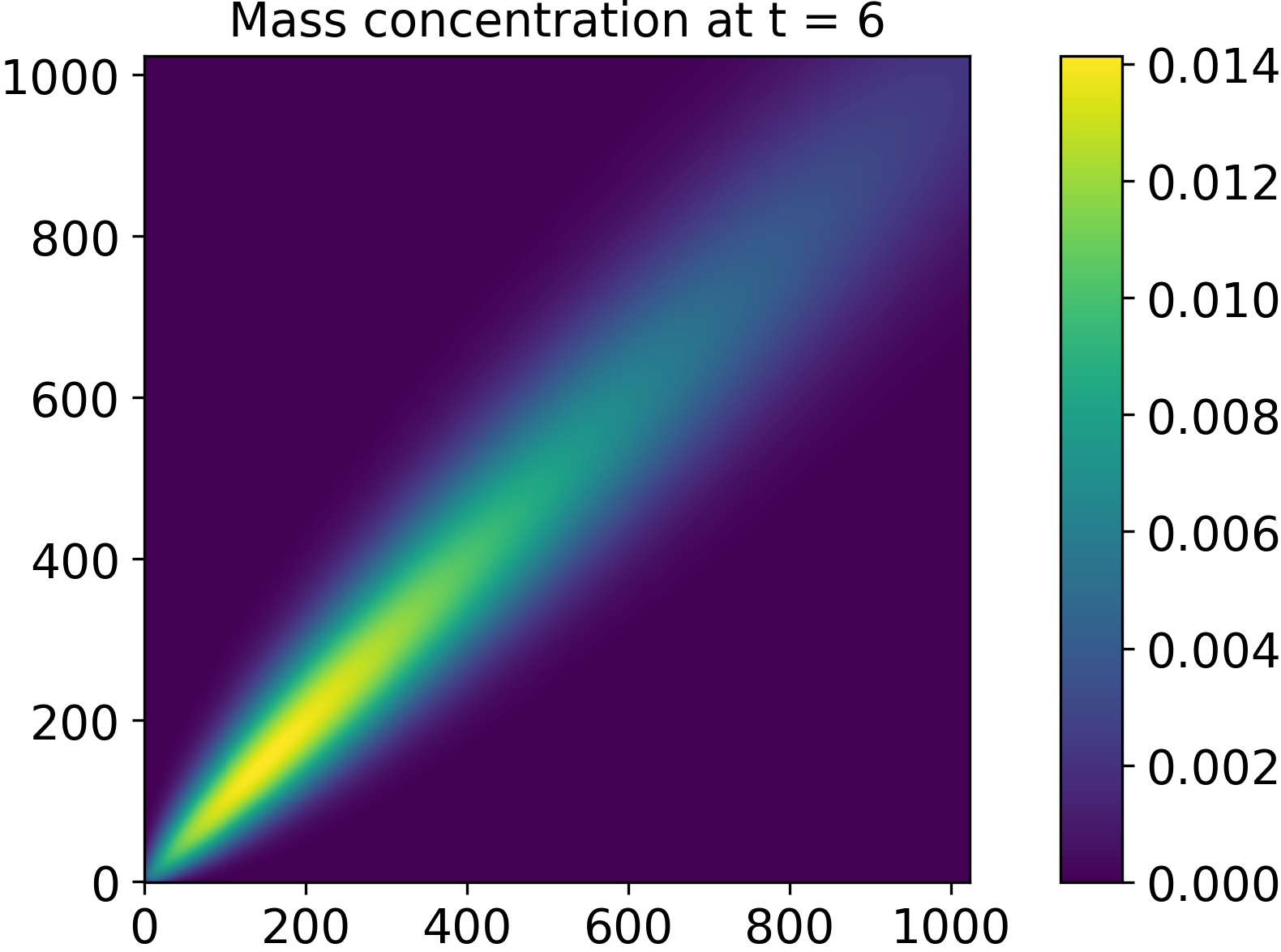}
	\caption{}
\end{subfigure}\hspace{0.05\textwidth}%
\begin{subfigure}[b]{0.45\textwidth}
\centering
	\includegraphics[width=0.8\textwidth]{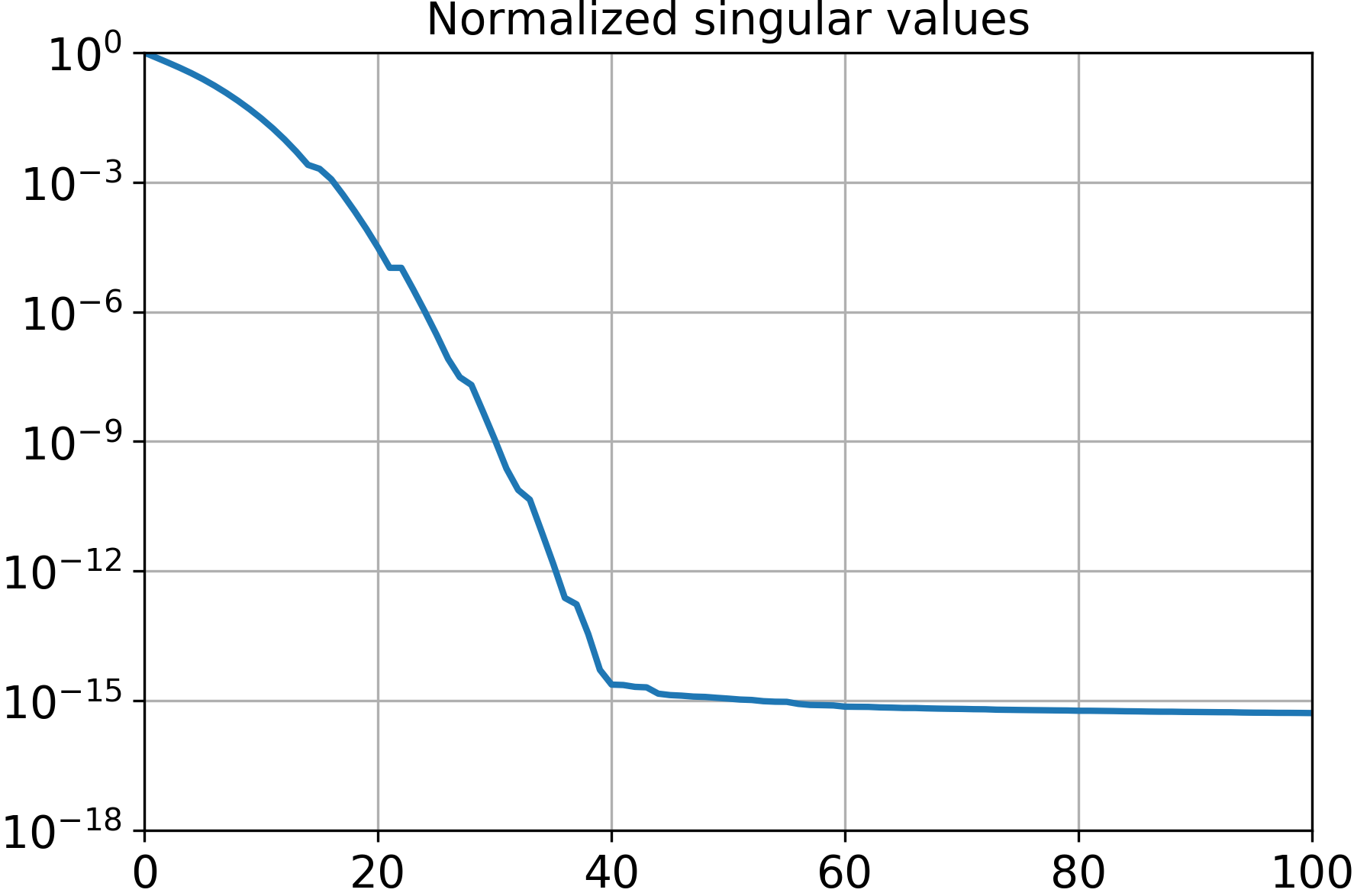}
	\caption{}
\end{subfigure}
\caption{Properties of the $1024 \times 1024$ solution of the Smoluchowski equation at $t = 6$: (a)~the solution itself and (b)~the normalized singular values.}
\label{num:fig:smolukh_data}
\end{figure}

For our numerical experiments, we set $K = 100$, $a = b = 1$, choose an equidistant rectangular grid with step $0.1$, and study rank-50 approximations of the $1024 \times 1024$ discretized mass-concentration function 
\begin{equation*}
    m(v_1, v_2, t) \equiv(v_1 + v_2) \cdot n(v_1, v_2, t)
\end{equation*}
corresponding to the solution \eqref{num:eq:Analytical} at $t = 6$. In Fig.~\ref{num:fig:smolukh_data}, we demonstrate the heatmap of $m(v_1, v_2, t)$ and the plot of its normalized singular values, which decay rapidly in agreement with \cite{matveev2016tensor}.

\begin{figure}[th]
\centering
\begin{subfigure}[b]{0.45\textwidth}
\centering
	\includegraphics[width=0.8\textwidth]{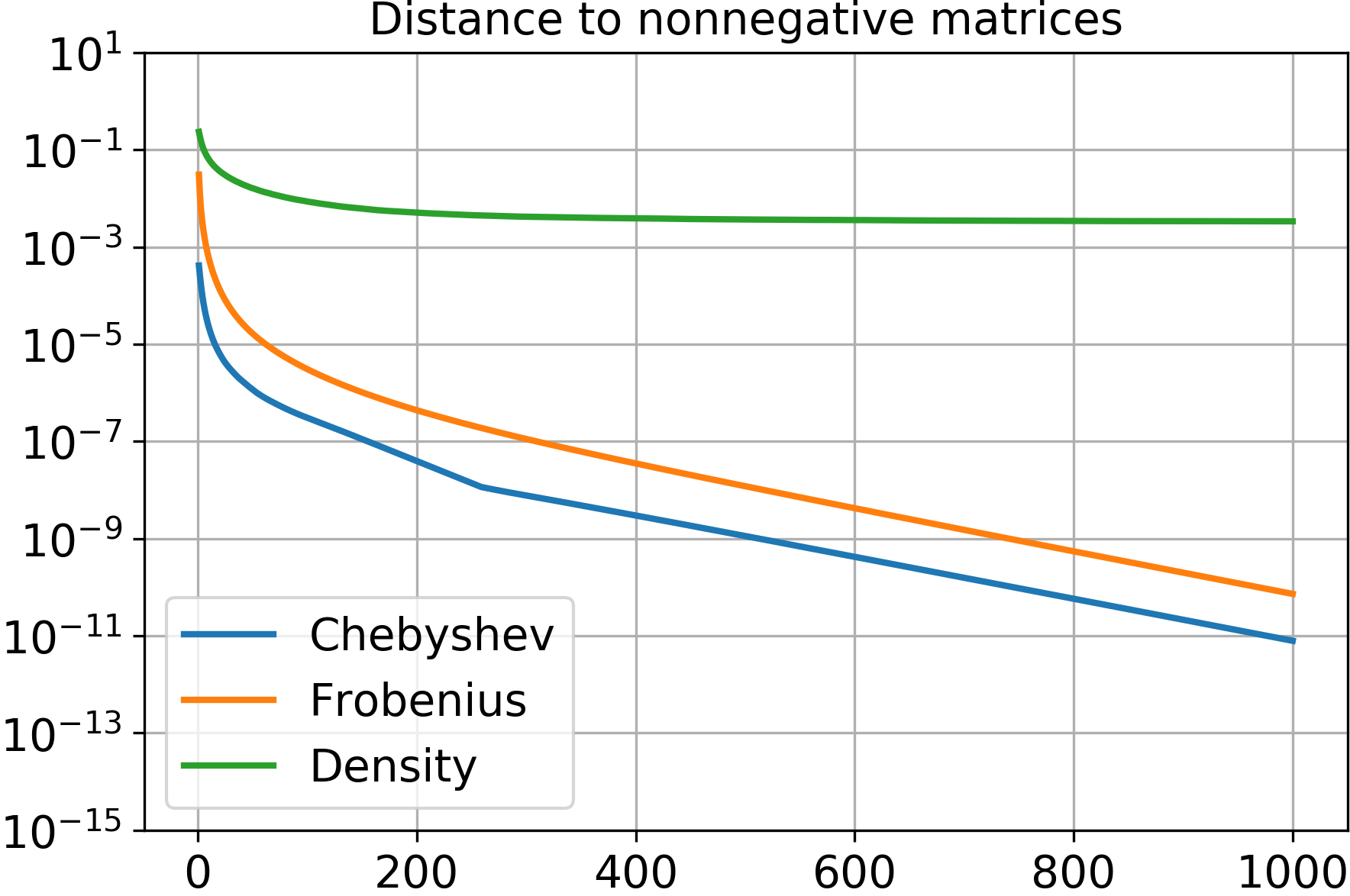}
	\caption{SVD}
\end{subfigure}\hspace{0.05\textwidth}%
\begin{subfigure}[b]{0.45\textwidth}
\centering
	\includegraphics[width=0.8\textwidth]{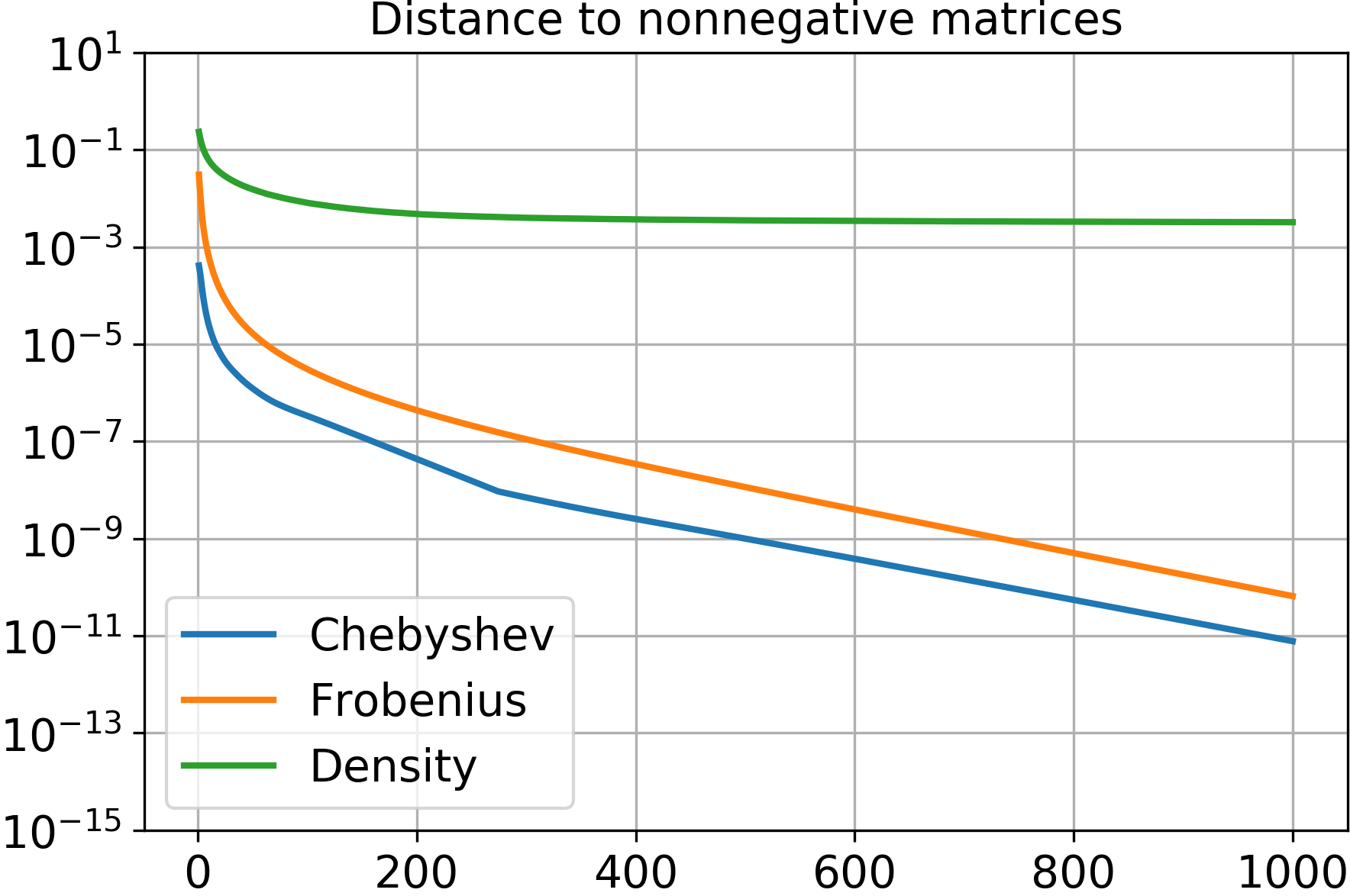}
	\caption{Tangent}
\end{subfigure}
\begin{subfigure}[b]{0.45\textwidth}
\centering
	\includegraphics[width=0.8\textwidth]{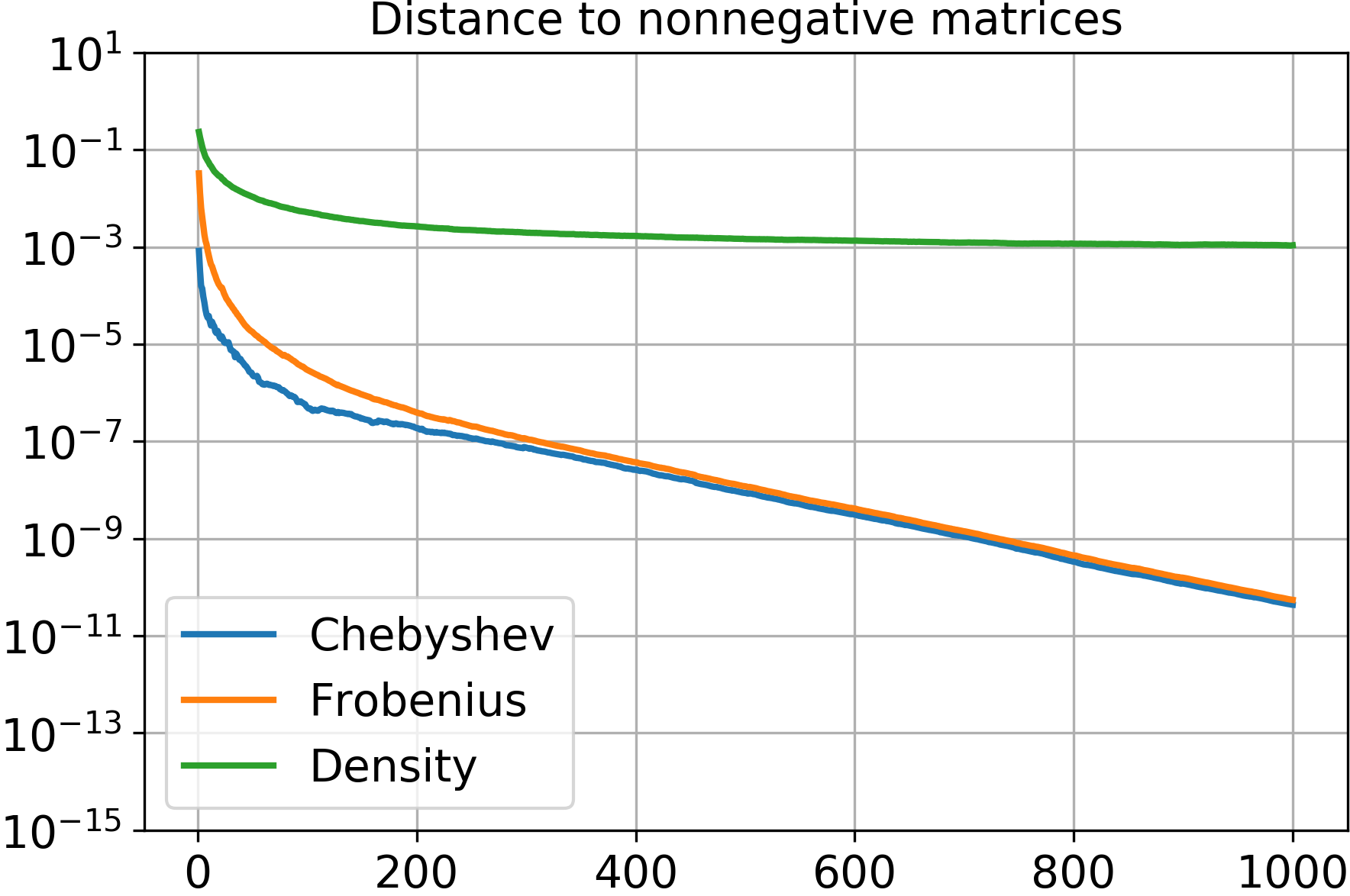}
	\caption{HMT(0, 15), $\mathrm{Rad}(0.2)$}
\end{subfigure}\hspace{0.05\textwidth}%
\begin{subfigure}[b]{0.45\textwidth}
\centering
	\includegraphics[width=0.8\textwidth]{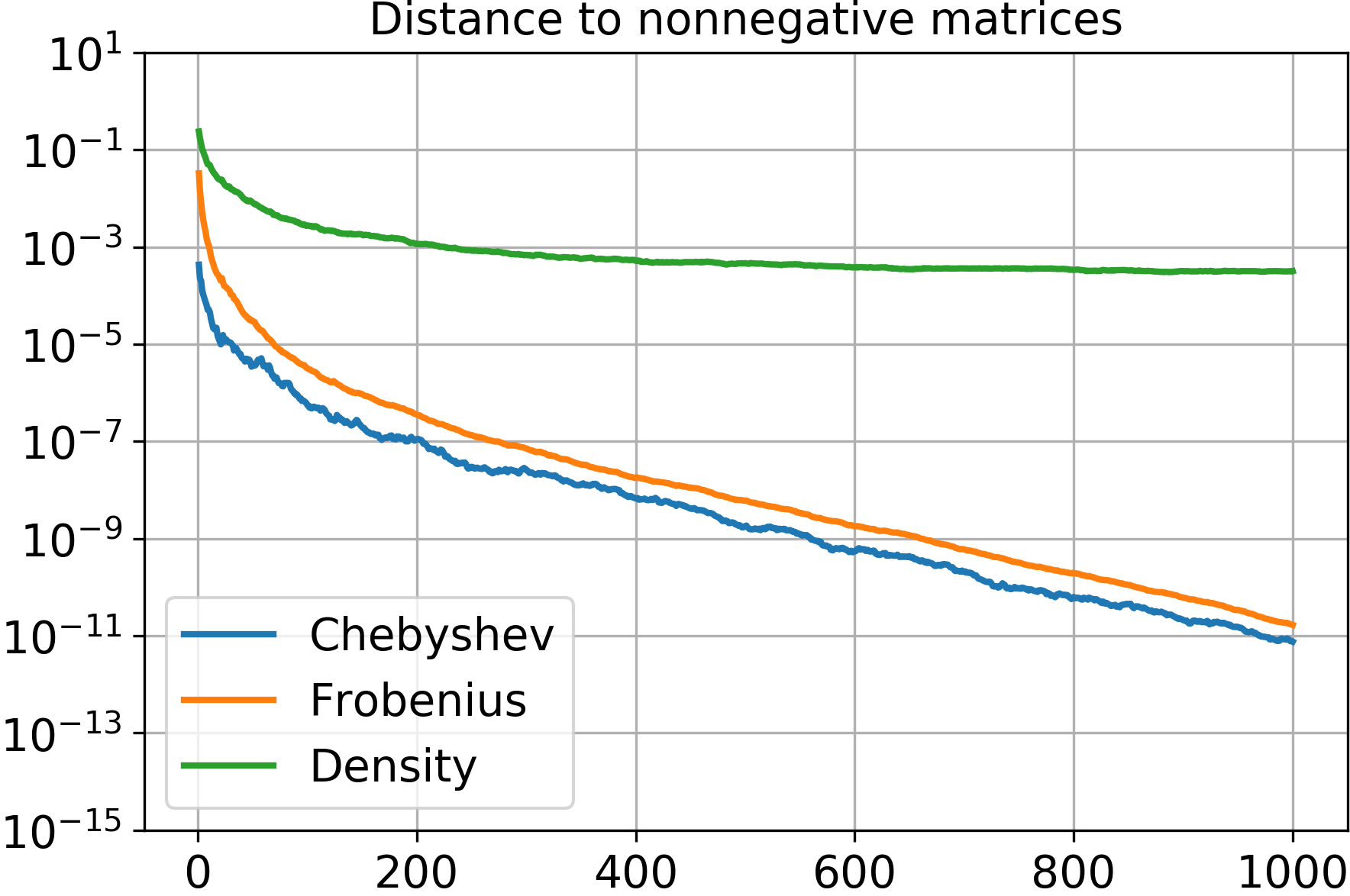}
	\caption{Tropp(15, 25), $\mathrm{Rad}(0.2)$}
\end{subfigure}
\begin{subfigure}[b]{0.45\textwidth}
\centering
	\includegraphics[width=0.8\textwidth]{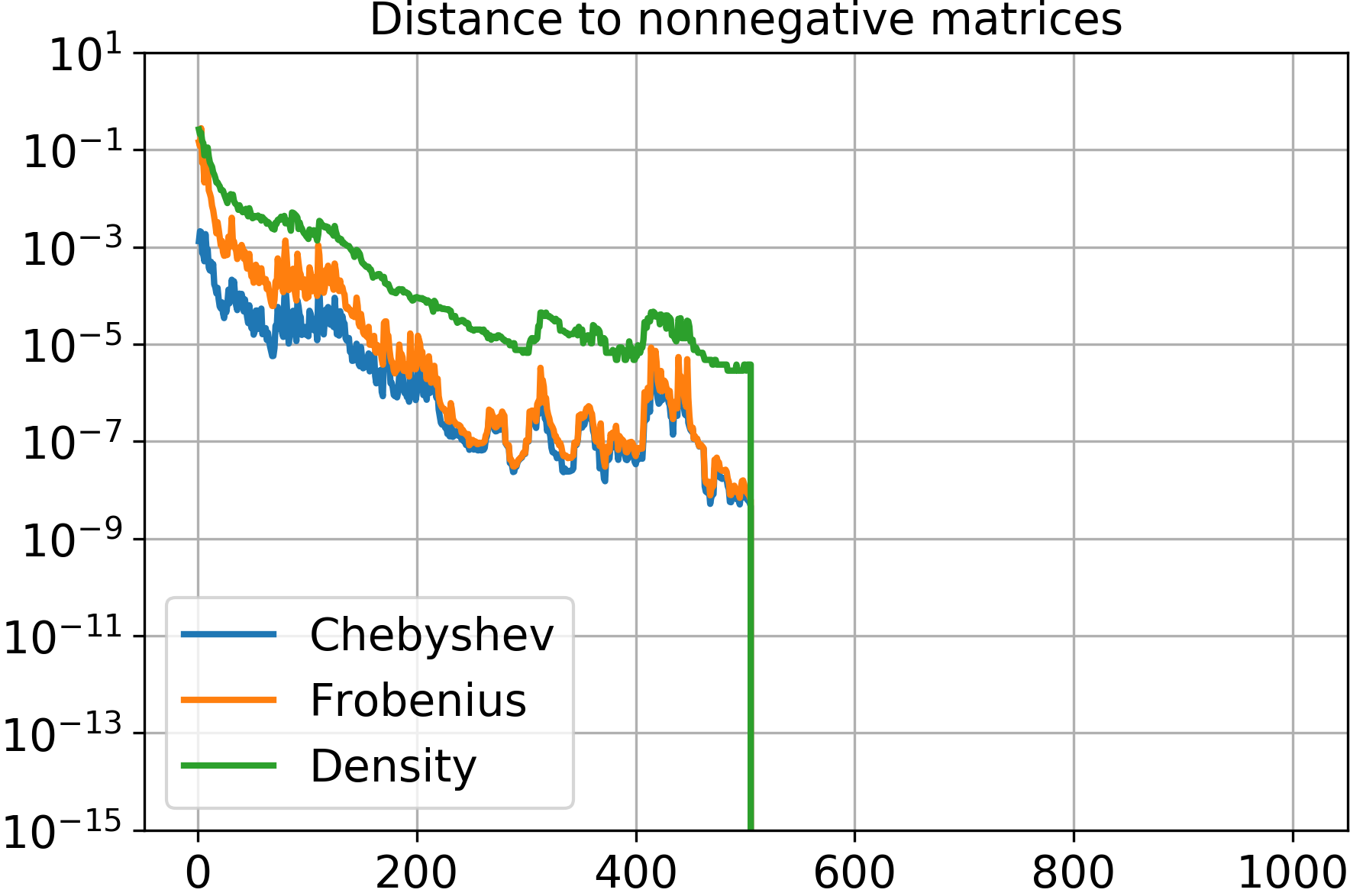}
	\caption{GN(40), $\mathrm{Rad}(0.2)$}
\end{subfigure}
\caption{Comparison of alternating projection methods for rank-$10$ nonnegative approximation of the $1024 \times 1024$ solution of the Smoluchowski equation: the Frobenius and Chebyshev norms of the negative part and the density of negative elements over 1000 iterations.}
\label{num:fig:smolukh_ap_distance}
\end{figure}

\begin{figure}[th]
\begin{subfigure}[b]{\textwidth}
\centering
	\includegraphics[width=0.3\textwidth]{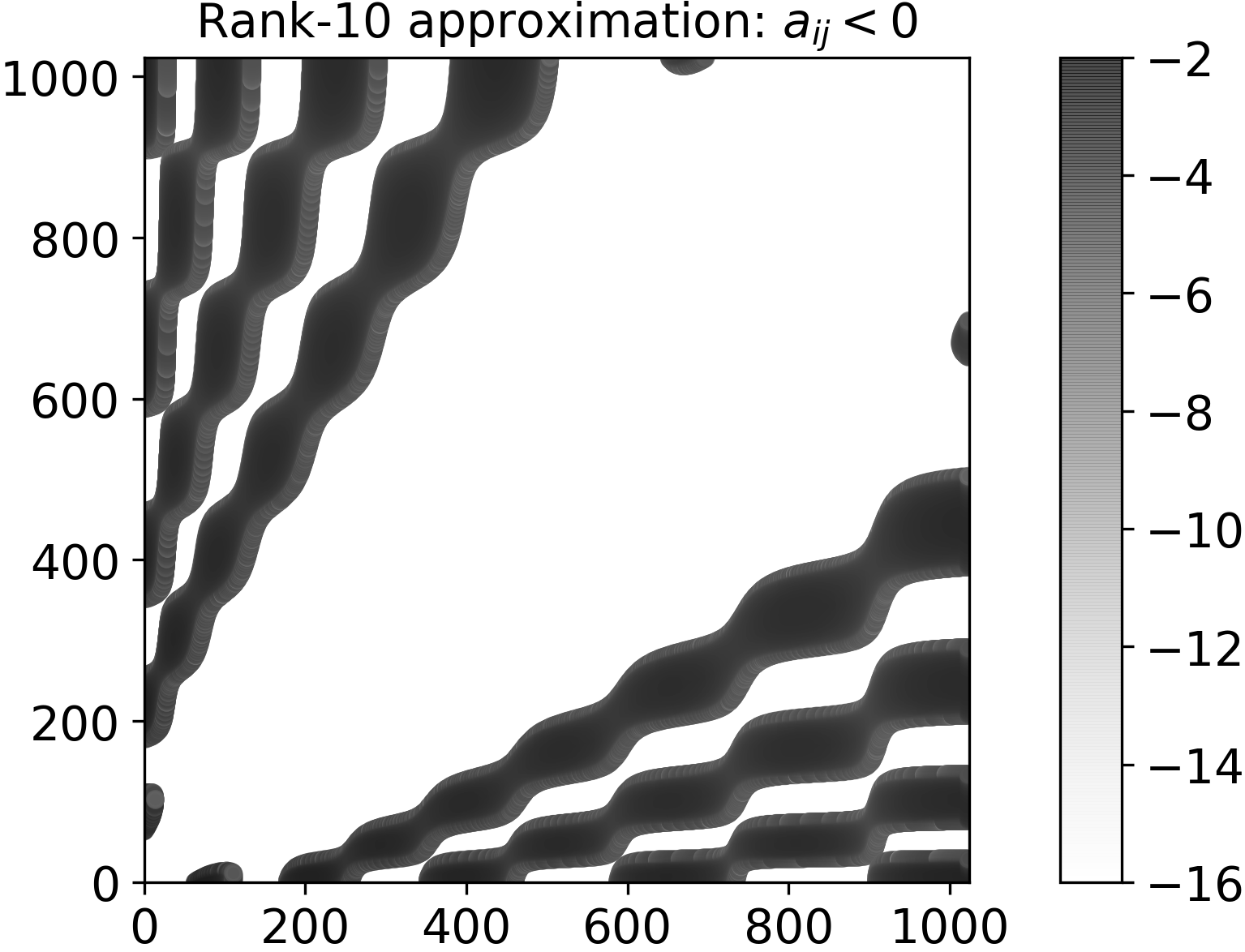}\hspace{0.05\textwidth}
	\includegraphics[width=0.3\textwidth]{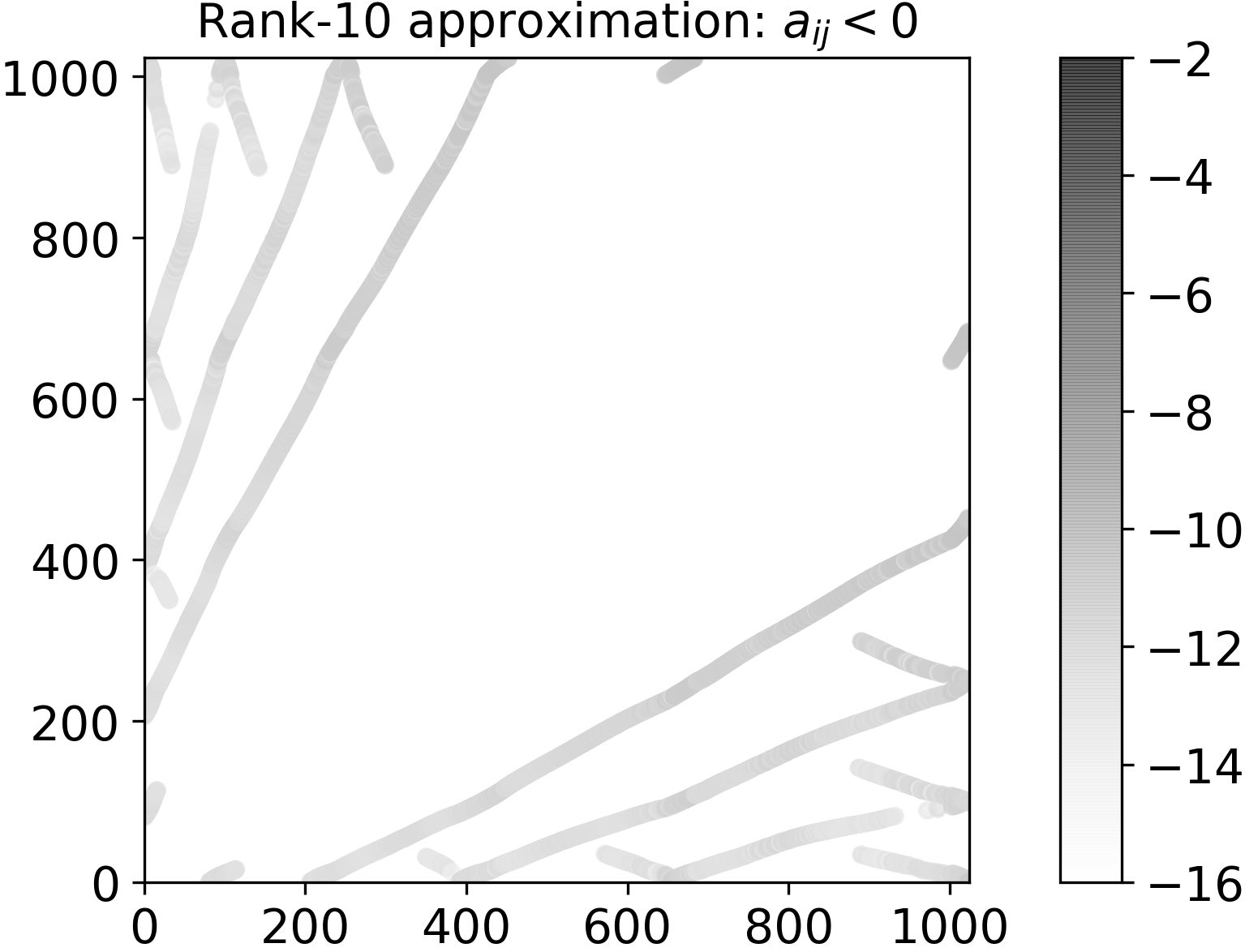}
	\caption{SVD: iterations 0 and 1000}
\end{subfigure}
\begin{subfigure}[b]{\textwidth}
\centering
	\includegraphics[width=0.3\textwidth]{pics/smolukh/smolukh_svd_negels_beg.png}\hspace{0.05\textwidth}
	\includegraphics[width=0.3\textwidth]{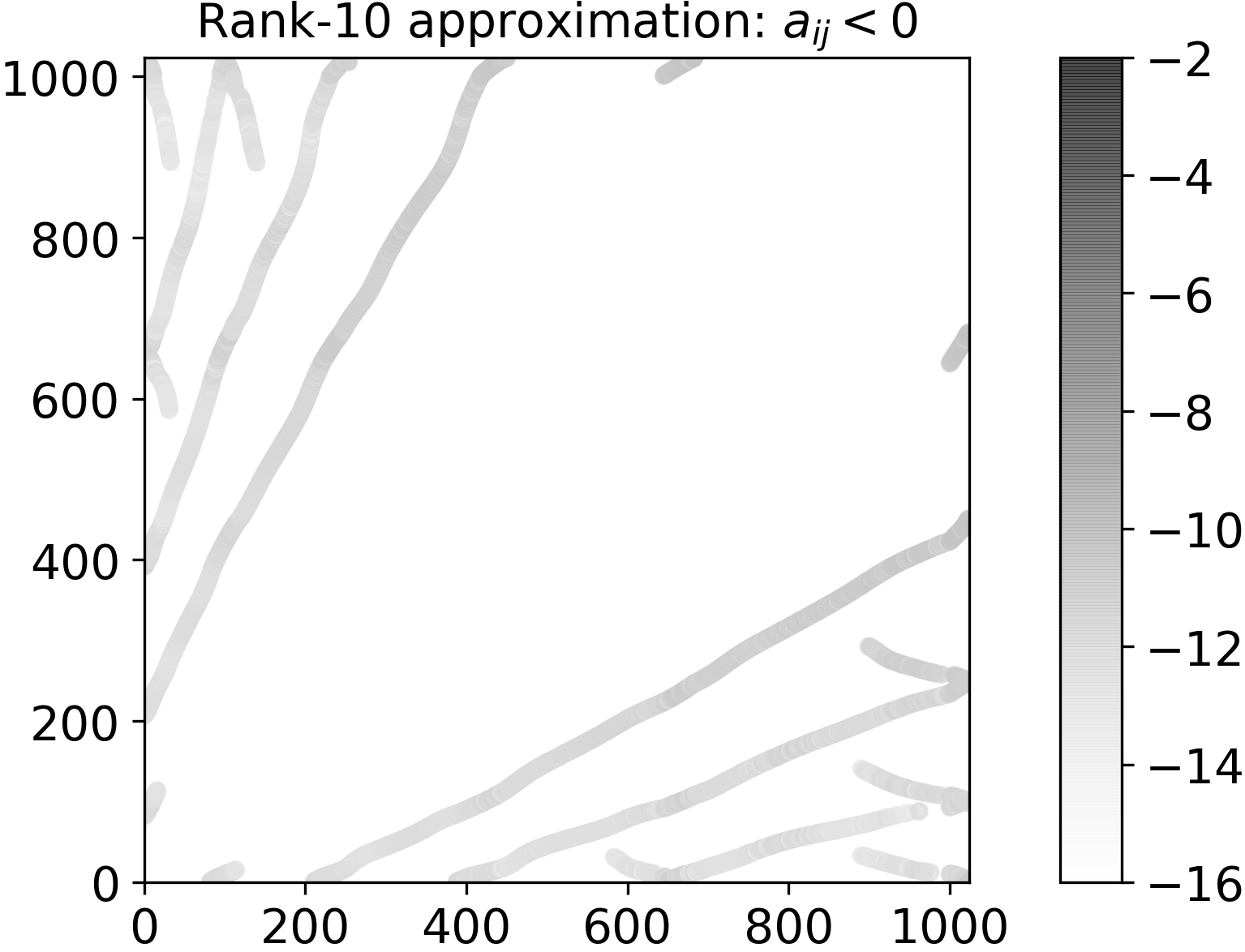}
	\caption{Tangent: iterations 0 and 1000}
\end{subfigure}
\begin{subfigure}[b]{\textwidth}
\centering
	\includegraphics[width=0.3\textwidth]{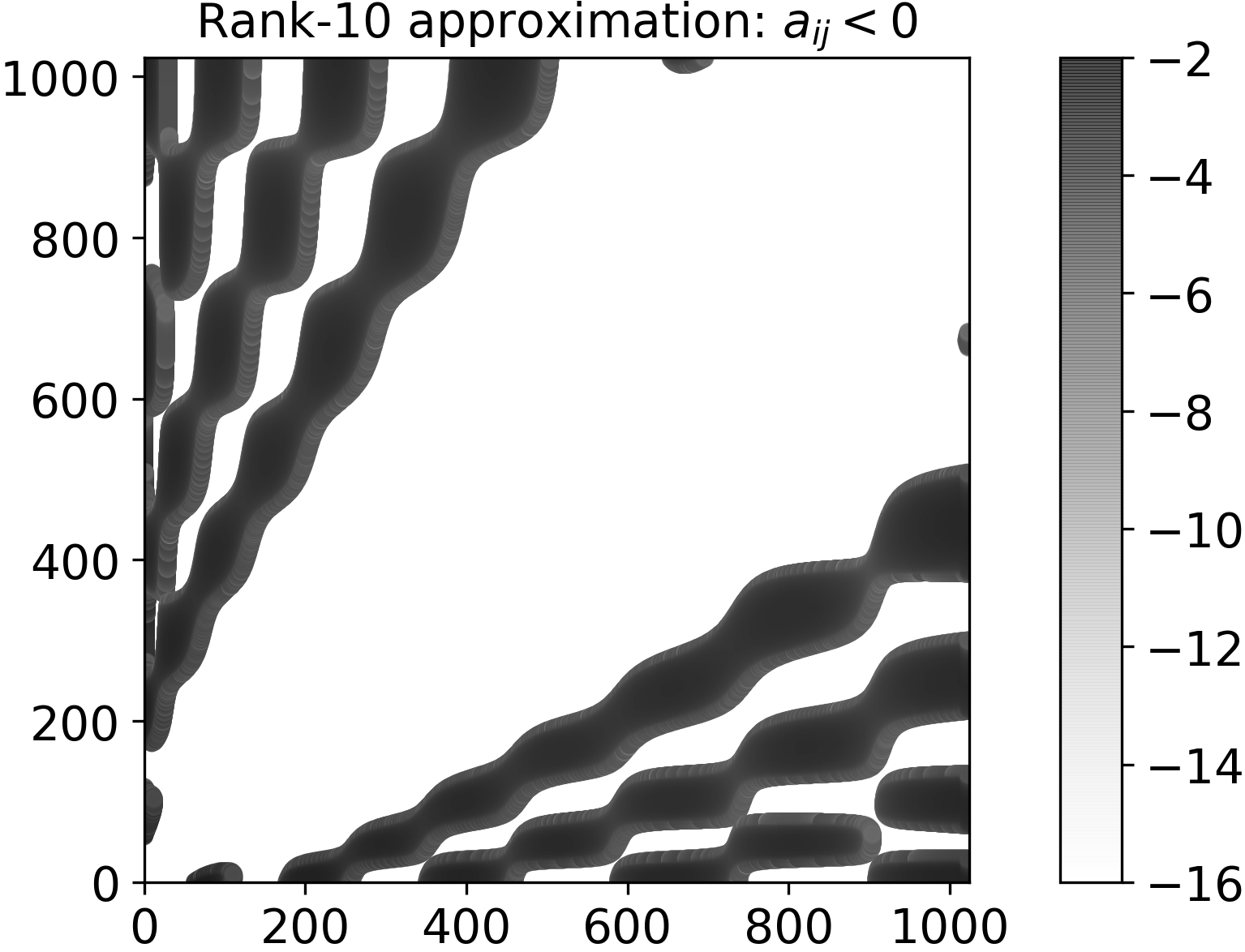}\hspace{0.05\textwidth}
	\includegraphics[width=0.3\textwidth]{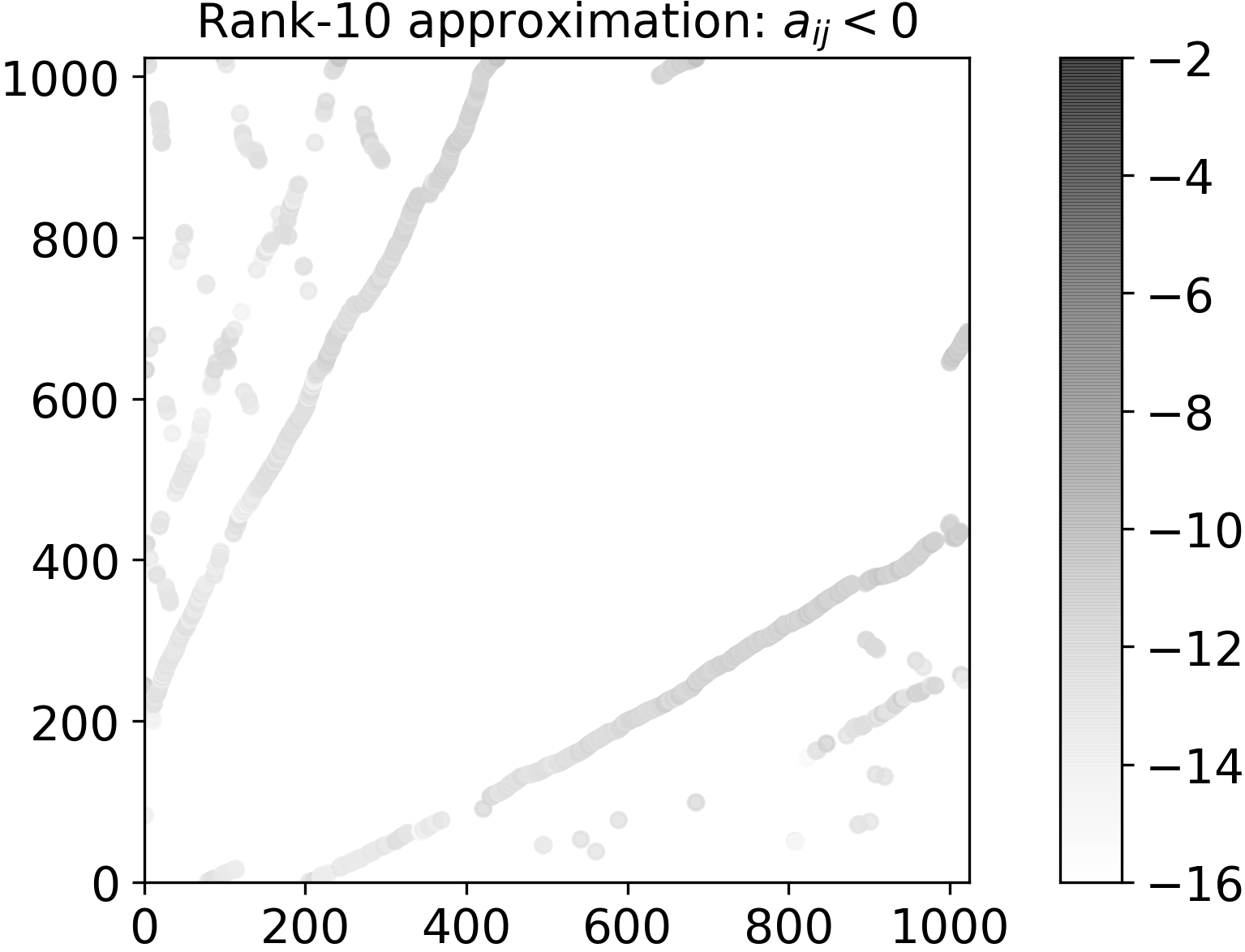}
	\caption{HMT(0, 15), $\mathrm{Rad}(0.2)$: iterations 0 and 1000}
\end{subfigure}
\begin{subfigure}[b]{\textwidth}
\centering
	\includegraphics[width=0.3\textwidth]{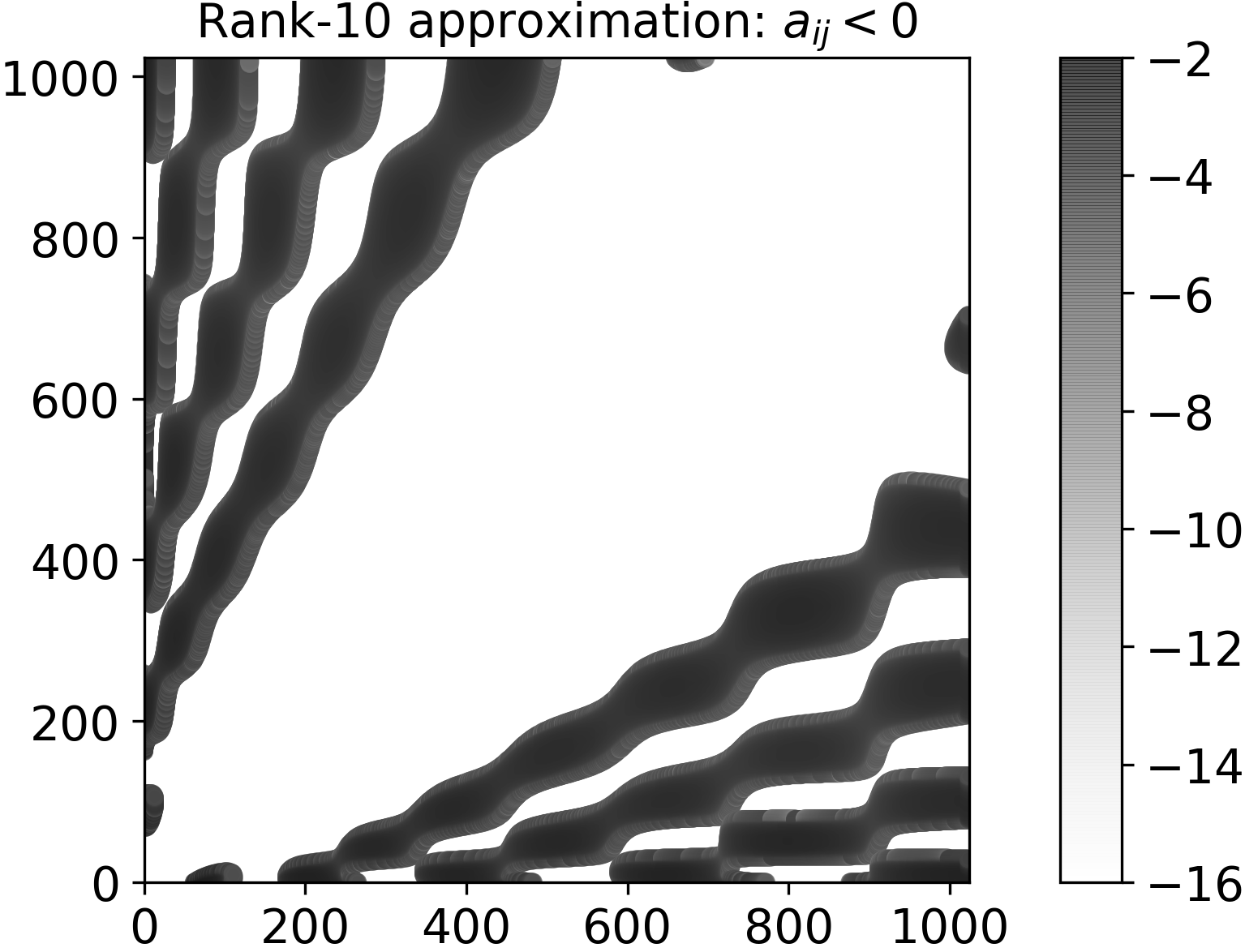}\hspace{0.05\textwidth}
	\includegraphics[width=0.3\textwidth]{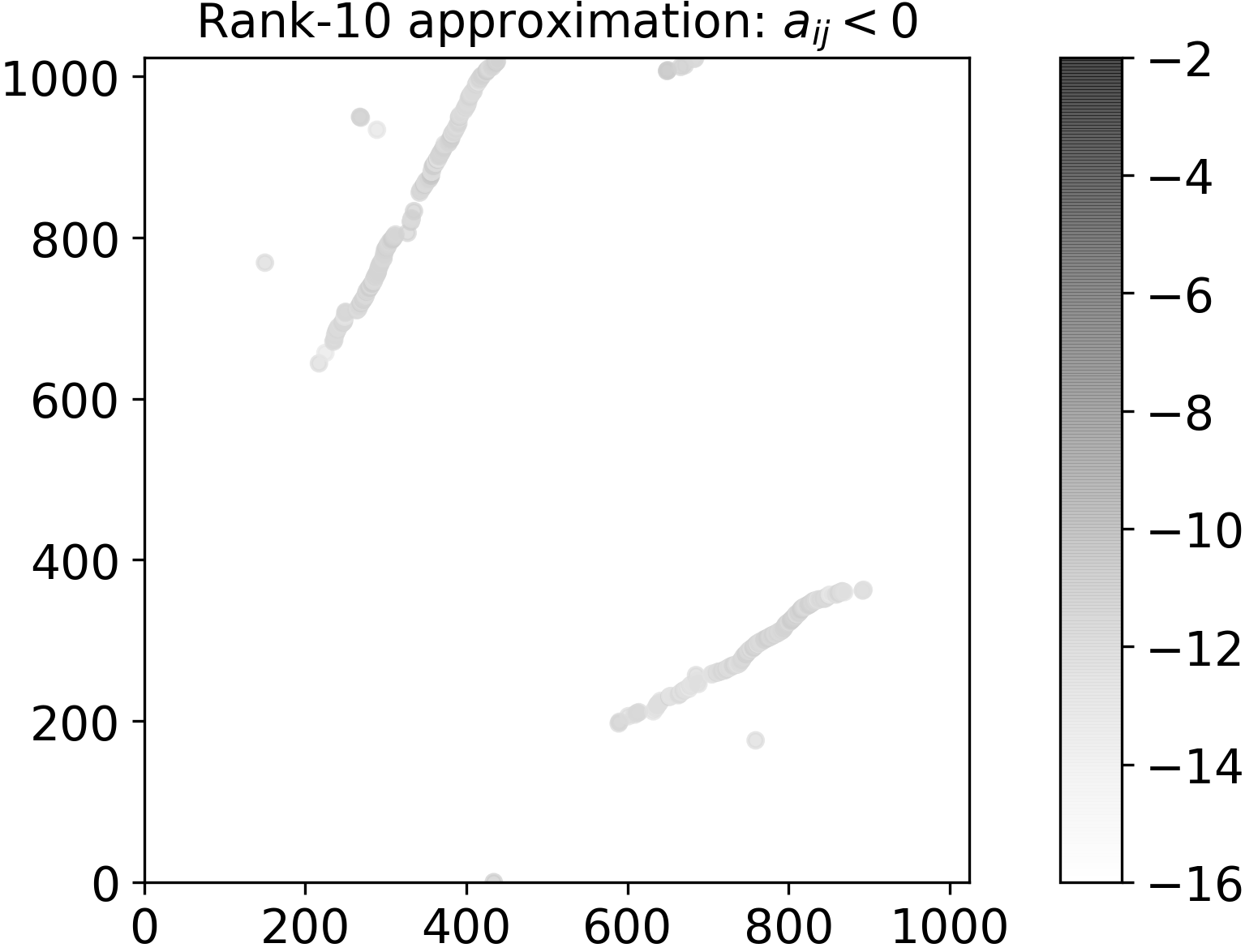}
	\caption{Tropp(15, 25), $\mathrm{Rad}(0.2)$: iterations 0 and 1000}
\end{subfigure}
\begin{subfigure}[b]{\textwidth}
\centering
	\includegraphics[width=0.3\textwidth]{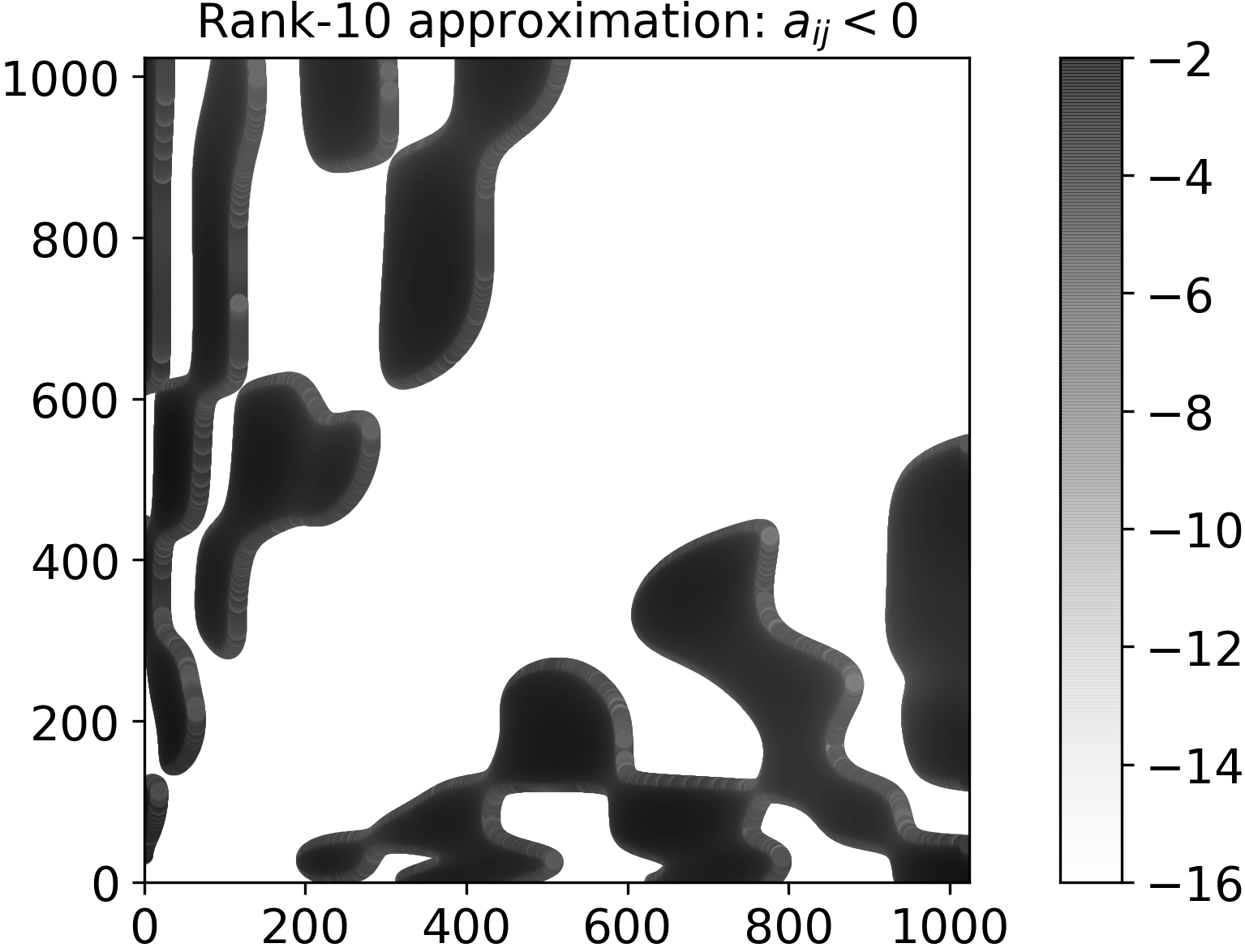}\hspace{0.05\textwidth}
	\includegraphics[width=0.3\textwidth]{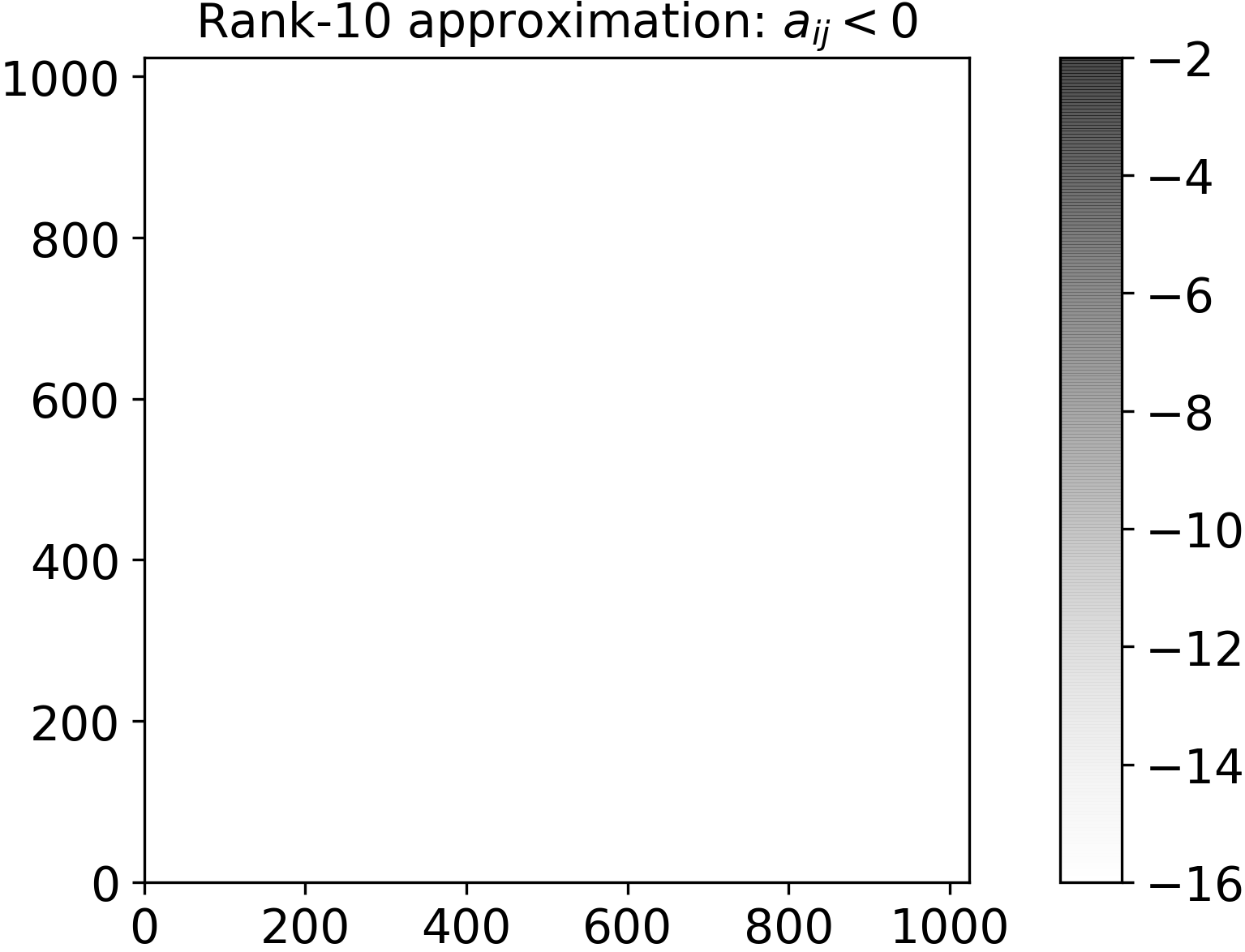}
	\caption{GN(40), $\mathrm{Rad}(0.2)$: iterations 0 and 1000}
\end{subfigure}
\caption{Comparison of alternating projection methods for rank-$10$ nonnegative approximation of the $1024 \times 1024$ solution of the Smoluchowski equation: the magnitude (in log scale) of the negative elements of the rank-$10$ approximation initially~(left) and after 1000 iterations~(right).}
\label{num:fig:smolukh_negels}
\end{figure}

Unlike the two previous examples, where we always started with the best low-rank approximation, here we use different initial low-rank approximations, according to the alternating projection method. In Tab.~\ref{num:tab:smolukh} and Fig.~\ref{num:fig:smolukh_ap_distance} we compare the performance of the discussed approaches: once again, randomized approaches are faster than deterministic ones and show similar convergence. The GN method eliminates the negative elements a lot sooner than the others, but its relative error in the Frobenius norm is 5 times higher. In Fig.~\ref{num:fig:smolukh_negels}, we show how the negative elements disappear after 1000 alternating projection iterations: HMT and Tropp leave the matrix with fewer negative elements than SVD and Tangent, and GN removes them completely. Also note how the initial low-rank approximation in GN has a distinct negative pattern.

\begin{table}[h]\centering
\begin{tabular}{@{}llccccc@{}}
    \toprule
    Method & Sketch & Flops (init/per iter) & Frobenius (init/res) & Chebyshev (init/res)\\
    \midrule
    SVD & N/A & $2.3 \cdot 10^{10} / 2.3 \cdot 10^{10}$ & $2.39 \cdot 10^{-2} / 2.70 \cdot 10^{-2}$ & $1.19 \cdot 10^{-1} / 1.48 \cdot 10^{-1}$ \\
    Tangent & N/A & $2.3 \cdot 10^{10} / 6.5 \cdot 10^{7}$ & $2.39 \cdot 10^{-2} / 2.72 \cdot 10^{-2}$ & $1.19 \cdot 10^{-1} / 1.49 \cdot 10^{-1}$ \\
    HMT(0, 15) & $\mathrm{Rad}(0.2)$ & $3.7 \cdot 10^{7} / 5.8 \cdot 10^{7}$ & $2.43 \cdot 10^{-2} / 2.75 \cdot 10^{-2}$ & $1.19 \cdot 10^{-1} / 1.57 \cdot 10^{-1}$ \\
    Tropp(15, 25) & $\mathrm{Rad}(0.2)$ & $1.2 \cdot 10^{7} / 3.3 \cdot 10^{7}$ & $2.43 \cdot 10^{-2} / 2.87 \cdot 10^{-2}$ & $1.33 \cdot 10^{-1} / 1.60 \cdot 10^{-1}$ \\
    GN(40) & $\mathrm{Rad}(0.2)$ & $1.2 \cdot 10^{7} / 3.3 \cdot 10^{7}$ & $8.47 \cdot 10^{-2} / 1.83 \cdot 10^{-1}$ & $1.54 \cdot 10^{-1} / 3.72 \cdot 10^{-1}$ \\
    \bottomrule\\
\end{tabular}
\caption{Comparison of alternating projection methods for rank-$10$ nonnegative approximation of the $1024 \times 1024$ solution of the Smoluchowski equation: their computational complexities and relative errors in the Frobenius and Chebyshev norms after 1000 iterations.}
\label{num:tab:smolukh}
\end{table}

\clearpage
\subsection{Images}
\begin{figure}[th]
\begin{subfigure}[b]{0.45\textwidth}
\centering
	\includegraphics[width=0.6\textwidth]{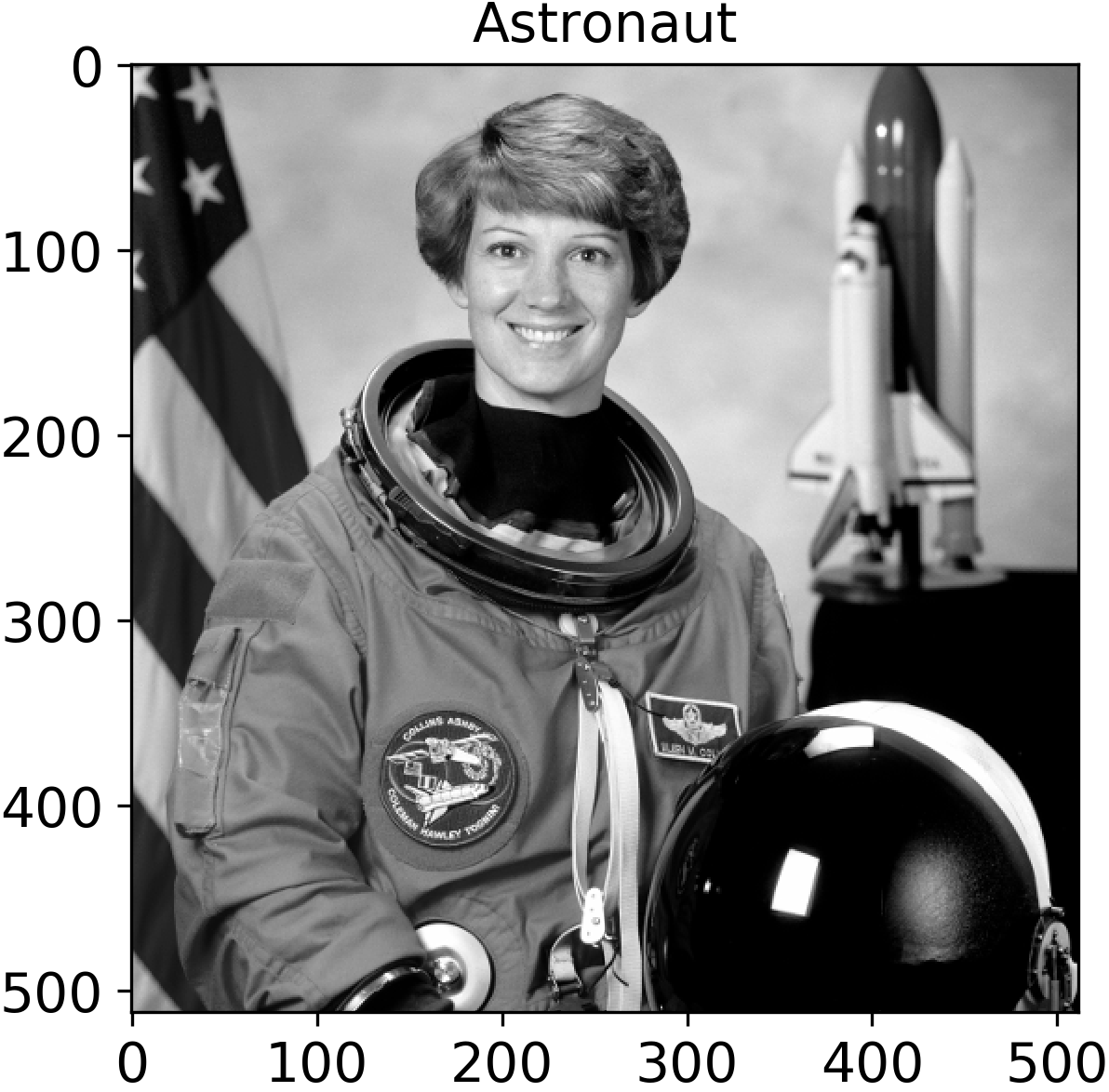}
	\caption{}
\end{subfigure}\hspace{0.05\textwidth}
\begin{subfigure}[b]{0.45\textwidth}
\centering
	\includegraphics[width=0.8\textwidth]{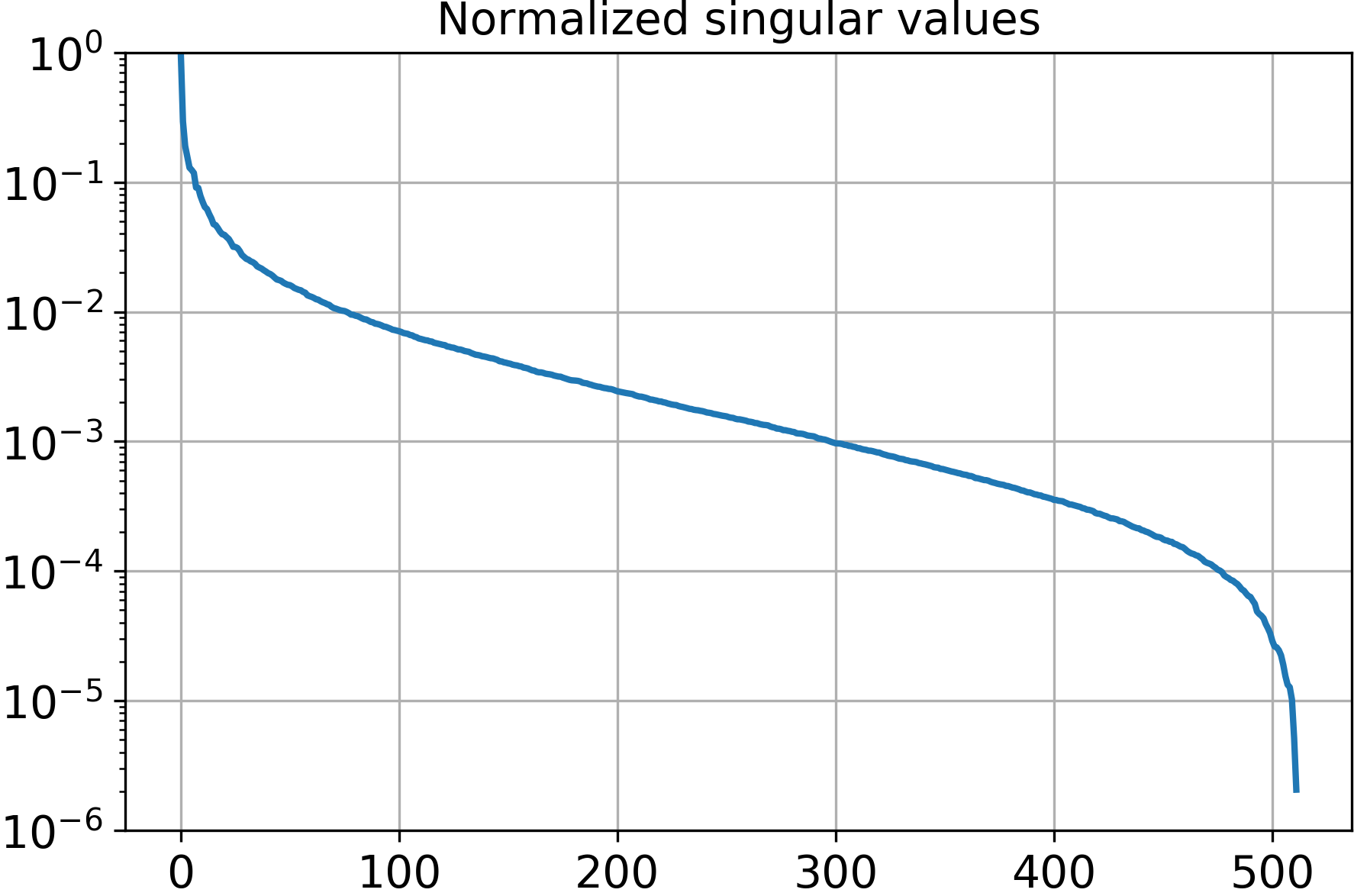}
	\caption{}
\end{subfigure}
\begin{subfigure}[b]{0.45\textwidth}
\centering
	\includegraphics[width=0.75\textwidth]{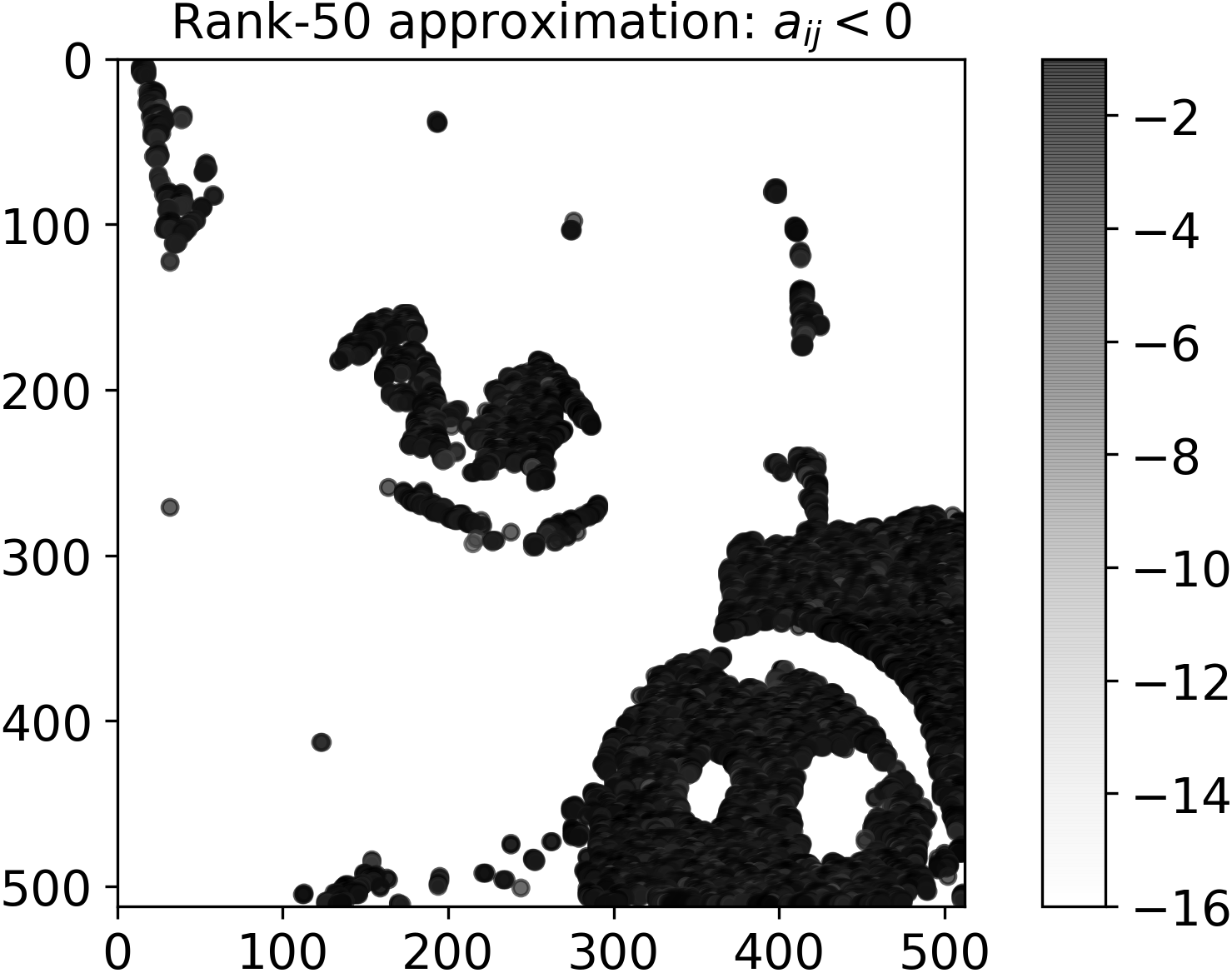}
	\caption{}
\end{subfigure}\hspace{0.05\textwidth}
\begin{subfigure}[b]{0.45\textwidth}
\centering
	\includegraphics[width=0.75\textwidth]{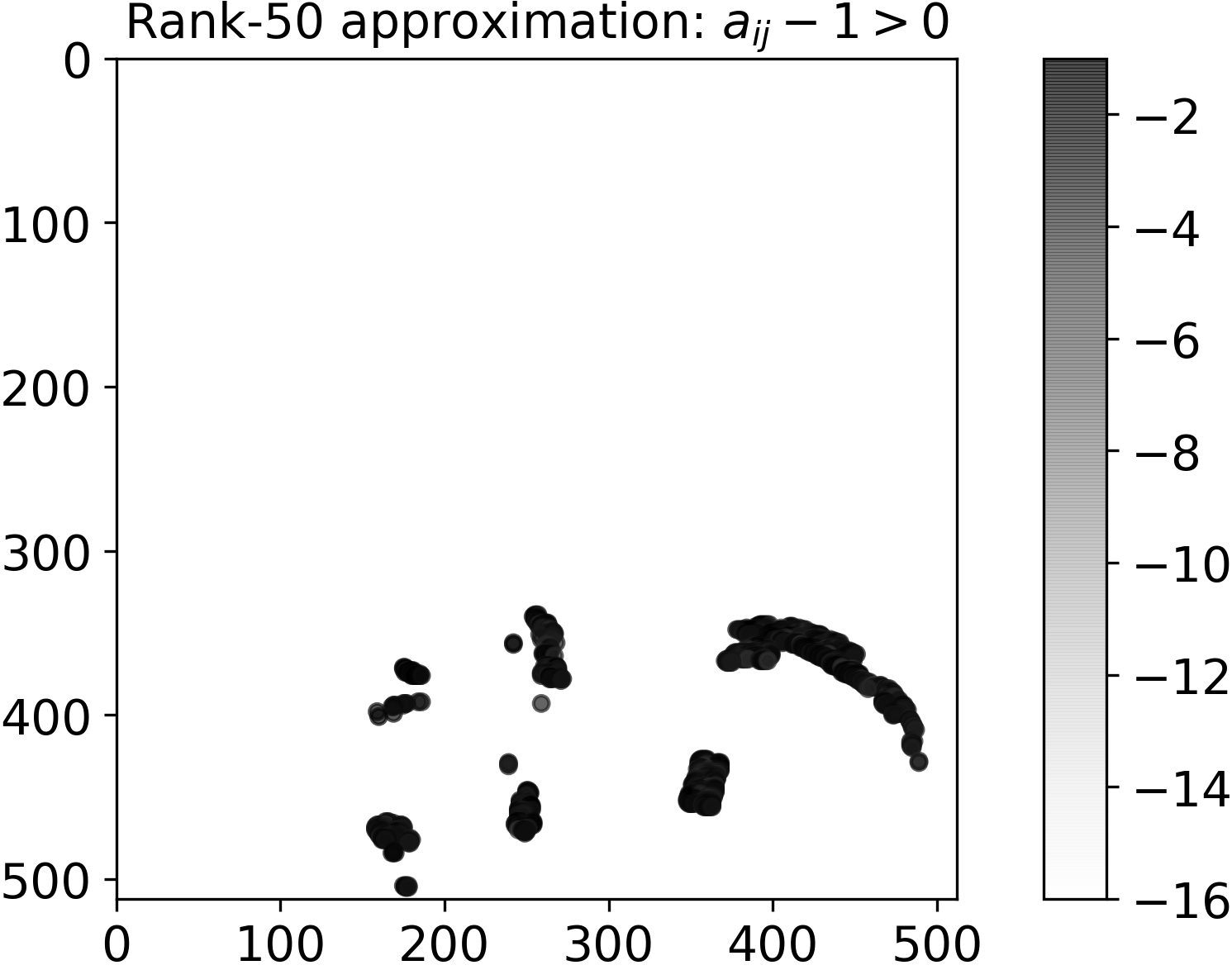}
	\caption{}
\end{subfigure}
\caption{Properties of the $512 \times 512$ Astronaut image: the image itself~(a), the normalized singular values~(b), the magnitude (in log scale) of the negative elements of its best rank-$50$ approximation~(c), and the elements greater than 1 (deviation from 1 in log scale) of its best rank-$50$ approximation~(d).}
\label{num:fig:astro_data}
\end{figure}

\begin{figure}[th]
\centering
\begin{subfigure}[b]{0.9\textwidth}
\centering
	\includegraphics[width=0.4\textwidth]{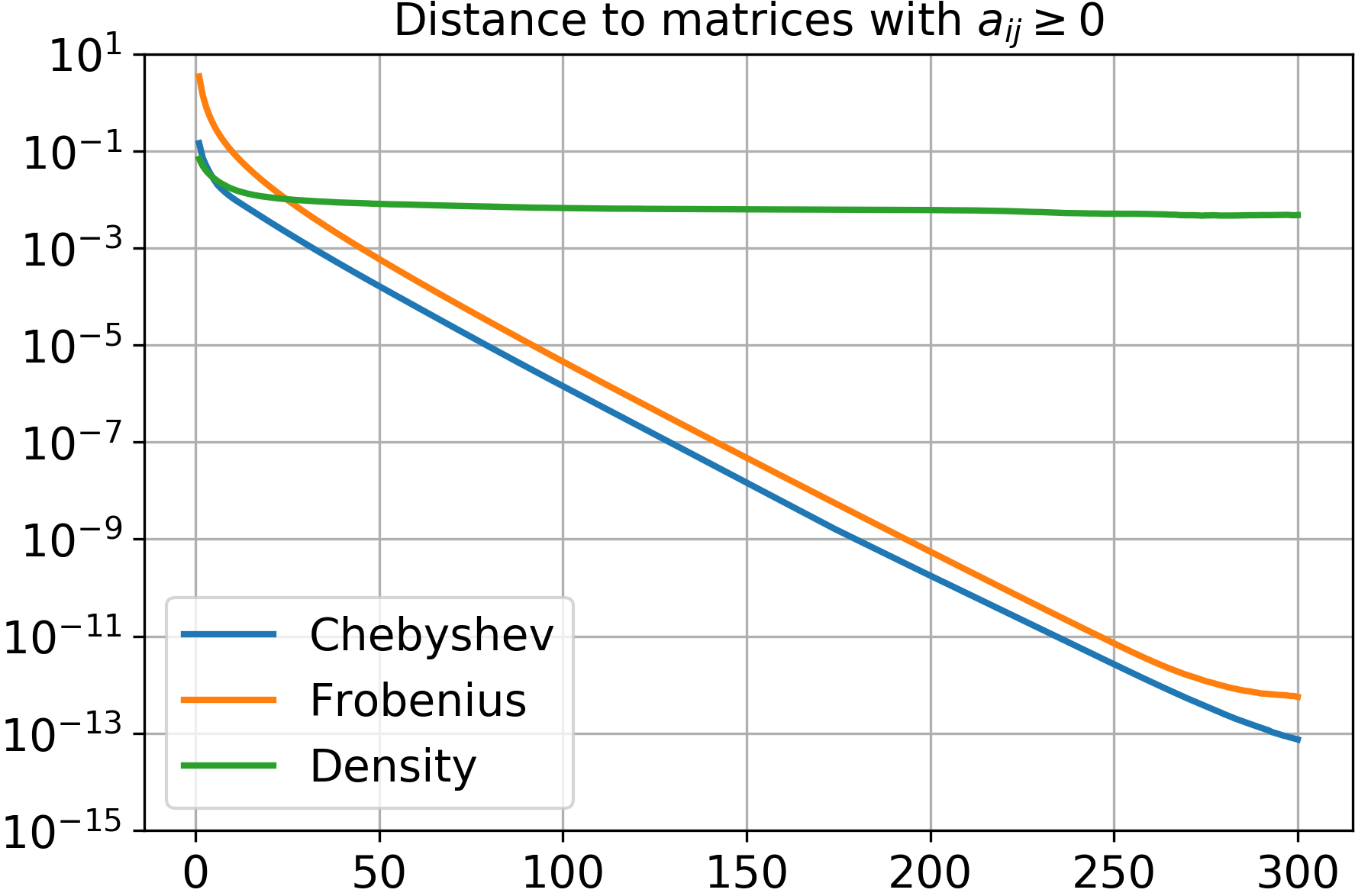}\hspace{0.05\textwidth}%
	\includegraphics[width=0.4\textwidth]{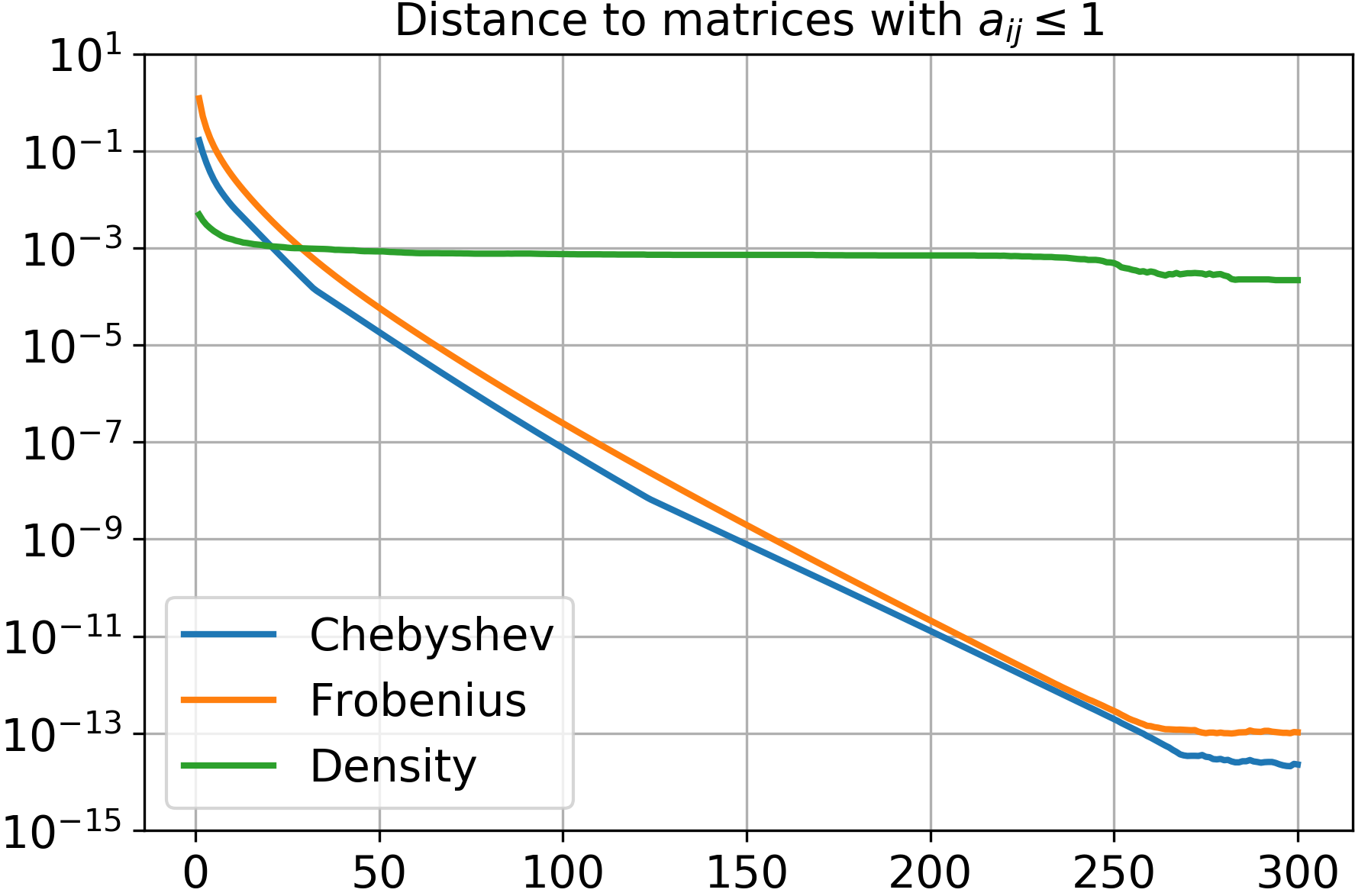}
	\caption{SVD}
\end{subfigure}
\begin{subfigure}[b]{0.9\textwidth}
\centering
	\includegraphics[width=0.4\textwidth]{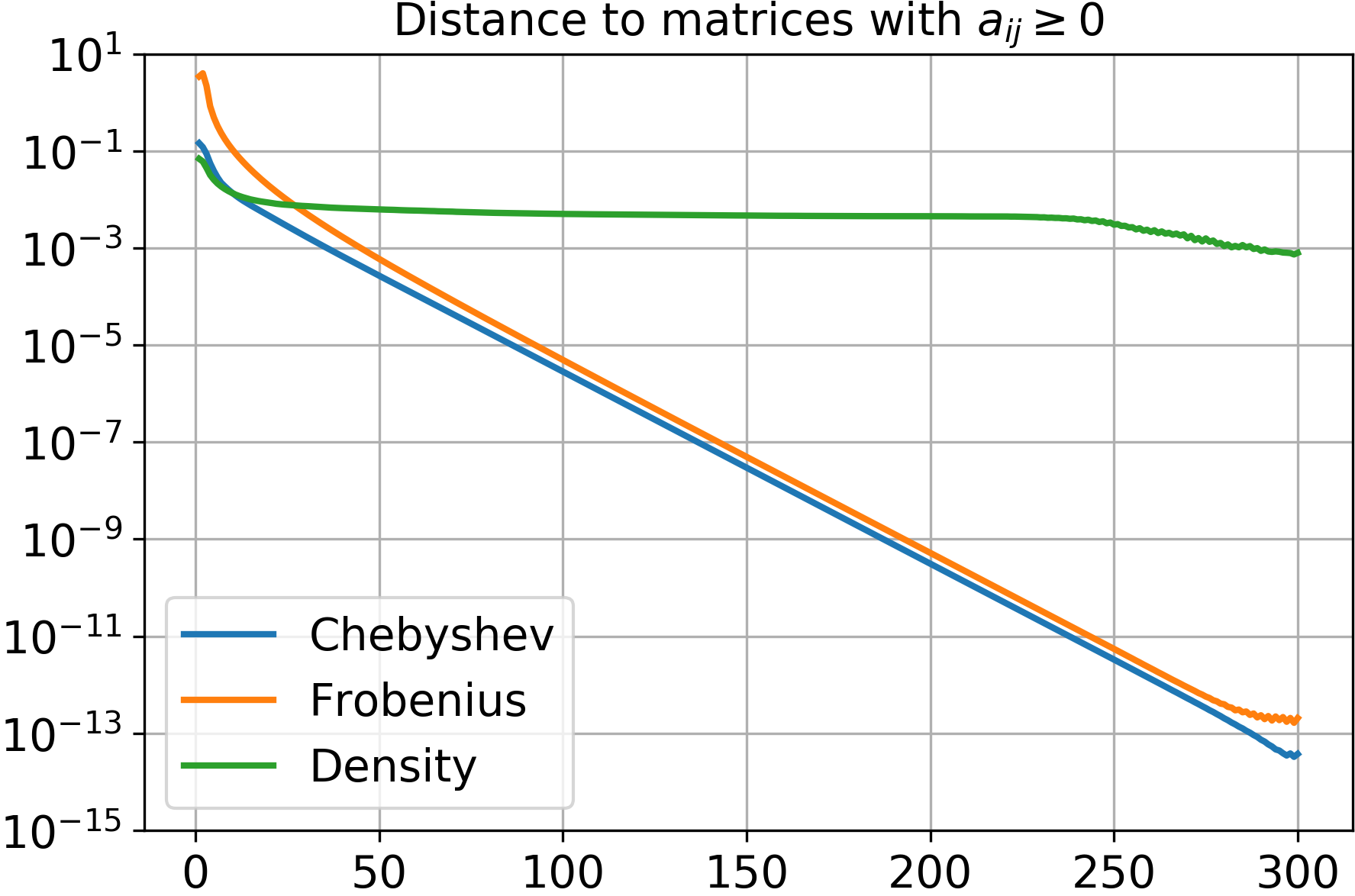}\hspace{0.05\textwidth}%
	\includegraphics[width=0.4\textwidth]{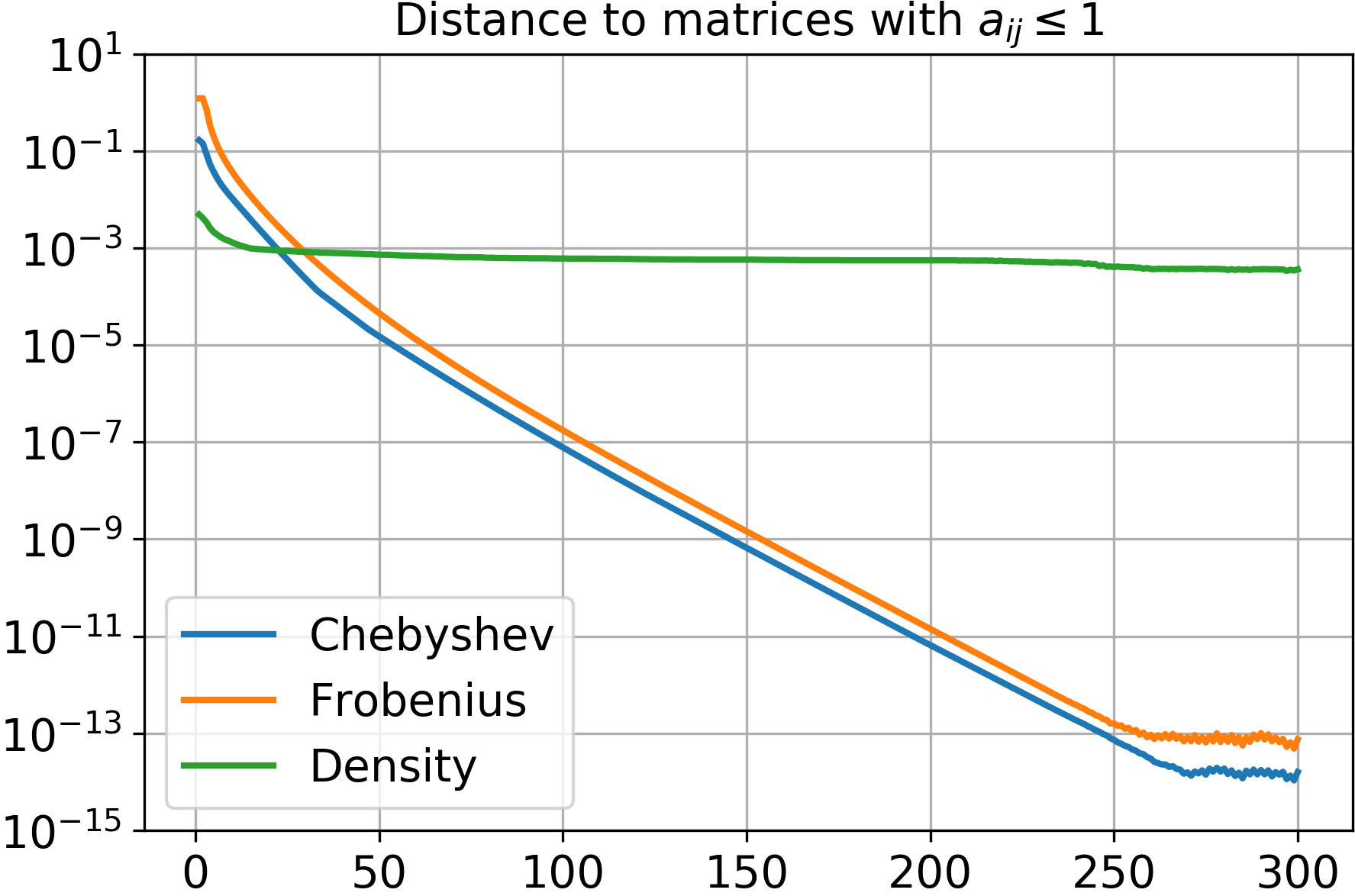}
	\caption{Tangent}
\end{subfigure}
\begin{subfigure}[b]{0.9\textwidth}
\centering
	\includegraphics[width=0.4\textwidth]{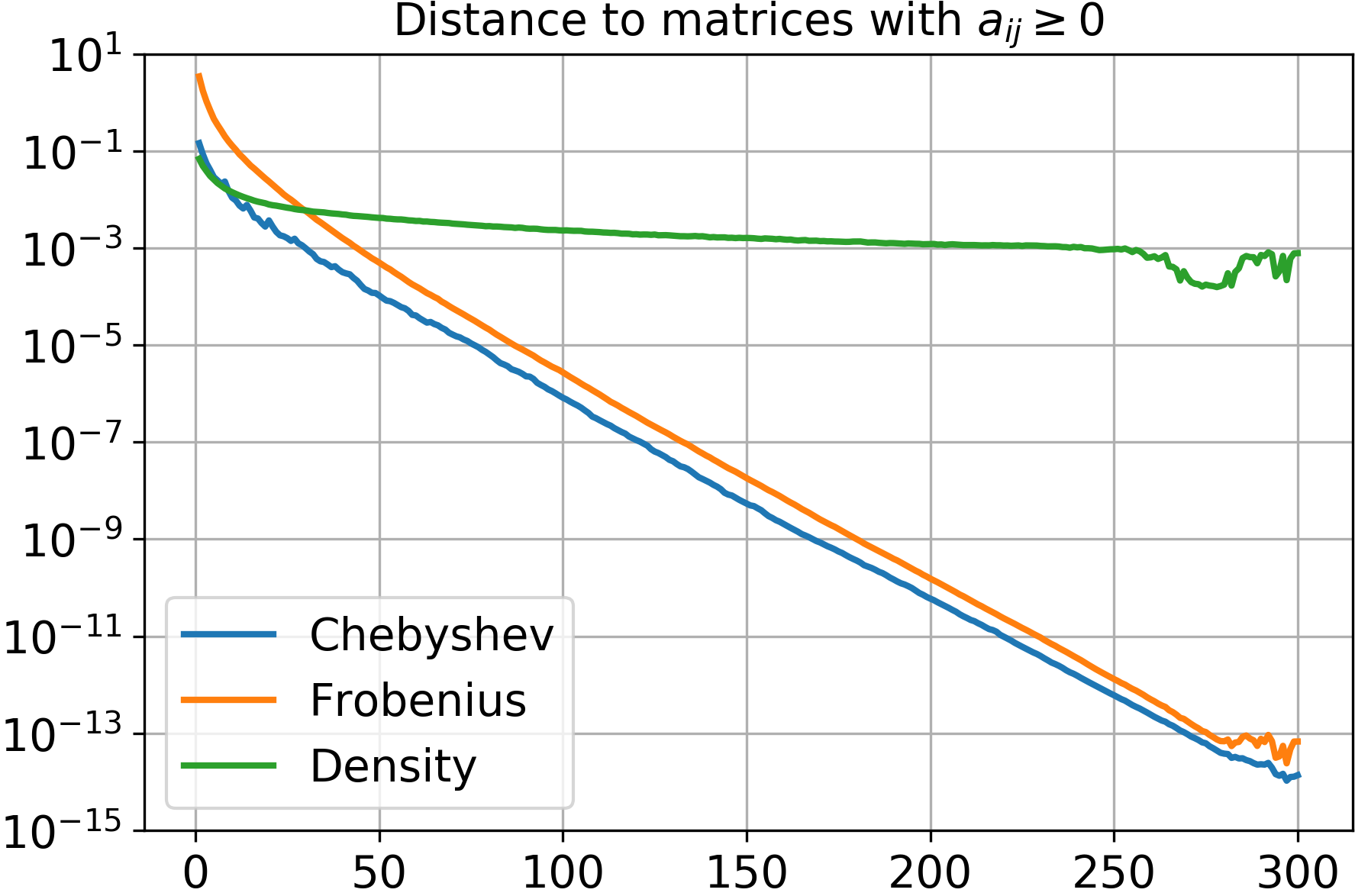}\hspace{0.05\textwidth}%
	\includegraphics[width=0.4\textwidth]{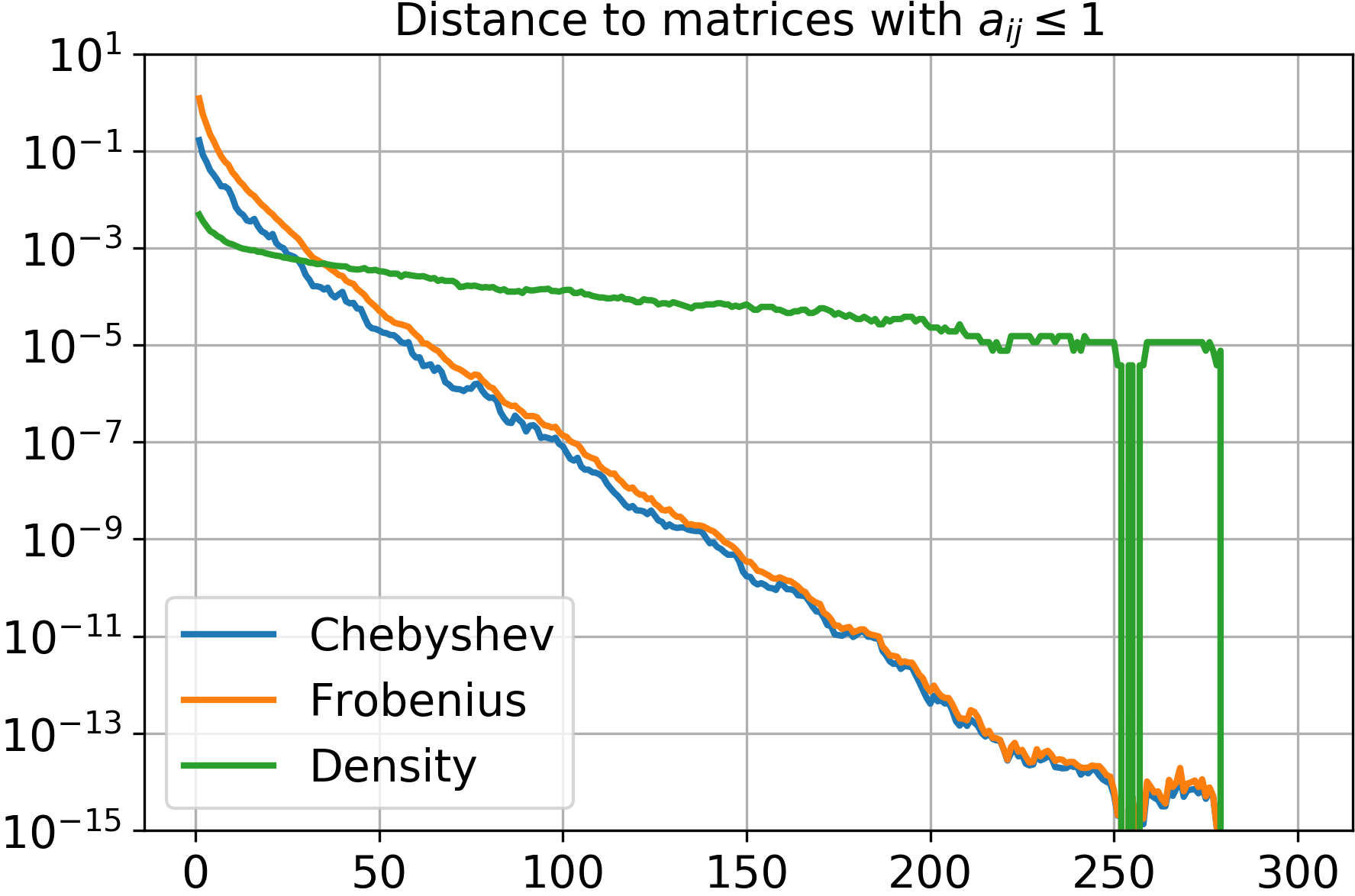}
	\caption{HMT(0, 60), $\mathrm{Rad}(0.2)$}
\end{subfigure}
\begin{subfigure}[b]{0.9\textwidth}
\centering
	\includegraphics[width=0.4\textwidth]{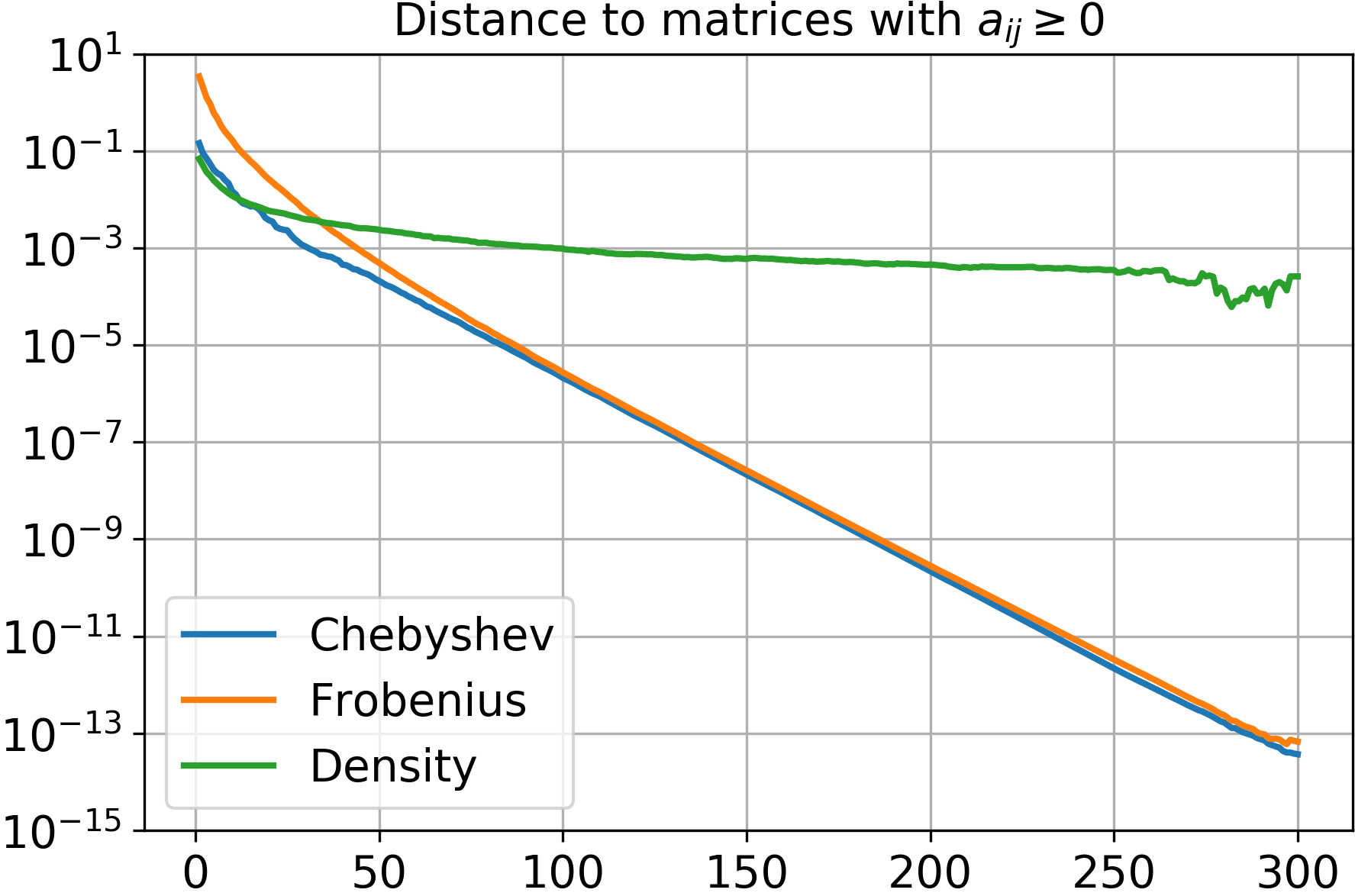}\hspace{0.05\textwidth}%
	\includegraphics[width=0.4\textwidth]{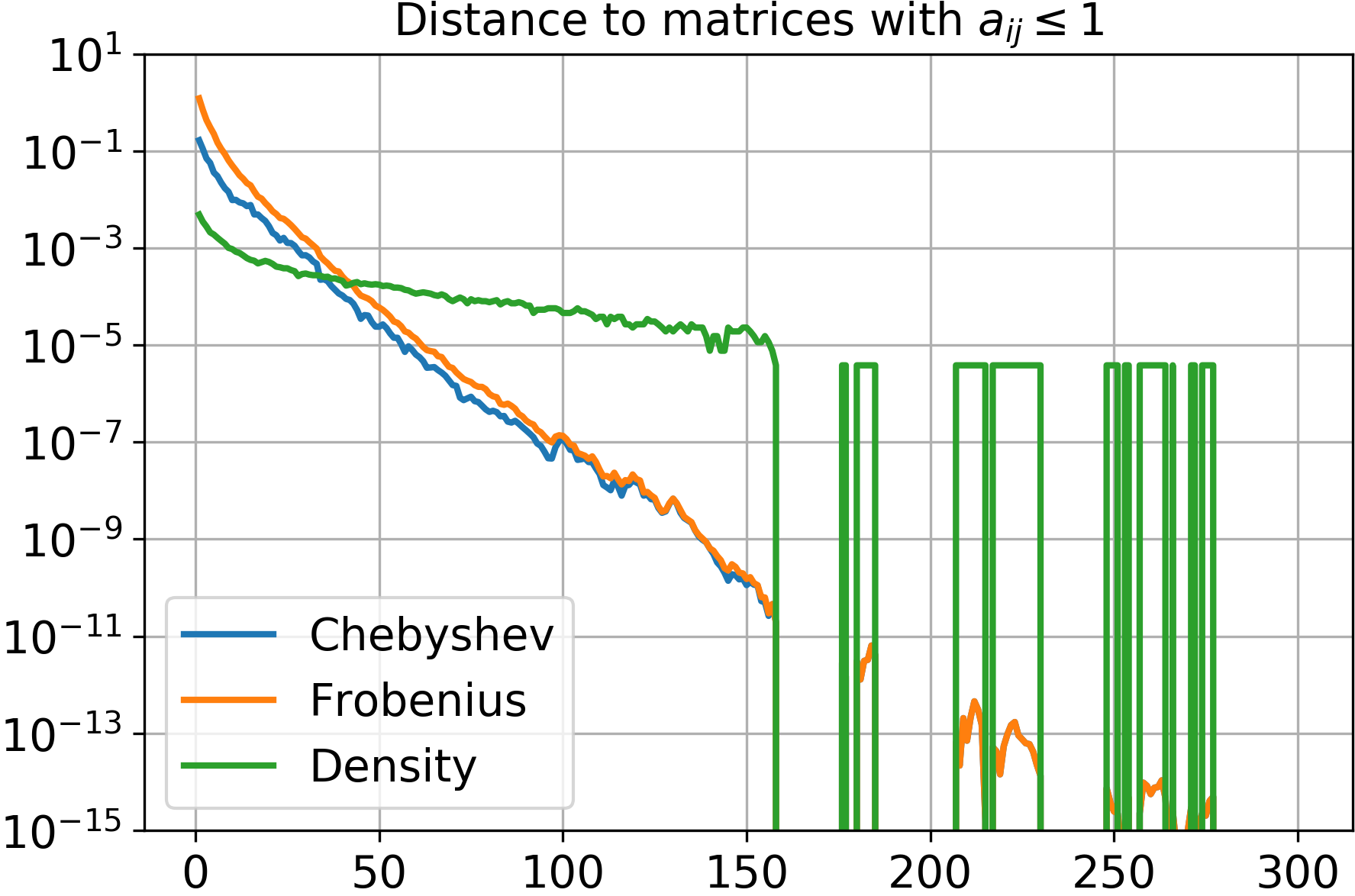}
	\caption{Tropp(60, 120), $\mathrm{Rad}(0.2)$}
\end{subfigure}
\begin{subfigure}[b]{0.9\textwidth}
\centering
	\includegraphics[width=0.4\textwidth]{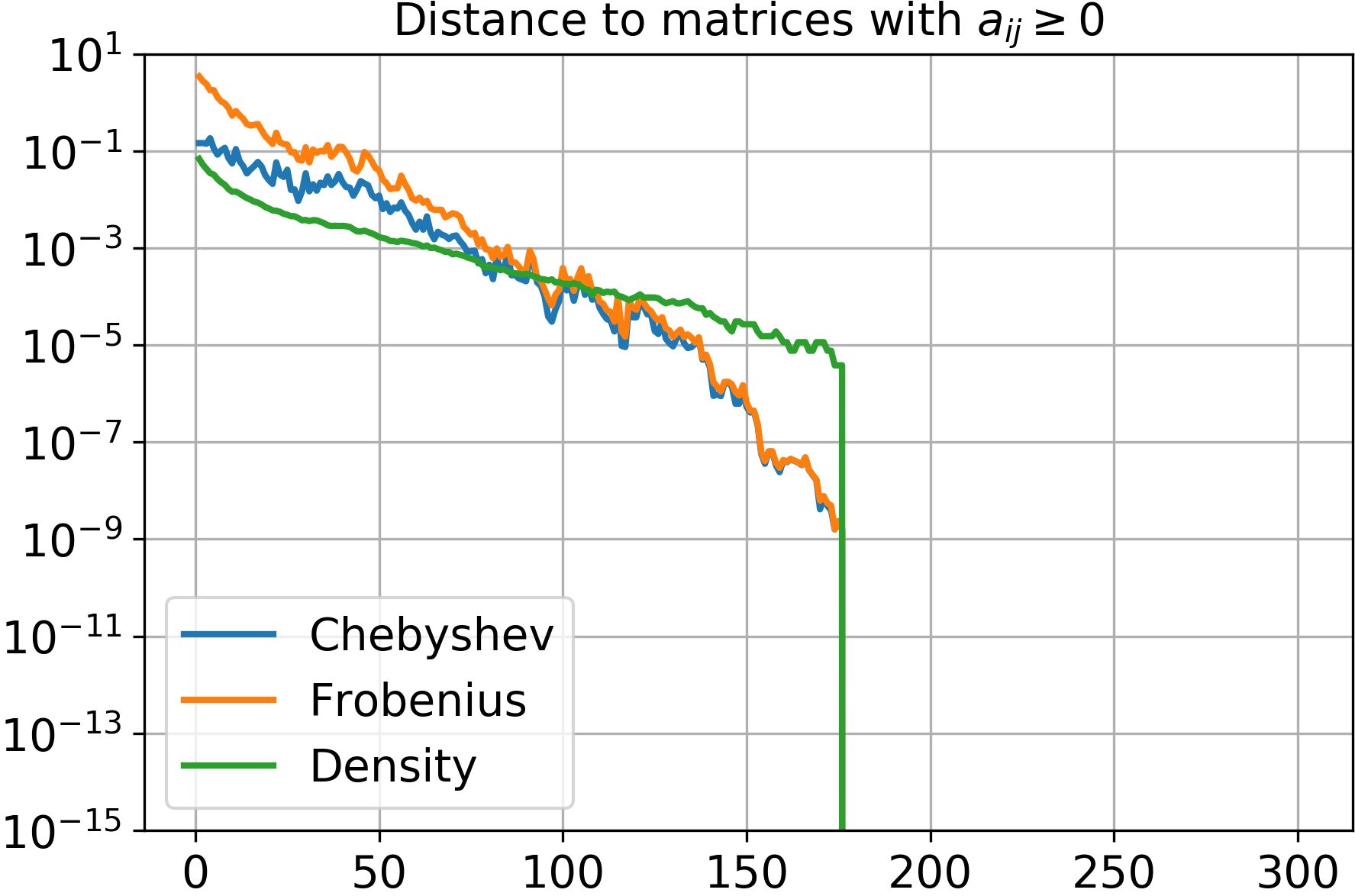}\hspace{0.05\textwidth}%
	\includegraphics[width=0.4\textwidth]{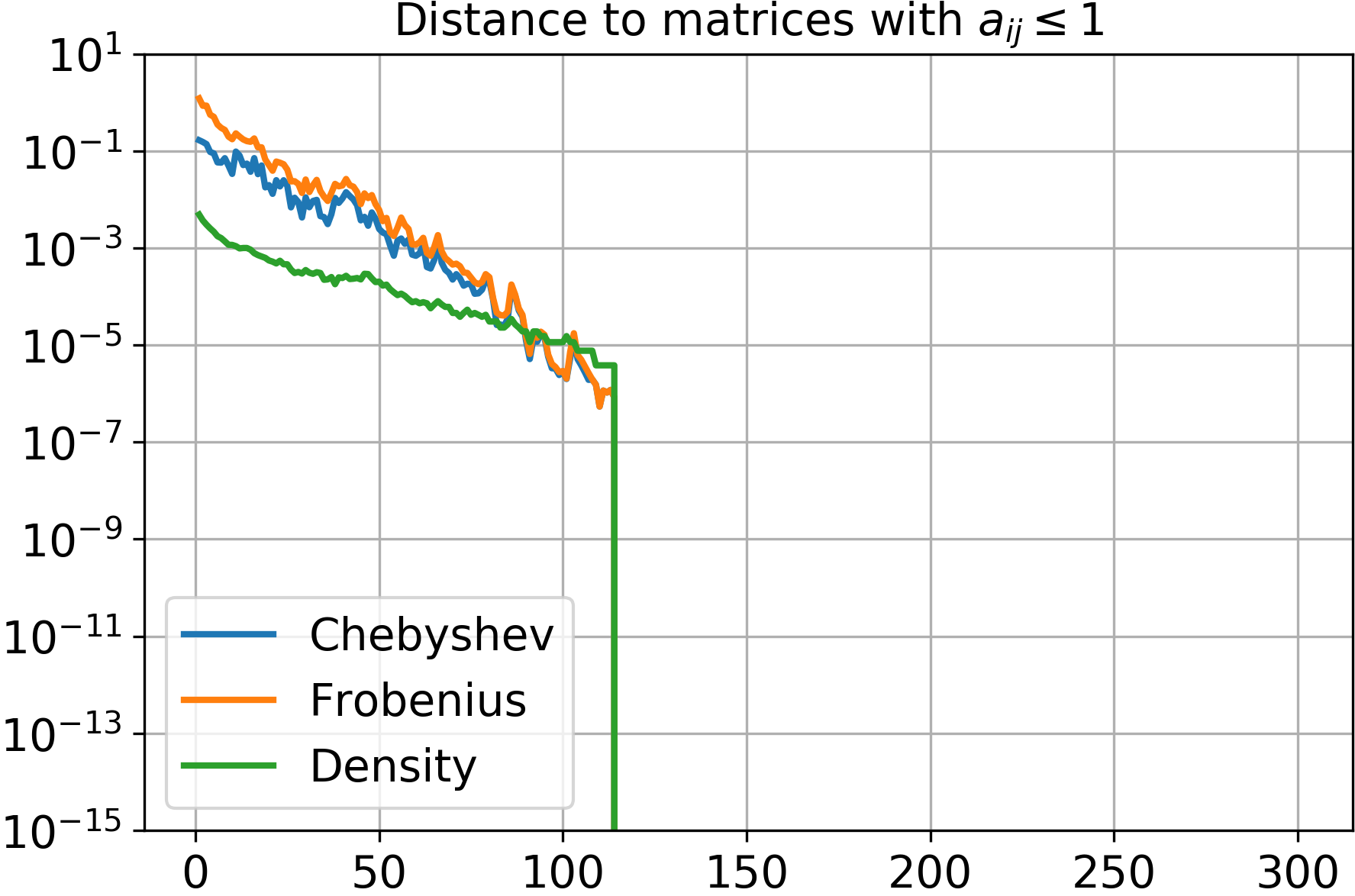}
	\caption{GN(340), $\mathrm{Rad}(0.2)$}
\end{subfigure}
\caption{Comparison of alternating projection methods for rank-$50$ nonnegative approximation of the Astronaut image: the Frobenius and Chebyshev norms of the negative part and the density of negative elements over 300 iterations~(left), same for the elements greater than 1~(right).}
\label{num:fig:astro_ap_comparison}
\end{figure}

\begin{figure}[th]
\begin{subfigure}[b]{0.3\textwidth}
\centering
	\includegraphics[width=0.9\textwidth]{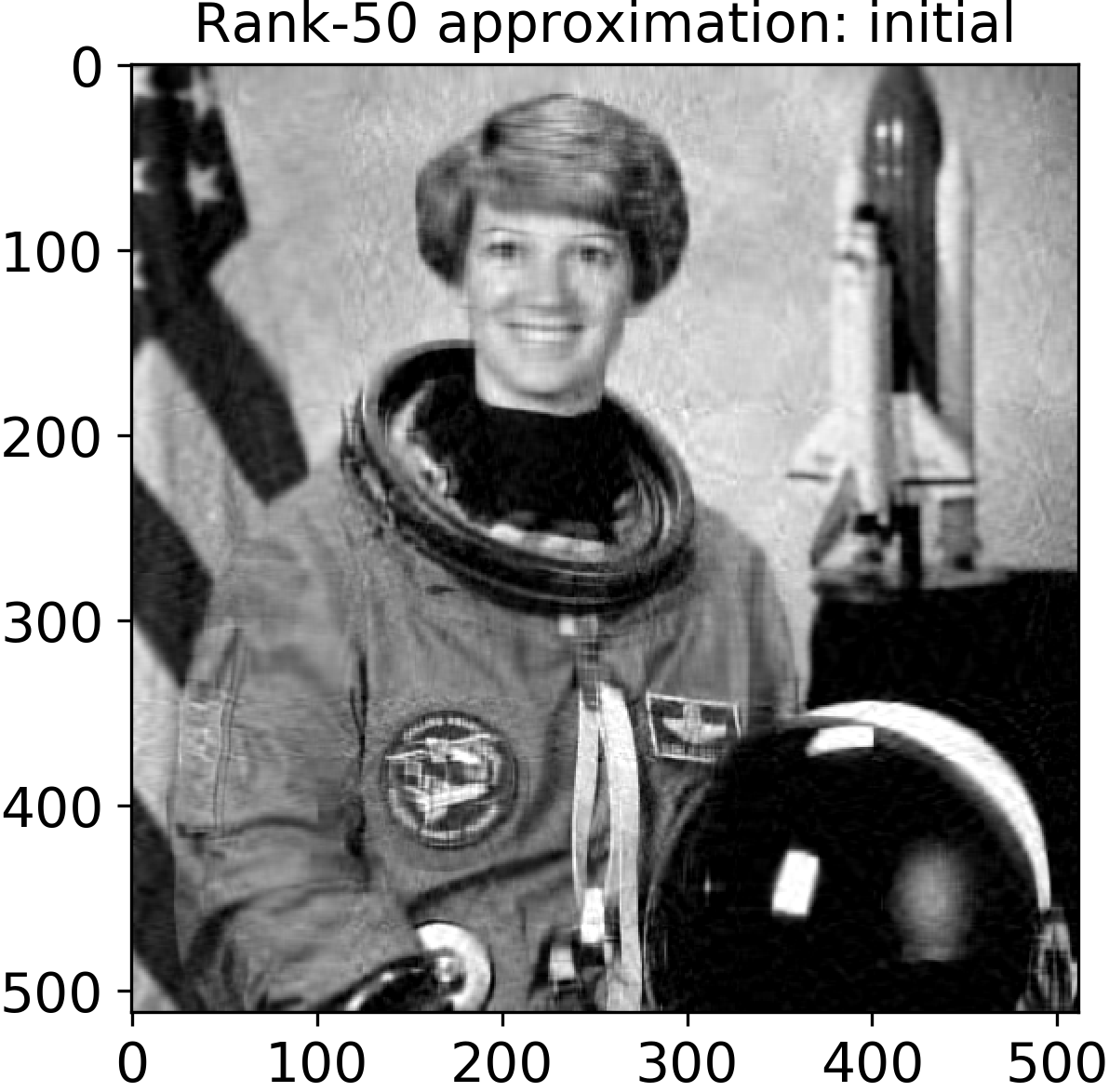}
	\caption{Initial SVD}
\end{subfigure}\hfill%
\begin{subfigure}[b]{0.3\textwidth}
\centering
	\includegraphics[width=0.9\textwidth]{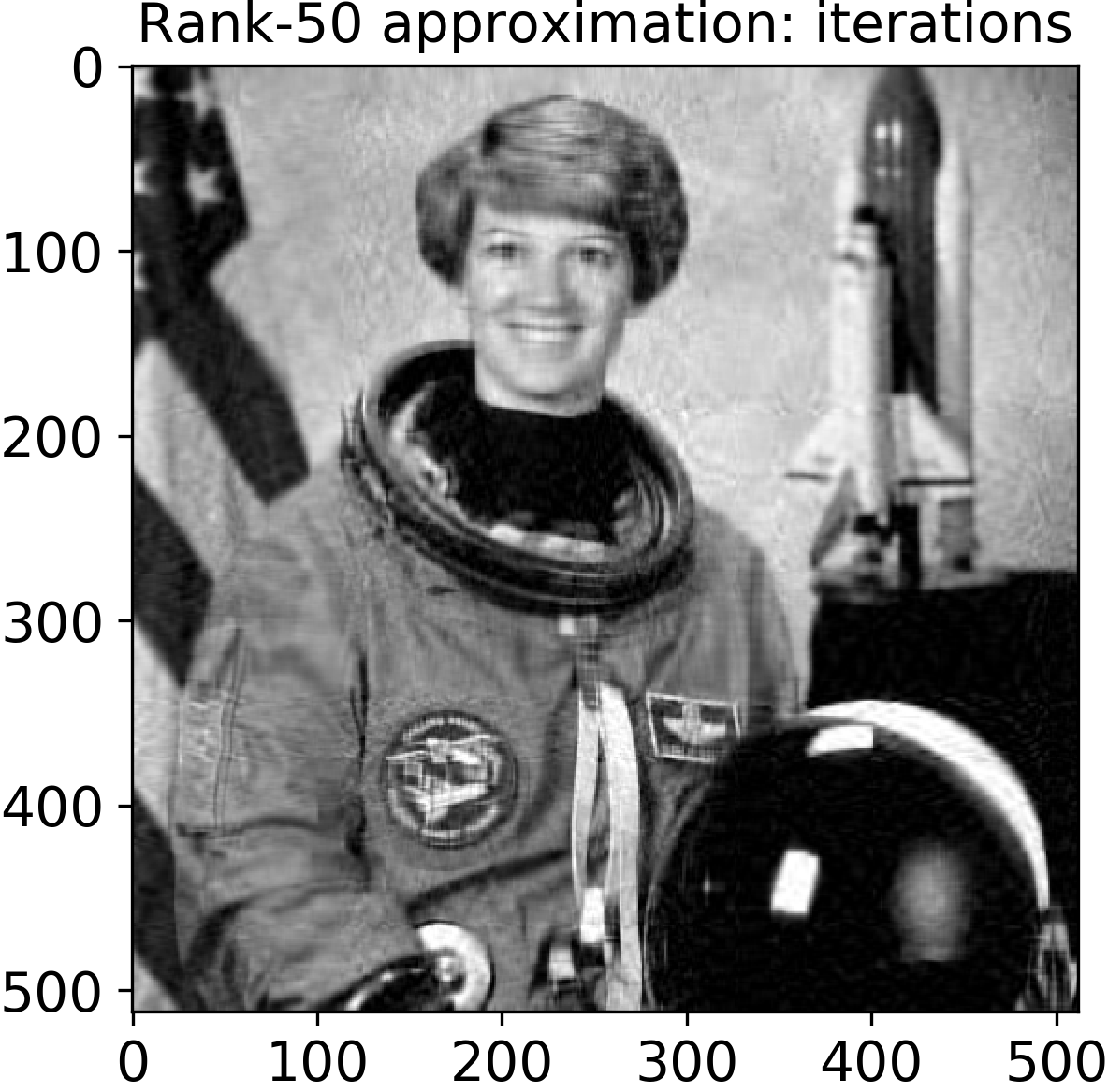}
	\caption{SVD}
\end{subfigure}\hfill%
\begin{subfigure}[b]{0.3\textwidth}
\centering
	\includegraphics[width=0.9\textwidth]{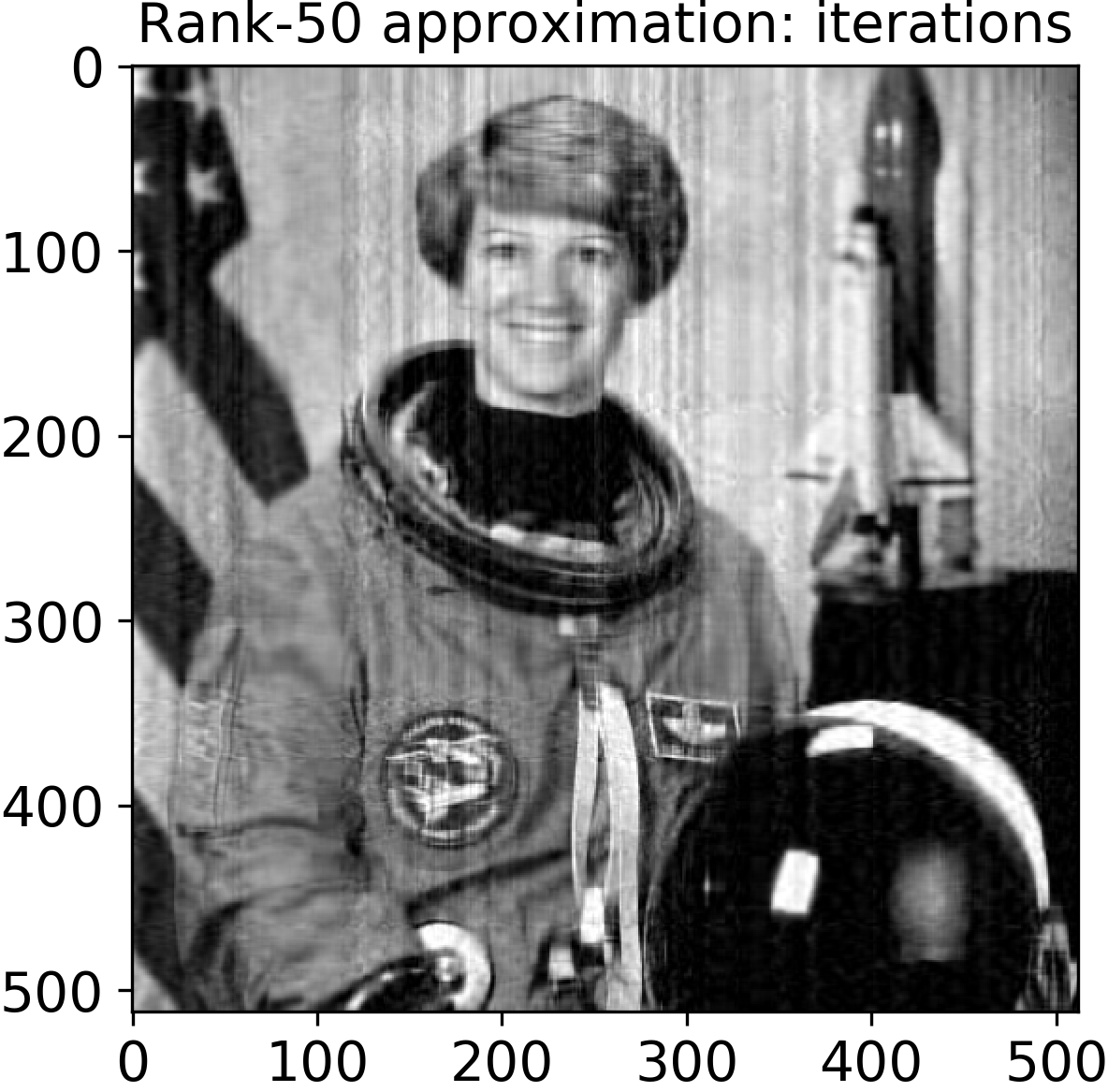}
	\caption{Tangent}
\end{subfigure}
\begin{subfigure}[b]{0.3\textwidth}
\centering
	\includegraphics[width=0.9\textwidth]{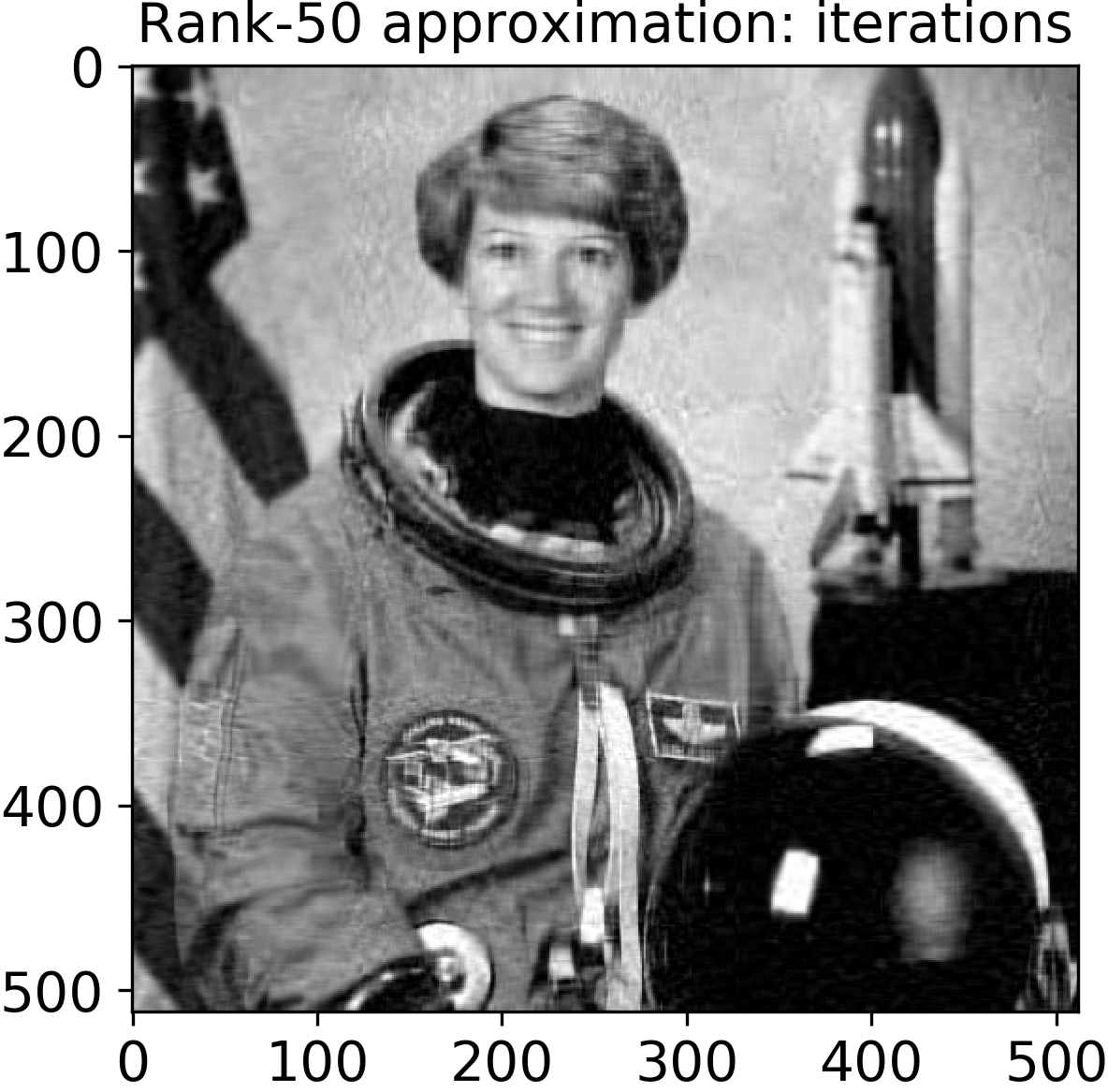}
	\caption{HMT(0, 60), $\mathrm{Rad}(0.2)$}
\end{subfigure}\hfill%
\begin{subfigure}[b]{0.3\textwidth}
\centering
	\includegraphics[width=0.9\textwidth]{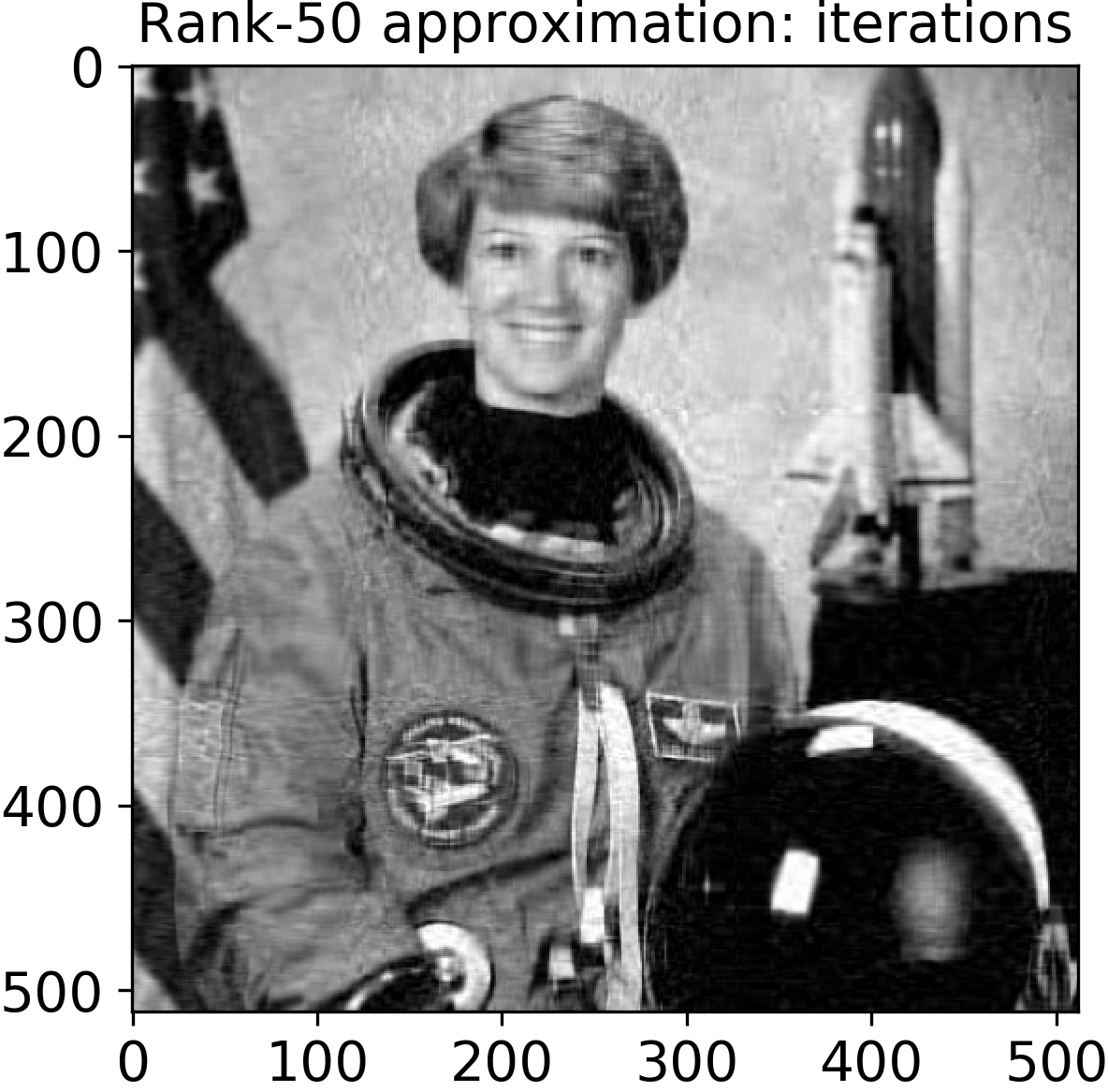}
	\caption{Tropp(60, 120), $\mathrm{Rad}(0.2)$}
\end{subfigure}\hfill%
\begin{subfigure}[b]{0.3\textwidth}
\centering
	\includegraphics[width=0.9\textwidth]{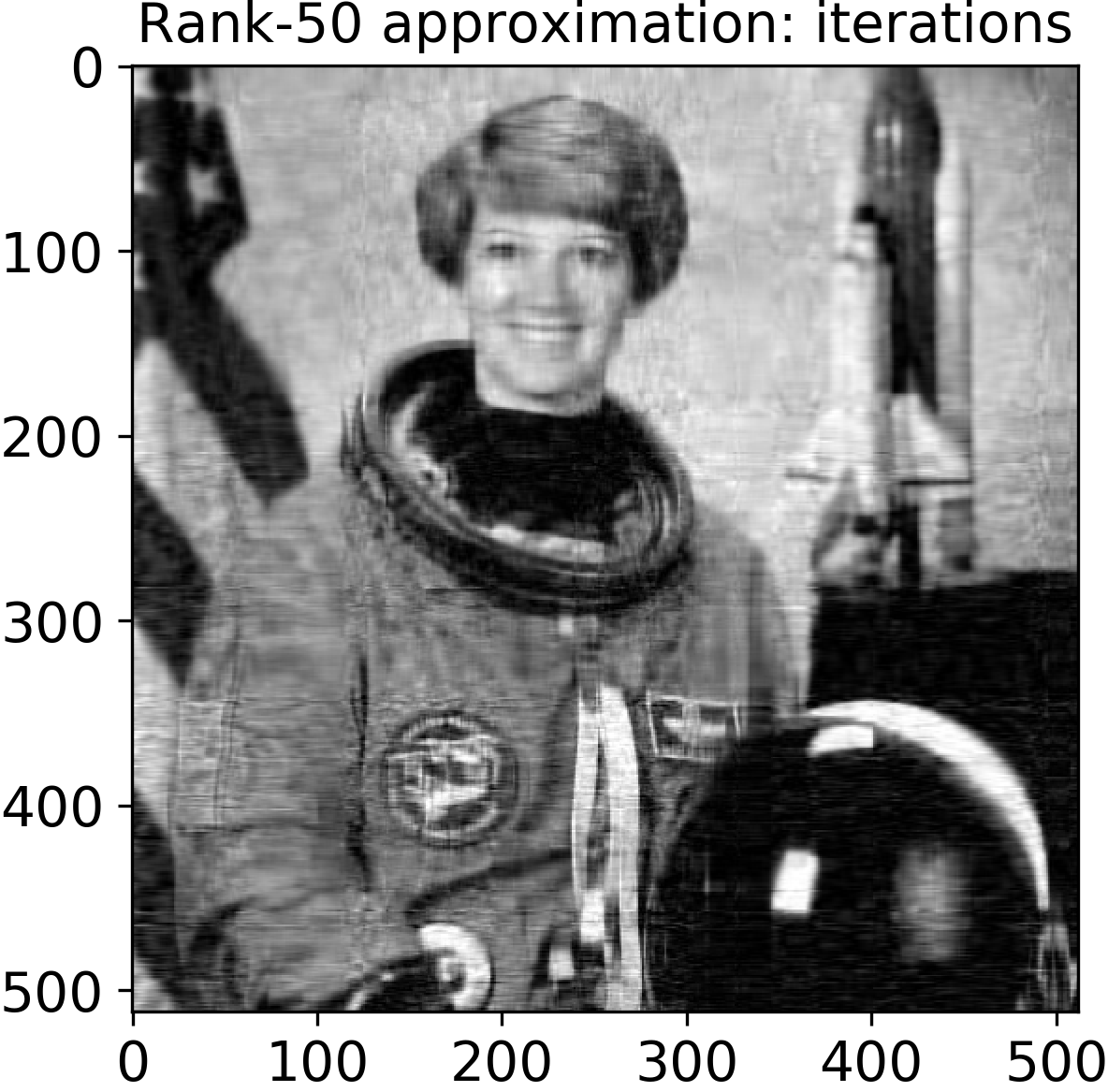}
	\caption{GN(340), $\mathrm{Rad}(0.2)$}
\end{subfigure}
\caption{Comparison of rank-$50$ nonnegative approximations of the Astronaut image.}
\label{num:fig:astro_ap_comparison_visual}
\end{figure}

The third example aims to show that alternating projections can be used to clip the values of a low-rank matrix to a prescribed range. We pick a $512 \times 512$ grayscale image and look for its rank-50 approximation, whose values lie in $[0, 1]$: this requires a simple modification of the algorithms. The best rank-50 approximation, that we refine with alternating projections, contains outliers both below $0$ and above $1$, as Fig.~\ref{num:fig:astro_data} shows. We see from Tab.~\ref{num:tab:astro} and Fig.~\ref{num:fig:astro_ap_comparison} that all methods converge and that randomized approaches are faster. By visually comparing the resulting approximations in Fig.~\ref{num:fig:astro_ap_comparison_visual}, we note that Tangent introduced vertical artifacts and GN lead to more disturbances than HMT and Tropp.

\begin{table}[h]\centering
\begin{tabular}{@{}llccccc@{}}
    \toprule
    Method & Sketch & Flops per iter & Frobenius & Chebyshev \\
    \midrule
    Initial $\mathrm{SVD}_r$ & N/A & $2.8 \cdot 10^9$ & $8.07 \cdot 10^{-2}$ & $4.94 \cdot 10^{-1}$ \\
    \midrule
    SVD & N/A & $2.8 \cdot 10^9$ & $8.30 \cdot 10^{-2}$ & $5.24 \cdot 10^{-1}$ \\
    Tangent & N/A & $1.3 \cdot 10^8$ & $1.04 \cdot 10^{-1}$ & $5.28 \cdot 10^{-1}$ \\
    HMT(0, 60) & $\mathrm{Rad}(0.2)$ & $8.7 \cdot 10^7$ & $8.50 \cdot 10^{-2}$ & $5.22 \cdot 10^{-1}$ \\
    HMT(0, 55) & $\mathrm{Rad}(0.2)$ & $8.0 \cdot 10^7$ & $8.82 \cdot 10^{-2}$ & $5.22 \cdot 10^{-1}$ \\
    Tropp(65, 110) & $\mathrm{Rad}(0.2)$ & $7.6 \cdot 10^7$ & $8.77 \cdot 10^{-2}$ & $5.49 \cdot 10^{-1}$ \\
    Tropp(60, 120) & $\mathrm{Rad}(0.2)$ & $7.1 \cdot 10^7$ & $8.92 \cdot 10^{-2}$ & $5.34 \cdot 10^{-1}$ \\
    GN(340) & $\mathrm{Rad}(0.2)$ & $7.1 \cdot 10^7$ & $1.16 \cdot 10^{-1}$ & $6.93 \cdot 10^{-1}$ \\
    GN(150) & $\mathrm{Rad}(0.2)$ & $4.8 \cdot 10^7$ & $1.31 \cdot 10^{-1}$ & $6.94 \cdot 10^{-1}$ \\
    \bottomrule\\
\end{tabular}
\caption{Comparison of alternating projection methods for rank-$50$ nonnegative approximation of the $512 \times 512$ Astronaut image: their computational complexities and relative errors in the Frobenius and Chebyshev norms after 300 iterations.}
\label{num:tab:astro}
\end{table}

\clearpage
\section{Conclusion}
In our work, we demonstrated that randomized sketching techniques can be successfully applied to solve the LRNMF problem via approximate alternating projections. By analyzing the computational complexities of randomized and deterministic approaches and evaluating them in three different numerical experiments, we showed that sketching can lead to more efficient algorithms than the tangent-space-based one and that they, nonetheless, exhibit similar convergence properties. We believe it is important to extend the LRNMF approaches to the multi-dimensional tensor case \cite{jiang2020nonnegative} and will attempt to do so in our future papers.

\begin{remark}
Since the first publication of this paper as a preprint, we have applied alternating projections with randomized sketching to compute low-rank nonnegative approximations of tensors in Tucker and tensor train formats. See \cite{sultonov2022low}.
\end{remark}

\section*{Acknowledgements}
We thank Dmitry Zheltkov, Nikolai Zamarashkin and Eugene Tyrtyshnikov for useful discussions. This work was supported by Russian Science Foundation (project 21-71-10072).

\bibliography{common/bib}
\bibliographystyle{ieeetr}

\appendix
\section{Computational complexity}
In this Appendix, we summarize the computational complexities of the different procedures encountered in the methods that we discussed and of the methods themselves. For reference, we use Golub's and Van Loan's classic \cite{GolubVanLoanMatrix2013}. The complexities will be listed in terms of flops, which are undrestood as in \cite[Sec.~1.1.15]{GolubVanLoanMatrix2013}.

\subsection{QR decomposition}
See \cite[Secs.~5.2.2 and 5.1.6]{GolubVanLoanMatrix2013}. The thin QR decomposition of a tall full-rank $m \times n$ matrix via Householder reflections requires
\begin{equation*}
    C_{QR}(m,n) = 4 mn^2 - \frac{4}{3} n^3
\end{equation*}
flops if the matrix $Q \in \Real^{m \times n}$ is explicitly formed from a product of $n$ reflections. If $Q$ is not formed, the cost of QR decomposition is halved and the subsequent matrix-vector products become cheaper.

\subsection{Singular value decomposition}
See \cite[Sec.~8.6.3]{GolubVanLoanMatrix2013}. The estimated number of flops needed to compute the economic SVD of a tall $m \times n$ matrix is
\begin{equation*}
    C_{SVD}(m,n) = 14 mn^2 + 8 n^3.
\end{equation*}
If $m \gg n$, a more efficient approach can be used to get
\begin{equation*}
    C_{SVD}(m,n) = 6 mn^2 + 20 n^3,
\end{equation*}
and if the matrix is square, $m = n$, its SVD can be computed in approximately
\begin{equation*}
    C_{SVD}(m,m) = 21 m^3.
\end{equation*}

\subsection{Sketches}
We mentioned three types of random sketching matrices: standard Gaussian, Rademacher, and sparse Rademacher. To generate any of these, one requires a pseudo-random number generator that produces samples from the uniform distribution on $[0, 1]$. A widespread choice is the Mersenne Twister \cite{MatsumotoNishimuraMersenne1998}, which does one floating-point division per sample.

To sample from the standard Gaussian distribution, one can apply the Box-Muller transform \cite{BoxMullerNote1958}. It takes two independent random variables $x_1$ and $x_2$ uniformly distributed on $[0, 1]$ and evaluates
\begin{equation*}
    y_1 = \sqrt{-2 \log(x_1)} \cos(2 \pi x_2), \quad y_2 = \sqrt{-2 \log(x_1)} \sin(2 \pi x_2)
\end{equation*}
that are independent standard Gaussian. The transform requires the computation of $\sqrt{\cdot}$, $\log$, $\sin$, and $\cos$, which are computationally demanding functions. The equivalent number of flops that they consume depends heavily on the particular processing unit, but for our needs we will crudely bound their contribution by assigning 150 flops to one application of the Box-Muller transform (two Mersenne-Twister samples included); hence, 75 flops per sample from the standard Gaussian distribution. Another popular class of algorithms are rejection sampling ones, exemplified by the Ziggurat algorithm \cite{MarsagliaTsangZiggurat2000}.

Sampling from the Rademacher distribution requires a single sample from the uniform distribution on $[0, 1]$. Generating an $m \times n$ matrix with random Rademacher entries then takes $mn$ flops, and it can be applied to a vector in $mn$ flops as well (twice faster than a standard matrix-vector product since no floating-point multiplications are done).

The construction of a sparse Rademacher sketching matrix is done in two steps: we generate a sparse mask and then sample the corresponding elements. Let $0 < \rho \leq 1$ be the desired density, i.e. the number of entries in the mask divided by the size of the matrix $mn$. To decide if an entry should be added to the mask, we need one sample from the uniform distribution on $[0,1]$. Then for the chosen $\rho m n$ entries, we generate Rademacher samples using $\rho m n$ flops. This gives $(1 + \rho)mn$ flops for the whole sparse Radmeacher sketching matrix. The multiplication with it is also fast: it takes only $\rho m n$ flops for a matrix-vector product.

\subsection{Low-rank nonnegative matrix approximation}
Here, for all the methods discussed in the paper we list their detailed computational complexities for one iteration: it starts with a nonnegative matrix $A_i$ (or a low-rank factorization of a matrix $B_i$) and results in an updated nonnegative matrix $A_{i + 1}$ (or in an updated low-rank factorization of $B_{i+1}$, respectively). We always assume that the matrix is tall, $m \geq n$, and that the target rank is $r$.
\begin{enumerate}
    \item Alternating projections \cite{SongNgNonnegative2020}:
    \begin{equation*}
        C_{SVD}(m,n) + 2mnr + nr
    \end{equation*}
    \item Tangent-space-based alternating projections \cite{SongEtAlTangent2020}:
    \begin{equation*}
        6mnr + 10mr^2 + 12nr^2 + nr + \left(165 + \frac{1}{3}\right) r^3
    \end{equation*}
    \item HMT \cite[Algs.~4.4 and 5.1]{HalkoEtAlFinding2011a} with co-range sketch size $k \geq r$ and $p$ iterations of the power method:
    \begin{enumerate}
        \item Gaussian sketching
        \begin{equation*}
            (4p+4)mnk + 2mnr + (4p+6)(m+n)k^2 + nr + 75nk + \frac{8}{3}(7 - p)k^3
        \end{equation*}
        \item Rademacher sketching
        \begin{equation*}
            (4p+3)mnk + 2mnr + (4p+6)(m+n)k^2 + nr + nk + \frac{8}{3}(7 - p)k^3
        \end{equation*}
        \item Sparse Rademacher sketching with density $\rho$
        \begin{equation*}
            (4p+2+\rho)mnk + 2mnr + (4p+6)(m+n)k^2 + nr + (1+\rho)nk + \frac{8}{3}(7 - p)k^3
        \end{equation*}
    \end{enumerate}
    \item Tropp \cite[Alg.~4]{TroppEtAlPractical2017} with co-range sketch size $k \geq r$ and range sketch size $l \geq k$:
    \begin{enumerate}
        \item Gaussian sketching
        \begin{equation*}
            2(r + k + l)mn + 2mkl + 5mk^2 + 7nk^2 + 2nkl + 75(nk + ml) + nr + \left( 17 + \frac{1}{3} \right) k^3 + 4lk^2
        \end{equation*}
        \item Rademacher sketching
        \begin{equation*}
            (2r + k + l)mn + mkl + 5mk^2 + 7nk^2 + 2nkl + (nk + ml) + nr + \left( 17 + \frac{1}{3} \right) k^3 + 4lk^2
        \end{equation*}
        \item Sparse Rademacher sketching with density $\rho$
        \begin{equation*}
            (2r + \rho k + \rho l)mn + \rho mkl + 5mk^2 + 7nk^2 + 2nkl + (1 + \rho)(nk + ml) + nr + \left( 17 + \frac{1}{3} \right) k^3 + 4lk^2
        \end{equation*}
    \end{enumerate}
    \item GN \cite[Alg.~2.1]{NakatsukasaFast2020} with range sketch size $l \geq r$:
    \begin{enumerate}
        \item Gaussian sketching
        \begin{equation*}
            (4r + 2l)mn + 4 nlr + mr^2 + 75(nr + ml) + 4 lr^2 - \frac{4}{3} r^3 
        \end{equation*}
        \item Rademacher sketching
        \begin{equation*}
            (3r + l)mn + 3 nlr + mr^2 + (nr + ml) + 4 lr^2 - \frac{4}{3} r^3 
        \end{equation*}
        \item Sparse Rademacher sketching with density $\rho$
        \begin{equation*}
            (2r + \rho r + \rho l)mn + (2+\rho) nlr + mr^2 + (1+\rho)(nr + ml) + 4 lr^2 - \frac{4}{3} r^3 
        \end{equation*}
    \end{enumerate}
\end{enumerate}

\newpage

\end{document}